\documentclass[pre,showpacs,amsmath,amssymb,floatfix,onecolumn]{revtex4-1}
\pdfoutput=1 
\usepackage{graphicx}
\usepackage{epsfig}
\usepackage{dcolumn}
\usepackage{bm}
\usepackage{color}
\usepackage{amsthm}
\usepackage{amsfonts}
\usepackage{flushend}
\usepackage{clrscode}
\usepackage{hyperref}

\usepackage[toc,page]{appendix}

\def\(({\left(}
\def\)){\right)}
\def\[[{\left[}
\def\]]{\right]}

\newcommand{\bx}{{\textbf {x}}}

\newcommand{\be}{\begin{equation}}
\newcommand{\ee}{\end{equation}}
\newcommand{\bea}{\begin{eqnarray}}
\newcommand{\eea}{\end{eqnarray}}

\newcommand{\sign}{\text{ sign}}

\newcommand{\erfc}{{\rm erfc}}

\begin{document}

\title{Constrained Low-rank Matrix Estimation: \\ Phase Transitions, Approximate
  Message Passing and Applications.} \author{Thibault Lesieur$^{1}$, Florent
  Krzakala$^{2}$, and Lenka
  Zdeborov\'a$^{1,*}$}

\affiliation{
$^1$ Institut de Physique Th\'eorique, CNRS, CEA, Universit\'e Paris-Saclay, F-91191, Gif-sur-Yvette, France.\\
$^2$ Laboratoire de Physique Statistique, Ecole Normale Sup\'erieure,  \& Universit\'e Pierre et Marie Curie, Sorbonne Universit\'es.\\ 
$^*$ To whom correspondence shall be sent: lenka.zdeborova@cea.fr\\
}

\begin{abstract}
  This article is an extended version of previous work of the authors \cite{lesieur2015phase,lesieur2015mmse}
  on low-rank matrix estimation
  in the presence of
  constraints on the factors into which the matrix is
  factorized. Low-rank matrix factorization is one of the basic
  methods used in data analysis for unsupervised learning of relevant
  features and other types of dimensionality reduction. We present a
  framework to study the constrained low-rank matrix estimation for a
  general prior on the factors, and a general output channel through
  which the matrix is observed. We draw
  a paralel with the study of vector-spin glass models -- presenting a
  unifying way to study a number of problems considered previously in
  separate statistical physics works. We present a number of
  applications for the problem in data analysis. We derive in detail a
  general form of the low-rank approximate message passing (Low-RAMP)
  algorithm, that is known in statistical physics as the TAP
  equations. We thus unify the derivation of the TAP equations for
  models as different as the Sherrington-Kirkpatrick model, the
  restricted Boltzmann machine, the Hopfield model or vector (xy,
  Heisenberg and other) spin glasses. The state evolution of the
  Low-RAMP algorithm is also derived, and is equivalent to the replica
  symmetric solution for the large class of vector-spin glass
  models. In the section devoted to result we study in detail phase
  diagrams and phase transitions for the Bayes-optimal inference in
  low-rank matrix estimation. We present a typology of phase
  transitions and their relation to performance of algorithms such as
  the Low-RAMP or commonly used spectral methods.
\end{abstract}

\date{\today}
\maketitle

\tableofcontents

\newpage
\section{Introduction}

\subsection{Problem Setting}
In this paper we study a generic class of statistical physics models
having  Boltzmann probability measure that can be written in one of
the two following forms:
\begin{itemize}

\item{{\bf Symmetric vector-spin glass model}:
\begin{equation}
P(X|Y) = \frac{1}{Z_{X}(Y)}\prod\limits_{1 \leq i \leq N} P_{X}(x_i)
\prod\limits_{1 \leq i < j \leq N} e^{g(Y_{ij}, x_i^{\top} x_j/\sqrt{N})} \,.
\label{Boltz_XX}
\end{equation}
Here $Y_{ij}\in {\mathbb R}^{N \times N}$ and $X \in {\mathbb R}^{N
  \times r}$ are real valued
matrices. In this case $Y_{ij}$ is a symmetric
matrix. In statistical physics $Y$ is called the quenched disorder. In the whole paper we denote by $x_i \in {\mathbb R}^{r}$
the {\it vector-spin} $i$ ($r$-dimensional column vector) that collects the elements of the $i$th row of the
matrix $X$, $x_i^{\top} x_j$ is the scalar product of the two vectors. $Z_X(Y)$ is the corresponding partition function playing
role of the normalization. 
}
\item{{\bf Bipartite vector-spin glass model}:
\begin{equation}
P(U,V|Y) = \frac{1}{Z_{UV}(Y)}\prod\limits_{1 \leq i \leq N}
P_{U}(u_i) \prod\limits_{1 \leq j \leq M} P_{V}(v_j) \prod\limits_{1
  \leq i \leq N,1 \leq j \leq M} e^{g( Y_{ij} , u_i^\top v_j/\sqrt{N})} \,.
\label{Boltz_UV}
\end{equation}
Defined as above, this time $Y_{ij}\in {\mathbb R}^{N \times M}$ and $U \in {\mathbb R}^{N
  \times r}$,  $V \in {\mathbb R}^{M  \times r}$. Again we denote by
$u_i$, $v_j$ the {\it vector-spins} of dimension $r$ that collect rows of
matrices $U$, $V$. In this case the graph of interactions between
spins is bipartite.
}
\end{itemize}

The main motivation on this work is twofold. On the one hand, the
above mentioned probability measures are posterior probability
measures of an important class of high-dimensional inference problems
known as constrained low-rank matrix estimation. In what follows we
give examples of applications of these matrix estimation problems in
data processing and statistics. On the other hand, our motivation from
the physics point of view is to present a unified formalism providing
the (replica symmetric) solution for a large class of mean-field
vectorial spin models with disorder.

The general nature of the present work stems from the fact that the
probability distributions $P_X$, $P_U$, $P_V$ and the function $g$ are
very generic (assumptions are summarize in section \ref{sec:universality}). These functions can even depend on the node $i$ or edge
$ij$. For simplicity we will treat site-independent functions $P_X$,
$P_U$, $P_V$ and $g$, but the theory developed here generalizes very
straightforwardly to the site or edge dependent case. From a
statistical physics point of view the terms $P_X$, $P_U$, $P_V$ play a
role of generic local magnetic fields acting on the individual
spins. Distributions $P_X$, $P_U$, $P_V$ describe the nature of the
vector-spin variables and the fields that act on them. The simplest
example is the widely studied Ising spins for which $r=1$ and
$P_X(x)= \rho \delta(x-1)+ (1-\rho) \delta(x+1)$, where $\rho$ here
would be related to the usual magnetic field $h$ and inverse
temperature $\beta$ as $\rho = e^{\beta h}/(2\cosh{\beta h})$. In this
paper we treat a range of other examples with $r \ge 1$ and elements
of $x$ being both discrete or continuous. This involves for instance
spherical spin models with $P_X$ being Gaussian, or Heisenberg spins
where $r=3$ and each $x$ is confined to the surface of a sphere.

Denoting 
\be
w_{ij} =x_i^\top x_j/\sqrt{N} \, , \quad  {\rm or} \quad w_{ij} = u_i^\top v_j/\sqrt{N}  \label{wij}
\ee
according to symmetric or bipartite context, the terms $g(Y,w)$ are
then interactions between pairs of spins that depend only on the
scalar product between the corresponding vectors. The most commonly
considered form of interaction in statistical physics is simply 
\be
g(Y,w)=\beta Yw \label{can_Ham} \ee
 with $\beta$ being a constant
called inverse temperature, leading to a range of widely considered
models with pair-wise interactions. We will refer to this form of
function as the {\it conventional Hamiltonian}.

In order to complete the setting of the problem we need to specify how
is the quenched disorder $Y$ chosen. We will consider two main cases
of the quenched disorder defined below. We note that even for problems
where the matrix $Y$ is not generated by either of the below our
approach might still be relevant, e.g. for the restricted Boltzmann
machine that is equivalent to the above bipartite model with $Y$ that
were learned from data (see for instance
\cite{gabrie2015training,tramel2016inferring} and the discussion below).
\begin{itemize}
   \item{{\bf Randomly quenched disorder:} In this case the matrix elements
      of $Y$ are chosen independently at random from some probability
      distribution $P(Y_{ij})$. In this paper we will consider this
      distribution to be independent of $N$ and in later parts for
      simplicity we will restrict its mean to be zero. This case of
      randomly quenched disorder will encompass many well known and
      studied spin glass models such as the Ising spin glass,
      Heisenberg spin glass or the spherical spin glass or the
      Hopfield model.} 
   \item{{\bf Planted models:}  Concerning applications in data science this is the more
       interesting case and most of this paper will focus on it. In
       this case we consider that there is some {\it ground truth}
       value of $X_0 \in {\mathbb R}^{N \times r_0}$ (or $U_0 \in
       {\mathbb R}^{N \times r_0}$, $V_0 \in {\mathbb R}^{M \times
         r_0}$) with rows that are generated
       independently at random from some probability distribution
       $P_{X_0}$ (or $P_{U_0}$, $P_{V_0}$). Then the disorder $Y$ is
       generated element-wise as a noisy observation of the
       product $w^0_{ij}= {x^{0,\top}_{i}} x^0_{j} /\sqrt{N}$ (or $w^0_{ij} = {u^{0,\top}_{i}} v^0_{j}/\sqrt{N}$)  via an output {\it
         channel} characterized by the output probability distribution 
        $ P_{\rm out}(Y_{ij} | w^0_{ij})$.}
\end{itemize}

\subsection{Preliminaries on the planted setting}

\subsubsection{Bayes optimal inference}
\label{sec:Bayes_optimal}

Many applications, listed and analyzed below, in which the planted setting is relevant,
concern problems where we aim to infer some ground truth matrices $X_0$, $U_0$,
$V_0$ from the observed data $Y$ and from the information we have
about the distributions $P_{X_0}$, $P_{U_0}$, $P_{V_0}$ and $P_{\rm
  out}$. The information-theoretically optimal way of doing inference
if we know how the data $Y$ and how the ground-truth matrices were generated is to
follow the Bayesian inference and compute marginals of the
corresponding posterior
probability distribution. According to the Bayes formula, the posterior
probability distribution for the symmetric case is 
\begin{equation}
P(X|Y) = \frac{1}{Z_{X}(Y)}\prod\limits_{1 \leq i \leq N} P_{X_0}(x_i) \prod\limits_{1 \leq i < j \leq N} P_{\rm out}\left(Y_{ij} \middle| \frac{x_i^{\top} x_j}{\sqrt{N}}\right) \,.
\label{Prob_XX}
\end{equation}
For the bipartite case it is
\begin{equation}
P(U,V|Y) = \frac{1}{Z_{UV}(Y)}\prod\limits_{1 \leq i \leq N} P_{U_0}(u_i) \prod\limits_{1 \leq j \leq M} P_{V_0}(v_j) \prod\limits_{1 \leq i \leq N,1 \leq j \leq M} P_{\rm out}\left(Y_{ij} \middle| \frac{u_i^\top v_j}{\sqrt{N}}\right) \,.
\label{Prob_UV}
\end{equation}
Making link with the Boltzmann probability measures (\ref{Boltz_XX})
and (\ref{Boltz_UV}) we see that the {\it Bayes optimal inference} of the
planted configuration is equivalent to the statistical physics of the above
vector-spin models with 
\bea
    P_{X_0} = P_X\, ,  \quad P_{U_0} = P_U \, ,  \quad P_{V_0} = P_V\, , \quad
    P_{\rm out}(Y|w) = e^{g(Y,w)}\, .   \label{Bayes_optimal}
\eea
This approach is optimal in the sense that the statistical estimator
$\hat X$ computed from the data $Y$ that minimizes the expected mean-squared
error between the estimator $\hat X$ and the ground truth $X_0$ is
given by the mean of the marginal of variable $x_i$ in the probability
distribution (\ref{Prob_XX})
\be
    \hat x_i(Y) = \int {\rm d} x \, x \,  \mu_i(x)\, , \quad {\rm
      where} \quad   \mu_i(x)  = \int   P(X|Y) \prod_{\{x_j\}_{j\neq
        i}} {\rm d}x_j\, .   
\ee
Analogously for the bipartite case. 

In the Bayes-optimal setting defined by conditions
(\ref{Bayes_optimal}) the statistical physics analysis of the problem
presents important simplifications known as the Nishimori conditions
\cite{NishimoriBook,zdeborova2015statistical}, which will be largely used in the present
paper. 
These conditions can be proven and stated without the usage of the
methodology developed below, they are a direct consequence of the Bayesian
formula for conditional probability and basic properties of
probability distributions.  

Assume Bayes-optimality of the output channel, that is $P_{\rm out} =e^{g(Y,w)}$.
First let us notice that every probability distribution has to be normalized
\begin{equation}
\forall w,\int {\rm d}Y P_{\rm out}(Y|w) = 1 \,.
\end{equation}
By deriving the above equation with respect to $w$ one gets.
\begin{eqnarray}
\forall w,\int {\rm d}Y P_{\rm out}(Y|w) \frac{\partial
  g(Y,w)}{\partial w} &=& \mathbb{E}_{ P_{\rm out}(Y|w)}\left[
  \frac{\partial g(Y,w)}{\partial w} \right]= 0  \, , \label{E_1st_der}
\\
\forall w,\int {\rm d}Y P_{\rm out}(Y|w) \left[ \left(\frac{\partial
      g(Y,w)}{\partial w} \right)^2 +  \frac{\partial^2
    g(Y,w)}{\partial w^2} \right] &=& \mathbb{E}_{ P_{\rm
    out}(Y,w)}\left[  \left(\frac{\partial g(Y,w)}{\partial w}
  \right)^2 +  \frac{\partial^2 g(Y,w)}{\partial w^2} \right] = 0 \, \label{E_2nd_der}
. 
\end{eqnarray}
Anticipating the derivation in the following we also define the
inverse Fisher information of an output channel $P_{\rm out}$ at $w=0$
as
\begin{equation}
\frac{1}{\Delta} = \mathbb{E}_{P_{\rm out}(Y|w=0)}\left[ \left(\frac{\partial g}{\partial w}\right)_{Y,w=0}^2 \right] \,, \label{Define_Delta}
\end{equation}

Let us now assume Bayes-optimality also for the prior distributions,
i.e. $r=r_0$ and $P_{X_0} = P_X$,  $P_{U_0} = P_U$,  $P_{V_0} =
P_V$. Then under  averages over the posterior distribution
(\ref{Prob_XX}) or (\ref{Prob_UV}) we can replace the random variables
$X$, $U$, $V$ for the ground truth $X_0$, $U_0$, $V_0$ and vice versa.
To prove this let us consider the symmetric case, the generalization
to the bipartite case is straightforward. Take
$X,X_1,X_2$ to be three independent samples from the posterior
probability distribution $P(X|Y)$, eq. (\ref{Prob_XX}). We then consider some
function $f({A},{B})$ of two configurations of the variables
${A},{B}$. Consider the following two expectations 
\bea
    {\mathbb E}\left[f(X_1,X_2)\right] &= &\int  f(X_1,X_2) P(Y) 
    P(X_1|Y) P(X_2|Y)  {\rm d}X_1 \,  {\rm d}X_2 \,  {\rm d}Y\, ,  \\
 {\mathbb E}\left[f(X_0,X)\right] &=& \int  f(X_0,X)   P(X_0,X)  {\rm
   d}X\,  {\rm d}X_0 = \int  f(X_0,X)   P(X_0,X,Y)  {\rm
   d}X\,  {\rm d}X_0 {\rm d}Y \nonumber \\ &=& \int  f(X_0,X)   P(X|Y,X_0) P_{\rm out}(Y | X_0)  P_{0}(X_0) {\rm
   d}X\,  {\rm d}X_0 {\rm d}Y \, .
\eea 
where we used the Bayes formula. We further observe that $P(X|Y,X_0)=P(X|Y)$ because~$X$ is
independent of $X_0$ when conditioned on~$Y$. In the
Bayes optimal case, i.e. when $P(X)= P_{0}(X)$ and
$P(Y|X)= P_{\rm out}(Y|X)$, we then obtain 
\be 
   {\rm Bayes \, \,  optimal:} \quad \quad {\mathbb E}\left[f(X_1,X_2)\right]=   {\mathbb E}\left[f(X_0,X)\right]\, , \label{Nish_gen}
\ee
meaning that under
expectations there is no statistical difference between the ground
truth assignment of variables $X_0$ and an assignment sampled
uniformly at random from the posterior probability distribution (\ref{Prob_XX}). This
is a simple yet important property that will lead to numerous
simplifications in the Bayes optimal case and it will be used in
several places of this paper, under the name {\it
  Nishimori condition}. 

From the point of view of statistical physics of disordered systems
the most striking property of systems that verify the Nishimori
conditions is that there cannot be any {\it replica symmetry breaking} in
the equilibrium solution of these systems
\cite{nishimori2001absence,NishimoriBook,zdeborova2015statistical}. This
simplifies considerably the analysis of the Bayes-optimal
inference. Note, however, that metastable (out-of-equilibrium)
properties of Bayes-optimal inference do not have to satisfy the
Nishimori conditions and replica symmetry breaking might be needed for
their correct description (this will be relevant in the cases of first
order phase transition described in section \ref{sec:1st_order}). 


\subsubsection{Principal component analysis}

In the above probabilistic inference setting the Bayesian approach of
computing the marginals is optimal. However, in a majority of the interesting cases
it is computationally intractable (NP hard) to compute the marginals
exactly. In all low-rank estimation problems the method that (for the bipartite case) comes to
mind as first when we look for a low-rank matrix close to an observe
matrix $Y$ is the {\it singular value decomposition} (SVD) where a rank $r$
matrix $\tilde Y$ that minimizes the squared difference in computed
\be
      {\rm argmin}_{\tilde Y}  \left[ \sum_{i,j} (Y_{ij} - \tilde Y_{ij})^2
      \right] = \sum_{s=1}^r u_s \lambda_s v^\top_s \, ,  \label{PCA}
\ee
where $\lambda_s$ is the $s$th largest singular value of $Y$, and $u_s \in {\mathbb
R}^N$, $v_s \in {\mathbb R}^M$ are the corresponding left-singular and right-singular  vectors of the
matrix $Y$. The above property is know as the Eckart-Young-Mirsky
theorem \cite{eckart1936approximation}. In the symmetric case $Y_{ij}=Y_{ji}$ we simply replace the singular values by
eigenvalues and the singular vectors by the eigenvectors. 
The above unconstrained low-rank approximation of the matrix
$Y$, eq. ~(\ref{PCA}), is also often referred to as {\it principal component analysis} (PCA), because
indeed when the matrix $Y$ is interpreted as $N$ samples of
$M$-dimensional data then the right-singular vectors $v_s$ are
directions of greatest variability in the data. 

PCA and SVD are methods of choice when the measure of performance is
the sum of square differences between the observe and estimated
matrices, and when there are no particular requirements on the
elements of the matrices $U$, $V$ or $X$. 

The methodology developed in this paper for the planted probabilistic
models, generalizes to arbitrary cost function that can be expressed
as a product of element-wise terms $e^{g(Y_{ij},w_{ij})}$ and to
arbitrary constraints on the rows of the matrices $U$, $V$, $X$ as long
as they can be described by row-wise probability distributions
$P_U$, $P_V$, $P_X$. Systematically comparing our results to the
performance of PCA is useful because PCA is well known and many
researcher have good intuition about what are its strengths and
limitations in various settings.

\subsection{The large size limit, assumptions and channel universality}
\label{sec:universality}

In this article we focus on the thermodynamic limit where $N, M\to
\infty$ whereas $r=O(1)$, and $\alpha \equiv M/N= O(1)$ and
all the elements of $Y$, $X$, $U$ and $V$ are of order 1. The
functions $P_X$, $P_U$, $P_V$ and $g$ do not depend on $N$
explicitly. In the planted model also the distribution $P_{X_0}$,
$P_{U_0}$, $P_{V_0}$ and $P_{\rm out}$ do not depend on $N$
explicitly. The only other requirement we impose on the distributions  $P_X$,
$P_U$, $P_V$ and $P_{X_0}$, $P_{U_0}$, $P_{V_0}$ is that they all have
a finite second moment. 

The factor $1/\sqrt{N}$ in the second argument of the function $g$
ensures that the behaviour of the above models is non-trivial and that
there is an interesting competition between the number $O(N)$ of local magnetic fields
$P_X$, $P_U$, $P_V$ and the number of $O(N^2)$ interactions. To
physics readership familiar with the Sherrington-Kirkpatrick (SK) model
this $1/\sqrt{N}$ factor will be familiar because in the SK model the
interaction between the Ising spins that lead to extensive free energy
are also of this order (with mean that is of order $1/N$). This is
compared to the ferromagnetic Ising model on a fully connected lattice for which the interactions
leading to extensive free energy scale as $1/N$.  

For readers interested in inference problems, i.e. the planted
setting, the $1/\sqrt{N}$ factor is the scaling of the signal-to-noise
ratio for which inference of $O(N)$ unknown from $O(N^2)$ measurements
is neither trivially easy nor trivially impossible. In the planted
setting $Y$ can be viewed as a random matrix with a rank-$r$
perturbation. The regime where the eigenvalues of dense random matrices with
low-rank perturbations split from the
bulk of the spectral density is precisely when the strength of the
perturbation is $O(1/\sqrt{N})$, see e.g. \cite{baik2005phase}. 

We are therefore looking at statistical physics models with $O(N^2)$
pairwise interactions where each of the interactions depend only
weakly on the configuration of the vector-spins. As a consequence,
properties of the system in the thermodynamic limit $N\to \infty$
depend only weakly on the details
of the interaction function $g(Y_{ij},w_{ij})$ with $w_{ij}$ given by
(\ref{wij}). The results of this paper hold for every function $g$ for
which the following Taylor expansion is well defined
\begin{equation}
e^{g(Y_{ij},w_{ij})} = e^{g(Y_{ij},0)} \left\{1+ \frac{\partial
    g(Y_{ij},w)}{\partial w}\Big|_{w=0} w_{ij} + \left[
    \left(\frac{\partial g(Y_{ij},w)}{\partial w}\Big|_{w=0}\right)^2
    + \frac{\partial^2 g(Y_{ij},w)}{\partial w^2}\Big|_{w=0}\right] \frac{w_{ij}^2}{2} + O(w_{ij}^3)\right\} \,.
\end{equation}
In order to simplify the notation in the following we denote
\begin{eqnarray}
S_{ij} &\equiv& \frac{\partial
    g(Y_{ij},w)}{\partial w}\Big|_{w=0} \, ,
\label{Define_S}
\\
R_{ij} &\equiv&  \left(\frac{\partial g(Y_{ij},w)}{\partial w}\Big|_{w=0}\right)^2
    + \frac{\partial^2 g(Y_{ij},w)}{\partial w^2}\Big|_{w=0} \,.
\label{Define_R}
\end{eqnarray}
We will refer to the matrix $S$ as the Fisher score matrix. 
The above expansion can now be written in a more compact way
\begin{equation}
e^{g(Y_{ij},w_{ij})}= e^{g(Y_{ij},0)}   \left[1+ S_{ij} w_{ij} + \frac{R_{ij} w_{ij}^2}{2} +
  O(w_{ij}^3)\right]  = e^{  g(Y_{ij},0) + S_{ij}w_{ij} + \frac{1}{2}(R_{ij}-S^2_{ij}) w_{ij}^2 + O(w_{ij}^3) }\,. \label{ExpandP_out}
\end{equation}

Let us now analyze the orders in this expansion. In the Boltzmann
measure (\ref{Boltz_XX}) and (\ref{Boltz_UV}) the terms $e^{g(Y_{ij},w_{ij})}$ appears
in a product over $O(N^2)$ terms and $w = O(1/\sqrt{N})$. At the same
time only terms of order $O(N)$ in the exponent of the Boltzmann
measure influence the leading order (in $N$) of the marginals (local
magnetizations), therefore all the terms that depend on 3rd and higher
order of $w$ are negligible. This means that the leading order of the
marginals depend on the function $g(Y,w)$ only trough the matrices of its first and
second derivatives at $w=0$, denoted $S$ and $R$
(\ref{Define_S}-\ref{Define_R}). This means in particular that in
order to understand the phase diagram of a model with general $g(Y,w)$
we only need to consider one more order than in the conventional Hamiltonian
considered in statistical physics (\ref{can_Ham}).  

In the sake of specificity let us state here the two examples of the
output channels $g(Y,w)$ considered most prominently in this paper and
their corresponding matrices $S$ and $R$. The first example
corresponds to observations of low-rank matrices with additive
Gaussian noise, we will refer to this as the Gaussian output channel
\begin{eqnarray}
{\rm \bf Gaussian \, \, \,  Noise \, \, \,  Channel:}\quad \quad
g(Y,w) = \frac{-(Y-w)^2}{2 \Delta} - \frac{1}{2} \log{2\pi \Delta}\, , \quad \quad 
S_{ij} = \frac{Y_{ij}}{\Delta}\, , \quad \quad 
R_{ij}= \frac{Y_{ij}^2}{\Delta^2} - \frac{1}{\Delta} \, ,
\label{gauss_out}
\end{eqnarray}
where $\Delta$ is the variance of the noise, for the specific case of
the Gaussian output channel $\Delta$ is also the Fisher information as
defined in eq.~(\ref{Define_Delta}).
The second example is the one most often considered in physics
given by eq. (\ref{can_Ham})
\begin{eqnarray}
{\rm \bf Conventional \, \, \,  Hamiltonian:}\quad \quad
g(Y,w) = \beta Y w \, , 
\quad \quad
S_{ij} =\beta  Y_{ij} \, , 
\quad \quad 
R_{ij} = \beta^2 Y_{ij}^2\, .  \label{can_Ham_Example}
\end{eqnarray}
with $\beta$ being a constant called inverse temperature.
Another quite different example of the output channel will be given to
model community detection in networks in section \ref{ParagraphCommunityDetection}.

\subsection{Examples and applications}
The Boltzmann measures (\ref{Boltz_XX}) and (\ref{Boltz_UV}) together
with the model for the disorder $Y$ can be used to describe a range of
problems of practical and scientific interest studied
previously in physics and/or in data sciences. In this section
we list several examples and applications for each of the four
categories -- the symmetric and bipartite case, and the randomly quenched and
planted disorder.

\subsubsection{Examples with randomly quenched disorder}

\paragraph{\bf Sherrington-Kirkpatrick (SK) model.} The SK model
\cite{SherringtonKirkpatrick75} stands at the roots of the theory of
spin glasses. It can be described by the symmetric Boltzmann measure
(\ref{Boltz_XX}) with the conventional Hamiltonian $g(Y,w)=\beta Yw$. 

The $x_i$ are Ising spins, i.e. $x_i \in \{\pm 1\}$, with distribution
\be
P_X(x_i) =\rho \delta(x_i-1) + (1-\rho) \delta(x_i+1).   \label{Ising}
\ee
 The parameter
$\rho$ is related to the inverse temperature $\beta$ and an external magnetic
field $h$ as $\rho = e^{\beta h}/(2\cosh{\beta h})$. Note that the
parameter $\rho$ could also be site-dependent and our approach would
generalize, but in this paper we work with site independent functions
$P_X$. 

The elements of the (symmetric) matrix $Y_{ij}$ are the quenched random
disorder, i.e. they are generated independently at random from some
probability distribution. Most usually considered distributions of
disorder would be the normal distribution $Y_{ij} \sim {\cal N}(0,1)$, or binary $Y_{ij}=1$ with probability $1/2$ and
$Y_{ij}=-1$ otherwise. 

The algorithm developed in this paper for the general case corresponds
to the Thouless-Anderson-Palmer \cite{ThoulessAnderson77} equations
for the SK model. The theory developed here correspond to the replica
symmetric solution of \cite{SherringtonKirkpatrick75}. Famously this
solution is wrong below certain temperature where effects of replica
symmetry breaking (RSB) have to be taken into account. In this paper we
focus on the replica symmetric solution, that leads to exact and
novel phase diagrams for the planted models. The RSB solution in the
present generic setting will be presented elsewhere. We
present the form of the TAP equations in the general case encompassing
a range of existing works. 

\paragraph{\bf Spherical spin glass.}

Next to the SK model, the spherical spin glass model \cite{kosterlitz1976spherical} stands behind
large fraction of our understanding about spin glass. Mathematically
much simpler than the SK model this model stands as a prominent case in the development in
mathematical physics. The spherical spin glass is formulated via the
symmetric Boltzmann measure (\ref{Boltz_XX}) with the conventional
Hamiltonian $q(Y,w)=\beta Y w$. The function $P_X(x_i) = e^{- x_i^2/2}$
with $x_i \in {\mathbb R}$
enforces (canonically) the spherical constraint $\sum_i x^2_i
=N$. External magnetic field can also be included in $P_X(x_i)$. 

The disorder $Y_{ij}$ is most commonly randomly quenched in physics
studies of the spherical spin glass model. 

\paragraph{\bf Heisenberg spin glass.}

In Heisenberg spin glass \cite{sommers1981theory} the Hamiltonian is again the conventional
symmetric one
with randomly quenched disorder. The spins are 3-dimensional vectors, $x_i \in
{\mathbb R}^3$, of unit length, $x_i^{\top} x_i=1$. Magnetic field
influences the direction of the spin so that 
\be
 P_X(x_i) = e^{ \beta h^{\top} x_i } \, , 
\ee
where $h \in {\mathbb R}^3$. The more general $r$-component model was also
astudied extensively in the spin glass literature \cite{gabay1981coexistence}. 

\paragraph{\bf Restricted Boltzmann Machine.}

Restricted Boltzmann machines (RBMs) are one of the triggers of the recent revolution
in machine learning called deep learning
\cite{hinton2010practical,hinton2006fast,lecun2015deep}. The way RBMs
are used in machine learning is that one considers the bipartite
Boltzmann measure (\ref{Boltz_UV}). In the training phase one searches a matrix $Y_{ij}$
such that the data represented as a set of configurations of the
$u$-variable have high likelihood (low energy). The $v$-variable are
called the hidden units and columns of the matrix $Y_{ij}$ (each
corresponding to one hidden unit) are often interpreted as features
that are useful to explain the structure of the data. 

The RBM is most commonly considered for the conventional Hamiltonian
$g(Y,w)=Yw$ and for binary variables $u_i\in \{0,1\}$ and $v_i \in
\{0,1\}$. But other distributions for both the data-variables $u_i$
and the hidden variables $v_i$ were considered in the literature and
the approach of the present paper applies to all of them.  

We note that the disorder $Y_{ij}$ that was obtained for an RBM
trained on real datasets does not belong to the classes for which the
theory developed in this paper is valid (training introduces involved
correlations). However, it was shown recently that the Low-RAMP equations as studied in the present
paper can be used efficiently for training of the RBM \cite{kappen1998boltzmann,gabrie2015training}. 

The RBM with Gaussian hidden variables is related to the well known
Hopfield model of associative memory \cite{hopfield1982neural}. Therefore the properties of
the bipartite Boltzmann measure (\ref{Boltz_UV}) with a randomly
quenched disorder $Y_{ij}$ are in one-to-one correspondence with the
properties of the Hopfield model. This relation in the view of the TAP equations
was studied recently in \cite{mezard2016mean}.

\subsubsection{Examples with planted disorder}

So far we covered examples where the disorder was randomly quenched
(or more complicated as in the RBM). The next set of examples involves
the planted disorder that is more relevant for applications in signal
processing or statistics, where the variables $X$, $U$, $V$ represent
some signal we aim to recover from its measurements $Y$. Sometimes it
is the low-rank matrix $w_{ij}$ that we aim to recover from its noisy
measurements $Y$. In the literature the general planted problem can be called
low-rank matrix factorization, matrix recovery, matrix denoising or
matrix estimation.

\paragraph{\bf Gaussian estimation}

The most elementary examples of the planted case is when the
measurement channel is Gaussian as in eq. (\ref{gauss_out}), and the
distributions $P_X$, $P_U$ and $P_V$ are also Gaussian i.e. 
\bea
P_X(x_i) &=& \frac{1}{\sqrt{{\rm Det}(2 \pi \sigma_X)}} e^{ - \frac{1}{2}
  (x_i - \mu_X)^\top \sigma_X^{-1}(x_i - \mu_X) } \,, \label{gauss_x} \\ 
P_U(u_i) &=& \frac{1}{\sqrt{{\rm Det}(2 \pi \sigma_U)}} e^{-\frac{1}{2}
  (u_i - \mu_U)^\top \sigma_U^{-1}(u_i - \mu_U) } \,  \label{gauss_u},
\\
P_V(v_i) &=& \frac{1}{\sqrt{{\rm Det}(2 \pi \sigma_V)}} e^{-\frac{1}{2}
  (v_i - \mu_V)^\top \sigma_V^{-1}(v_i - \mu_V) }\, , \label{gauss_v}
\eea 
where $\mu_X,\mu_U,\mu_V \in  {\mathbb R}^r$ are the means of the
distributions, $\sigma_X,\sigma_U,\sigma_V \in {\mathbb R}^{r\times
  r}$ are the covariance matrices, $Y_{ij}, w_{ij} \in {\mathbb R}$
with $w_{ij}$ being given by (\ref{wij}).

We speak about the estimation problem as Bayes-optimal Gaussian estimation if the disorder
$Y_{ij}$ was generated according to 
\be
               P_{\rm out}(Y_{ij}| w^0_{ij}) = e^{g(Y_{ij},w^0_{ij})}
               \, ,
\ee 
where $g(Y,w)$ is given by eq. (\ref{gauss_out}), and 
\be
        w^0_{ij} ={x^0}_i^\top x^0_j/\sqrt{N} \, , \quad  {\rm or}
        \quad w^0_{ij} = {u^0_i}^\top v^0_j/\sqrt{N} \, .
\ee
with $X_0$, $U_0$, and $V_0$ being generated from probability
distributions $P_{X_0}=P_X$, $P_{U_0}=P_U$, $P_{V_0}=P_V$. The goal is
to estimate matrices $X_0$, $U_0$, and $V_0$ from $Y$. 

\paragraph{\bf Gaussian Mixture Clustering}

Another example belonging to the class of problems discussed in
this paper is the model for Gaussian mixture clustering. In this case the
spin variables $u_i$ are such that 
\be
   P_U(u_i) = \sum_{s=1}^r    n_s   \delta(u_i - e_s) \, , \label{Gauss_mix_prior}
\ee
where $r$ is the number of clusters, and $e_s$ is a unit
$r$-dimensional vector with all components except $s$ equal to zero,
and the $s$-component equal to 1, e.g. for $r=3$ we have
$e_1=(1,0,0)^\top$, $e_2=(0,1,0)^\top$, $e_s=(0,0,1)^\top$. Having
$u_i = e_s$ is interpreted as data points $i$ belongs to cluster
$s$. We have $N$ data points. 

The columns of the matrix $V$ then represent centroids of each of
the $r$ clusters in the $M$-dimensional space. The distribution $P_V$
can as an example take the Gaussian form (\ref{gauss_v}) with the covariance
$\sigma_V$ being an identity and the mean $\mu_V$ being zero. 
The output channel is Gaussian as in (\ref{gauss_out}). All together
this means the $Y_{ij}$ collects positions of $N$ points in $M$
dimensional space that are organized in $r$ Gaussian clusters. 
The goal is to estimate the centers of the clusters and the cluster
membership from $Y$. 

Standard algorithms for data clustering incluse those based on the
spectral decomposition of the matrix $Y$ such as principal component
analysis \cite{hastie2005elements,wasserman2013all}, or Loyd's
$k$-means \cite{lloyd1982least}. Works on Gaussian mixtures that are
methodologically closely related to the present paper include
application of the replica method to the case of two clusters $r=2$ in
\cite{watkin1994optimal,barkai1994statistical,biehl1994statistical} or
the AMP algorithm of \cite{NIPS2013_5074}. Note that for two clusters
with the two centers being symmetric around the origin, the resulting
Boltzmann measure of the case with randomly quenched disorder is
equivalent to the Hopfield model as treated e.g. in
\cite{mezard2016mean}. 

Note also that there are interesting variants of the Gaussian mixture
clustering such as subspace clustering \cite{parsons2004subspace}
where only some of the $M$ directions are relevant for the
clustering. This can be modeled by a prior on the vectors $v_i$ that have a non-zero weight of $v_i$ being the null vector. 

The approach described in the present paper on the Bayes-optimal
inference in the Gaussian mixture clustering problem has been used in
a work of the authors with other collaborators in the work presented
in \cite{lesieur2016phase}.

\paragraph{\bf Sparse PCA}

Sparse  Principal  Component  Analysis  (PCA) \cite{johnstone2004sparse,zou2006sparse} is  a
dimensionality  reduction  technique  where  one  seeks  a  low-rank  representation  of  a  data  matrix  with  additional  sparsity
constraints  on  the  obtained  representation. The motivation is that
the sparsity helps to interpret the resulting representation. 
Formulated within the above setting sparse PCA corresponds to the
bipartite case (\ref{Boltz_UV}). The variables $U$ is considered unconstrained, as an example one often
considers a Gaussian prior on $u_i$ (\ref{gauss_u}). The variables $V$ are such that many
of the matrix-elements are zero. 

In the literature the sparse PCA problem was mainly considered in the
rank-one case $r=1$ for which a series of intriguing results was
derived. The authors of \cite{johnstone2004sparse} suggested an
efficient algorithm called diagonal thresholding that solves the sparse PCA
problem (i.e. estimates correctly the position and the value of the
non-zero elements of $U$) whenever the number of data samples is $N >
C K^2 \log{M}$ \cite{amini2008high}, where $K$ is the number of
non-zeros and $C$ is some constant. More recent works show existence
of efficient algorithm that only need $N > \hat C K^2$ samples \cite{deshpande2014sparse}.
For very sparse systems, i.e. small $K$, this is a striking
improvement over the  conventional PCA that would need $O(M)$
samples. This is why sparsity linked together with PCA brought
excitement into data processing methods. At the same time, this result
is not as positive as it may seem, because by searching exhaustively
over all positions of the non-zeros the correct  support can be
discovered with high probability with number of samples $N > \tilde C K \log{M}$. 

Naively one might think that polynomial algorithms that need less
thank $O(K^2)$ samples might exist and a
considerable amount of work was devoted to their search without
success. However, some works suggest that perhaps polynomial algorithm
that solve the sparse PCA problems for number of samples $N < O(K^2)$
do not exist. Among the most remarkable one is the work
\cite{krauthgamer2015semidefinite} showing that the SDP algorithm,
that is otherwise considered rather powerful, fails in this
regime. The work of \cite{berthet2013computational} goes even further
showing that if sparse PCA can be solved for $N < O(K^2)$ then also a
problem known as the planted clique problem can be solved in a regime
that is considered as algorithmically hard for already several
decades.

The problem of sparse PCA is hence one of the first examples of a
relatively simple to state problem that
currently presents a wide barrier between computational and
statistical tractability. Deeper
description of the origin of this barrier is likely to shed light on
our understanding of typical algorithmic complexity in a broader scope. 

The above works consider the scaling when $N\to \infty$ and $K$ is
fixed (or growing slower than $O(N)$). A regime natural to many
applications is when $K = \rho N$ where $\rho = O(1)$. This regime was
considered in \cite{DeshpandeM14} where it was shown that for $\rho>
\rho_0 \approx 0.04139$ an efficient algorithm that achieves the information
theoretical performance exists. This immediately bring a question of
what exactly happens for $\rho < \rho_0$ and how does the barrier
described above appear for $K \ll N$? This question was illuminated in
a work by the present authors \cite{lesieur2015phase} and will be developed
further in this paper.  

We consider sparse PCA in the bipartite case, with Gaussian $U$ (\ref{gauss_v})
and sparse $V$
\be
      P_V(v_i) = \rho \delta(v_i-1) + (1-\rho) \delta(v_i) \, ,  \label{Bernoulli_V}
\ee
as corresponds to the formulation of
\cite{johnstone2004sparse,zou2006sparse,amini2008high,berthet2013computational} and others. 
This probabilistic setting of sparse PCA was referred to as {\it spiked
Wishart model} in \cite{DeshpandeM14}, notation that we will adopt in
the present paper. This model is also equivalent to the one studied
recently in \cite{monasson2015estimating} where the authors simply
integrate over the Gaussian variables. 

In \cite{DeshpandeM14} the authors also considered a symmetric variant
of the sparse PCA, and refer to it as the {\it spiked Wigner model}. The
spiked Wigner model is closer to the planted clique problem, that
can be formulated using (\ref{Boltz_XX}) with $X$ having many zero
elements. In the present work we will consider several models for the
prior distribution $P_X$. The
{\it Bernoulli model} as in \cite{DeshpandeM14} 
where 
\be
      {\rm \bf Bernoulli \, \, \, model:} \quad \quad  P_X(x_i) =  \rho \delta(x_i-1) + (1-\rho) \delta(x_i) \, .  \label{Bernoulli_X}
\ee

The spiked Bernoulli model can also be interpreted as a problem of submatrix
localization where a submatrix of size $\rho N \times \rho
N$ of the matrix $Y$ has a larger mean than a randomly chosen
submatrix. The submatrix localization is also relevant in the
bipartite case, where it has many potential applications. The most
striking ones being in gene expression where large-mean submatrices of the
matrix of gene expressions of different patients may correspond to
groups of patients having the same type of disease \cite{madeira2004biclustering,cheng2000biclustering}. 

In this paper we will also consider the {\it spiked Rademacher-Bernoulli model} with 
\be
    {\rm \bf Rademacher-Bernoulli \, \, \, model:} \quad \quad P_X(x_i) =  \frac{\rho}{2} \left[ \delta(x_i-1)  + \delta(x_i +1)\right]+ (1-\rho) \delta(x_i) \, ,  
    \label{P_X_Rademacher}
\ee
as well as the {\it spiked Gauss-Bernoulli}
\be
     {\rm \bf Gauss-Bernoulli \, \, \, model:} \quad \quad P_X(x_i) =  \rho {\cal N}(x_i,0,1) + (1-\rho) \delta(x_i) \, ,  
\ee
where ${\cal N}(x_i,0,1)$ is the Gaussian distribution with zero mean
and unit variance. 

So far we discussed the sparse PCA problem in the case of rank one,
$r=1$, but the case with larger rank is also interesting, especially
in the view of the question of how does the algorithmic barrier depend
on the rank. To investigate this question in \cite{lesieur2015phase} we
also considered the jointly-sparse PCA, where the whole $r$-dimensional lines of $X$
are zero at once, the non-zeros are Gaussians of mean $\vec{0}$ and
covariance being the identity. Mathematically, $x_i \in {\mathbb R}^r$ with
\begin{equation}
P_X(x_i) = \frac{\rho}{(2 \pi)^{r/2}}e^{\frac{-x^\top x}{2}} + (1-\rho)\delta(x_i)\,. \label{P_X_JointedGaussBernoulli}
\end{equation}

Another example to consider is the independently-sparse
PCA where each of the $r$ components of the lines in $X$ is taken
independently from the Gauss-Bernoulli distribution, for $x_i \in
{\mathbb R}^r$ we have then
\begin{equation}
P_X(x_i) = \prod\limits_{1 \leq k \leq r} \left[
\frac{\rho}{\sqrt{2 \pi}}e^{\frac{-x_{ik}^2}{2}} + (1-\rho)\delta(x_{ik})\right] \,.\label{P_X_GaussBernoulli}
\end{equation}

\paragraph{\bf Community detection}
\label{ParagraphCommunityDetection}
Detection of communities in networks is often modeled by the
stochastic block model (SBM) where pairs of nodes get connected with
probability that depends on the indices of groups to which the two
nodes depend. Community detection is a widely studied problem, see
e.g. the review \cite{fortunato2010community}. Studies of statistical and computationally
barriers in the SBM recently became very active in mathematics and statistics
starting perhaps with a series of statistical physics conjectures
about the existence of phase transitions in the problem
\cite{decelle2011inference,decelle2011asymptotic}. The particular
interest of those works is that they are set in the sparse regime
where every node has $O(1)$ neighbors. 

Also the dense regime where every node has $O(N)$ neighbors is
theoretically interesting when the difference between probabilities to
connect depends only weakly on the group membership. The relevant
scaling is the same as in the Potts glass model
studied in \cite{GrossKanter85}. In fact the dense community detection
is exactly the planted version of this Potts glass model. In the
setting of the present model (\ref{Boltz_XX}) the community detection
was already considered in \cite{deshpande2016asymptotic} for two symmetric groups, in \cite{lesieur2015mmse},
and in \cite{barbier2016mutual}. In the present paper we detail the results reported
briefly in \cite{lesieur2015mmse} and in \cite{barbier2016mutual}. We
consider the case with a general number of equal sized groups, the
symmetric case. And also a case with two groups of different sizes, but
such that the average degree in each of the groups is the same.

To set up the dense community detection problem we consider a network
with $N$ nodes. Each node $i$ belongs to
a community indexed by $t_i \in \{1, \cdots ,r \}$.
For each pair $(i,j)$ we create an edge with probability $C_{t_i
  t_j}$. Where $C$ is an $r \times r$ matrix called the connectivity
matrix. In the examples of this paper we will consider two special cases of
the community detection problem.

One example with $r$ symmetric equally
sized groups where 
for each pair of nodes $(i,j)$ we create an edge between the two nodes with probability $p_{\rm in}$ if they are in the same group and with probability $p_{\rm out}$ if not:
\begin{eqnarray}
C &=& \begin{pmatrix}
p_{\rm in} & p_{\rm out} & \cdots & p_{\rm out}\\
p_{\rm out} & \ddots & \ddots & \vdots \\
\vdots & \ddots &\ddots & p_{\rm out}\\
p_{\rm out} & \cdots & p_{\rm out} & p_{\rm in}
\end{pmatrix} = p_{\rm in}I_r + p_{\rm out}(J_r - I_r)
\label{ConnectivityMatrix} \,,
\end{eqnarray}
where $I_r$ is a $r$-dimensional identity matrix, and $J_r$ is a
$r$-dimensional matrix filled with unit elements. 
The scaling we consider here is 
\begin{eqnarray}
p_{\rm out}&=& O(1) \, , \quad p_{\rm in} = O(1)\, ,
\label{Community_Pin_Pout_1}
\\
|p_{\rm in}-p_{\rm out}| &=& \frac{\mu}{\sqrt{N}} \, , 
\quad \mu = O(1) \label{Community_Pin_Pout_2} \, ,
\end{eqnarray}
so that the average degree in the graph is extensive. Note, however,
that by rather generic correspondence between diluted and dense
models, that has been made rigorous recently
\cite{deshpande2016asymptotic,caltagirone2016recovering}, the results derived in this case
hold even for average degrees that diverge only very mildly with $n$. 
The goal is to infer the group index to which each node belongs purely from
the adjacency matrix of the network (up to a permutation of the indices).
This  problem is transformed into the low-rank matrix factorization
problem through the use of the following prior probability
distribution
\begin{equation}
P_{X}(x_i) = \frac{1}{r}\sum\limits_{s = 1}^r \delta(x_i- {e}_s) \,. \label{SBM_P_X}
\end{equation}
where $e_s \in {\mathbb R}^r$ is the vector with $0$ everywhere except
a $1$ at position $s$. Eq. (\ref{SBM_P_X}) is just a special case of (\ref{Gauss_mix_prior}).
The output channel that describes the process of creation of the graph is 
\begin{eqnarray}
P_{\rm out}(Y_{ij} = 1 | x^\top_i x_j /\sqrt{N}) = p_{\rm out} + \frac{\mu
  x^\top_i x_j}{\sqrt{N}} \, ,\label{SBM_out1}
\\
P_{\rm out}(Y_{ij} = 0 | x^\top_i x_j /\sqrt{N}) = 1 - p_{\rm out} -
\frac{\mu x^\top_i x_j}{\sqrt{N}} \,.  \label{SBM_out2}
\end{eqnarray}
Next to the conventional Hamiltonian (\ref{can_Ham_Example}) and the Gaussian noise
(\ref{gauss_out}), the SBM output (\ref{SBM_out1}-\ref{SBM_out2}) is a
third example of an output channel that we consider in this article. It
will be used to illustrate the simplicity that arises due to the
channel universality, as also considered in \cite{deshpande2016asymptotic} and
\cite{lesieur2015mmse}. 
Here, we obtain for the output matrices 
\begin{eqnarray}
S_{ij}(Y_{ij}=1) &= \frac{\mu}{p_{\rm out}},\quad S_{ij}(Y_{ij}=0) =& \frac{-\mu}{1-p_{\rm out}} \,,
\\
R_{ij}(Y_{ij}=1) &= 0,\quad  \:
R_{ij}(Y_{ij}=0) =& 0 \,.
\end{eqnarray}
Here $\mu$ is parameter that can be used to fix the signal to noise ratio.


Another example of community detection is the one with two balanced communities,
i.e. having different size but the same average degree. In that setting there are two communities of size $\rho n$ and $(1-\rho)n$ with $\rho \in [0;1]$. The connectivity matrix of this model is given by
\begin{equation}
C = \begin{pmatrix}
p_{\rm out} & p_{\rm out} \\
p_{\rm out} & p_{\rm out}
\end{pmatrix}+
\frac{\mu}{\sqrt{N}} 	\begin{pmatrix}
\frac{1-\rho}{\rho} & -1 \\
-1 & \frac{\rho}{1-\rho}
\end{pmatrix}
\label{ConnectivityMatrix_Balanced_Intro} \,.
\end{equation}
This can be modelled at the symmetric matrix factorization with rank $r=1$ and
the prior given as
\begin{eqnarray}
P_{X}(x) &=& \rho\delta\left(x-\sqrt{\frac{1-\rho}{\rho}}\right) +
             (1-\rho)\delta\left(x+\sqrt{\frac{\rho}{1-\rho}}\right)
             \, .
\label{Proba_Cliques}
\end{eqnarray}
The values in $C$ are chosen so that each community has an average
degree of $p_{\rm out}N$.
The fact that in both of these cases each community has the same average degree means that one can not hope to just use the histogram of degrees to make the distinction between the communities.
The output channel here is identical to the one given in (\ref{SBM_out1}-\ref{SBM_out2})

A third example of community detection is locating one denser community
in a dense network, as considered in \cite{montanari2015finding}
(specifically the large degree limit considered in that paper). We
note that thanks to the output channel universality (Sec.~\ref{sec:universality}) this
case is equivalent to the spiked Bernoulli model of symmetric sparse PCA. 

As a side remark we note that the community detection setting as
considered here is also relevant in the bipartite matrix factorization
setting where it
becomes the problem of biclustering
\cite{madeira2004biclustering,cheng2000biclustering}. The analysis developed in this paper can be straightforwardly extended to the
bipartite case.

\subsection{Main results and related work}
\label{Main_results}
The present paper is built upon two previous shorter papers by the
same authors \cite{lesieur2015phase,lesieur2015mmse}. it focuses on
the study of the general type of models described by probability
measures (\ref{Boltz_XX}) and (\ref{Boltz_UV}). On the one hand, these
represent Boltzmann measures of vectorial-spin systems on fully
connected symmetric or bipartite graphs. Examples of previously
studied physical models that are special cases of the setting
considered here would be the Sherrington-Kirkpatrick model
\cite{SherringtonKirkpatrick75}, the Hopfield model
\cite{hopfield1982neural,mezard2016mean}, the inference (not learning)
in the restricted Boltzmann machine
\cite{gabrie2015training,mezard2016mean,tubiana2016emergence}.  Our
work hence provides a unified replica symmetric solution and TAP
equations for generic class of prior distributions and Hamiltonians.

On the other hand, these probability distributions (\ref{Boltz_XX})
and (\ref{Boltz_UV}) represent the posterior probability distribution
of a low-rank matrix estimation problem that finds a wide range of
applications in high-dimensional statistics and data analysis. The
thermodynamic limit $M/N=\alpha$ with $\alpha = O(1)$ whereas
$N,M\to \infty$ was widely considered in the studies of spin glasses,
in the context of low-rank matrix estimation this limit correspond to
the challenging high-dimensional regime, whereas traditional
statistics considers the case where $M/N \gg 1$.  We focus on the
analysis of the phase diagrams and phase transitions of low-rank
matrix estimation corresponding to the Bayes optimal matrix
estimation. We note that because we assume the data $Y$ were generated
from random factors $U$, $V$, $X$ we obtain much tighter, including
the constant factors, control of the high-dimensional behaviour
$N,M\to \infty$, $\alpha= M/N = O(1)$,  than some traditional bounds
in statistics that aim not to assume any generative model but instead
craft proofs under verifiable conditions of the observed data $Y$.

We note at this point that methodologically closely related series of
work on matrix factorization is concerned with the case of high rank,
i.e. when $r/N = O(1)$
\cite{krzakala2013phase,parker2014bilinear,kabashima2014phase}. While
that case also has a set of important application (among then the
learning of overcomplete dictionaries) it is different from the
low-rank case considered here. The theory developed for the high rank
case requires the prior distribution to be separable
component-wise. The high rank case also does not present any known output
channel universality, the details of the channel enter explicitly the
resulting state evolution. Whereas the low-rank case can be viewed as
a generalization of a spin model with pairwise interaction, in the
graphical model for the high-rank case the interactions involve $O(N)$
variables.

No attempt is made at mathematical rigor in the present article.
Contrary to 
the approach traditionally taken in
statistics 
most of the analysis in this paper is non-rigorous, and rather relies
of the accepted assumptions of the replica and cavity methods.  It it
worth, however, mentioning that for the case of Bayes-optimal
inference, a large part of the results of this paper were proven
rigorously in a recent series of works
\cite{rangan2012iterative,javanmard2013state,DeshpandeM14,krzakala2016mutual,deshpande2016asymptotic,barbier2016mutual,LelargeMiolane16,miolane2017fundamental}. These
proofs include the mutual information (related to the replica free
energy) in the Bayes-optimal setting and the corresponding minimum
mean-squared-error (MMSE), and the rigorous establishment that the
state evolution is indeed describing asymptotic evolution of the
Low-RAMP algorithm. The study out of the Bayes-optimal conditions
(without the Nishimori conditions) are move involved.

It has become a tradition in related literature
\cite{zdeborova2015statistical} to conjecture that the performance of
the Low-RAMP algorithm cannot be improved by other polynomial
algorithms. We do analyze here in detail the cases where Low-RAMP does
not achieve the MMSE, and we remark that since effects of replica
symmetry breaking need to be taken into account when evaluating the
performance of the best polynomial algorithms, the conjecture of the Low-RAMP
optimality among the polynomial algorithms deserves further detailed
investigation. 

This section gives a brief summary of our main results and their
relation  to existing work. 
\begin{itemize}
\item {\bf Approximate Message Passing} for Low-Rank matrix estimation
  ( {\bf Low-RAMP}): In section \ref{Sec:Low-RAMP} we derive and detail the
  approximate message passing algorithm to estimate marginal
  probabilities of the probability measures (\ref{Boltz_XX}) and
  (\ref{Boltz_UV}) for general prior distribution, rank and
  Hamiltonian (output channel). We describe various special case of
  these equations that arise due to the Nishimori conditions or due
  to self-averaging. In the physics literature this would be the TAP
  equations \cite{ThoulessAnderson77} generalized to vectorial spins
  with general local magnetic fields and generic type of pairwise
  interactions. The Low-RAMP equations encompass as a special case the
  original TAP equations for the Sherrington-Kirkpatrick model, TAP
  equations for the Hopfield model
  \cite{MezardParisi87b,mezard2016mean}, or the restricted Boltzmann
  machine
  \cite{gabrie2015training,tramel2016inferring,mezard2016mean,tubiana2016emergence}. Within
  the context of low-rank matrix estimation, the AMP equations were
  discussed in
  \cite{rangan2012iterative,NIPS2013_5074,DeshpandeM14,lesieur2015phase,lesieur2015mmse,
    deshpande2016asymptotic,lesieur2016phase}. Recently the Low-RAMP
  algorithm was even generalized to spin-variables that are not real
  vectors but live on compact groups \cite{perry2016message}.

     AMP type of algorithm is a promising alternative to gradient
     descent type of methods to minimize the likelihood. One of the
     main advantage of AMP is that it provides estimates of uncertainty
     which is crucial for accessing reliability and interpretability
     of the result. Compared to other Bayesian algorithms, AMP tends
     to be faster than Monte-Carlo based algorithms and more precise
     than variational mean-field based algorithms. 

     We distribute two open-source Julia and Matlab versions of
     LowRAMP at
     \url{http://krzakala.github.io/LowRAMP/}. We
     strongly encourage the reader to download, modify, and 
     improve on it.

   \item In physics, message passing equations are always closely
     linked with the  {\bf Bethe free energy} who's stationary points are the
     message passing fixed point equations. In the presence of
     multiple fixed points it is the value of the Bethe free energy
     that decides which of the fixed points is the correct one. In
     section \ref{PlefkaSection} and appendix \ref{app:Plefka} we
     derive the Bethe free energy on a single instance of the low-rank matrix
     estimation problem.  The form of free energy that we
     derive has the convenience to be variational in the sense that in
     order to find
     the fixed point we are looking for a maximum of the free energy,
     not for a saddle. Corresponding free energy for the compressed
     sensing problem was derived in \cite{krzakala2014variational,DBLP:journals/corr/abs-1301-6295} and can be used to guide
     the iterations of the Low-RAMP algorithm \cite{DBLP:journals/corr/VilaSRKZ14}.

   \item In section \ref{StateEvolutionSE} we derive the general form of the  {\bf state evolution} (under the
     replica symmetric assumption) of the Low-RAMP algorithm,
     generalizing previous works, e.g. \cite{rangan2012iterative, DeshpandeM14,lesieur2015phase,lesieur2015mmse}. We present simplifications for the Bayes-optimal
     inference and for the conventional form of the Hamiltonian. We
     also give the corresponding expression for the free energy. We
     derive the state evolution and the free energy using both the
     cavity and the replica method. 

     For the Bayes-optimal setting the
     replica Bethe free energy is up to a simple term related to the
     mutual information from which one can deduce the value of
     the  {\bf minimum information-theoretically achievable mean squared
     error}. Specifically, the MMSE correspond to the global maximum of
     the replica free energy (defined here with the opposite sign than
     in physics), the performance of Low-RAMP correspond to the
     maximum of the replica free energy that has the highest MSE. 
  
   We stress here that Low-RAMP algorithm belongs to the same class of
     approximate Bayesian inference algorithms as generic type of
     Monte Carlo Markov chains of variational mean-field methods. Yet
     Low-RAMP is very particular compared to these other two because of the
     fact that on a class of random models considered here its performance can be
     analyzed exactly via the state evolution and (out of
     the hard region) Low-RAMP asymptotically matches the performance
     of the Bayes-optimal estimator. Study of AMP-type of algorithms
     hence opens a way to put the variational mean field algorithms
     into more theoretical framework.

  \item We discuss the  {\bf output channel universality} as known in special
    cases in statistical physics (replica solution of the SK model
    depends only on the mean and variance of the quenched disorder
    not on other details of the distribution) and statistics
    \cite{deshpande2016asymptotic} (for the two group stochastic block
    model). The general form of this
    universality was first put into light for the Bayes-optimal estimation in \cite{lesieur2015mmse},
    proven in \cite{krzakala2016mutual}, in this paper we discuss this
    universality out of the Bayes-optimal setting. 

  \item{In section \ref{sec:PCA_analysis} we show that the state evolution with a Gaussian prior can be
      use to analyse the asymptotic  {\bf performance of spectral algorithms}
    such as PCA (symmetric case) or SVD (bipartite case) and derive
    the corresponding spectral mean-squared errors and phase
    transitions as studied in the random matrix literature \cite{baik2005phase}. For a
    recent closely related discussion see \cite{perry2016optimality}.}

   \item{In section \ref{sec:uniform} and \ref{sec:1st_order} we discuss the  {\bf typology of phase transition
     and phases} that arise in Bayes-optimal low-rank matrix
     estimation. We provide sufficient criteria for existence of phases where
     estimation better than random guesses from the prior distribution
     is not information-theoretically possible. We analyze linear
     stability of the fixed point of the
     Low-RAMP algorithm related to this phase of undetectability, to
     conclude that the threshold where Low-RAMP algorithm starts to
     have better performance than randomly guessing from the prior
     agrees with the spectral threshold known in the literature on
     low-rank perturbations of random matrices.  We also provide sufficient criteria
     for existence of first order phase transitions related to
     existence of phases where information-theoretically optimal
     estimation is not achievable by existing algorithms. And analyze
     the three thresholds $\Delta_{\rm Alg}$, $\Delta_{\rm IT}$ and
     $\Delta_{\rm Dyn}$ related to the first order phase transition.}

\item{In section \ref{sec:Phase_Diagrams} we give a number of examples of  {\bf phase diagrams}
  for the following models: rank-one $r=1$ symmetric $XX^\top$ case with prior
  distribution being Bernoulli, Rademacher-Bernoulli,
     Gauss-Bernoulli, corresponding to balanced 2-groups. For generic
     rank $r \ge 1$ we give the phase diagram for the symmetric
     $XX^\top$ case for the jointly-sparse PCA, and for the symmetric
     community detection. }

   \item Section \ref{sec:small_rho} and appendix \ref{AppendixSmallRho} is devoted to  {\bf small $\rho$ analysis} of the
     above models. This is motivated by the fact that in most existing
     literature with sparsity constraints the number of non-zeros is
     usually considered to be a 
     vanishing fraction of the system size. In our analysis the
     number of non-zeros is a finite fraction $\rho$ of the system
     size, we call the the regime of {\it linear sparsity}. We investigate whether the $\rho \to 0$ limit corresponds to previously
     studied cases.  What concerns the information-theoretically
     optimal performance and related threshold $\Delta_{\rm IT}$ our
     small $\rho$ limit agrees with the results known for sub-linear sparsity. Concerning the performance of efficient
     algorithms, from our analysis we conclude that for linear
     sparsity in the leading order the existing algorithms do not beat
     the threshold of the naive spectral method. This correspond to
     the known results for the planted dense subgraph problem (even
     when the size of the planted subgraph is sub-linear). However,
     for the sparse PCA problem with sub-linear sparsity algorithms
     such as covariance thresholding are known to beat the naive
     spectral threshold \cite{deshpande2014sparse}. In the regime of linear sparsity we do not
     recover such behaviour, suggesting that for linear sparsity,
     $\rho=O(1)$, efficient algorithms that take advantage of the
     sparsity do not exist. 

   \item In section \ref{sec:large_r} and appendix
     \ref{appendix:AppendixSpinodaleNetwork} we discuss analytical
     results that we can obtain for the community
     detection, and joint-sparse PCA models in the  {\bf limit of large
     rank $r$}. These large rank results are matching the
     rigorous bounds derived for these problems in
     \cite{banks2016information}. 
\end{itemize}

\section{Low-rank approximate message passing, Low-RAMP}
\label{Sec:Low-RAMP}

In this section we derive the approximate message passing for the
low-rank matrix estimation problems, i.e. estimation of marginals of
probability distributions (\ref{Boltz_XX}) and (\ref{Boltz_UV}). We
detail the derivation for the symmetric case  (\ref{Boltz_XX}) and
state the resulting equations for the bipartite case
(\ref{Boltz_UV}). We derive the algorithm in several stages. 
We aim to derive a tractable algorithm that under certain conditions (replica
symmetry, and absence of phase coexistence related to a first order
phase transitions) estimates marginal probabilities of the
probability distribution (\ref{Boltz_XX}) exactly in the large size
limit. We start with the belief propagation (BP) algorithm, as well
explained e.g. in \cite{Yedidia:2003}. 

It should be noted at this point that the BP algorithm is guaranteed
to be exact for a general system only when
the underlying graphical model is a tree. Yet it is often applied to
locally-tree like graphical models. On tree-like graphical models
usage of BP is justified when correlations decay fast enough so that
far-away boundary of a rooted tree typically does not influence the
root. This is formalized as point-to-set correlations
\cite{montanari2006rigorous}. In the present case the graphical model is
fully connected and hence a derivation based on belief propagation
might seem puzzling at the first sight. The main assumptions done in
belief propagation is that the incoming messages are probabilistically
independent when conditioned on the value of the root. For this to be
true the underlying graphical model does not have to be tree-like it
can also be because the interactions are sufficiently weak. In the
present case the factor $1/\sqrt{N}$ in the Hamiltonian
(\ref{Boltz_XX}-\ref{Boltz_UV}) makes the interactions sufficiently
weak so that the assumption of independence of incoming messages is
correct in the large size limit. This is quite a common situation
known from the cavity method derivation of the solution for
the Sherrington-Kirkpatrick spin glass model
\cite{mezard1986sk}. The resulting equations can be
derived also by other means, e.g. the Plefka expansion \cite{0305-4470-15-6-035} that we use
later to derive the Bethe free energy. We prefer the derivation from
BP, because that is the only one we know of that provides time indices
that lead to a converging algorithm (note that the TAP equations with
the "naive" time indices were notoriously known not to converge even
in the paramagnetic phase).

Dealing with continuous
variables makes the BP algorithm quite unsuitable to implementation,
as one would need to represent whole probability distributions. One
takes advantage of the fact that BP assumes incoming messages to be
statistically independent and since there are many incoming messages
due to central limit theorem only means and variances of the products
give the desired result called {\it relaxed-belief propagation} (in analogy
to similar algorithm for sparse linear estimation \cite{rangan2010estimation}). 

Next step is analogous to what is done for the TAP equations in the SK
model \cite{ThoulessAnderson77} where one realizes that the algorithm
can be further simplified and instead of dealing with one message for
every directed edge, we can work only with node-dependent
quantities. This gives rise to the so-called Onsager terms. Keeping
track of the correct time indices under iteration
\cite{zdeborova2015statistical} in order to preserve convergence of
the iterative scheme, we end up with the Low-RAMP algorithm that was
derived in various special forms in several previous works. In the
early literature on the Hopfield model these equations were written in
various forms \cite{MezardParisi87b}  without paying close attention
to the time iteration indices that cause convergence problems
\cite{zdeborova2015statistical}. For more recent revival of interest in
this approach that was motivated mainly by
application of similar algorithmic scheme to the problem of linear
estimation \cite{DonohoMaleki09,bayati2011dynamics} and by the need of new algorithms for constrained low-rank estimation see
\cite{rangan2012iterative,NIPS2013_5074,DeshpandeM14,lesieur2015phase,mezard2016mean}. 
We finish this section with simplifications that can be done due to
self-averaging, or due to the Nishimori condition in the Bayes-optimal
case of due to the particular form of the Hamiltonian. 

We note that the Low-RAMP algorithm is related but different from the
variational mean field equations for this problem. To make the
difference apparent, we state the the variational mean field
algorithms in the notation of this paper in Appendix
\ref{Mean_Field_Appendix}.

\subsection{From BP to relaxed BP}

\begin{figure}
\begin{center}
\def\svgwidth{1.1	\linewidth}
\hspace{2cm}
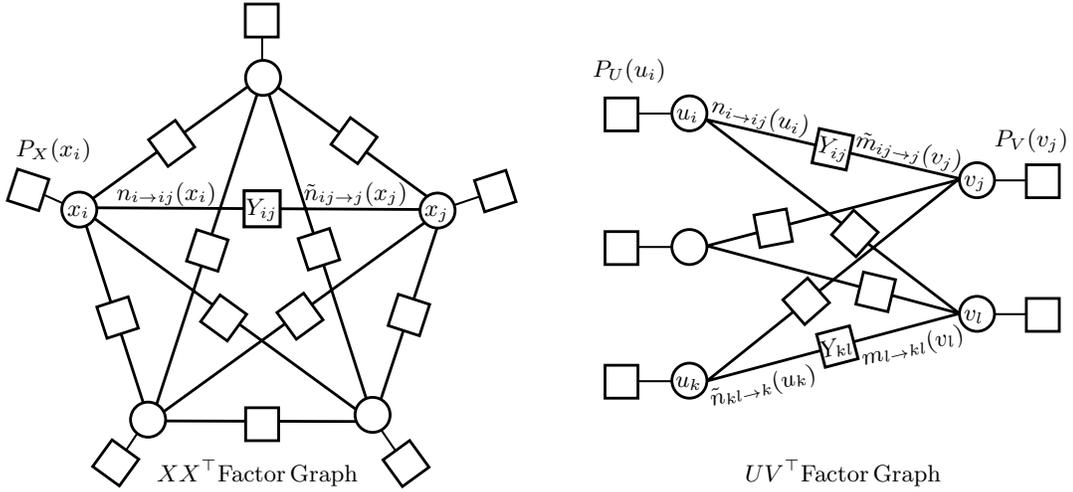
\end{center}
\caption{
This is the factor graph in the symmetric, $XX^\top$, and bipartite,
$UV^\top$, matrix factorization. The squares are factors (or interaction terms), the circles represent variables.
This factor graph allows us to introduce messages $n_{i \rightarrow ij}(x_i)$ and $\hat{n}_{ij \rightarrow j}(x_j)$ for the $XX^\top$ case. These are messages from variables to factors and from factors to variables.
For the $UV^\top$ we introduce the four kinds of messages.   $\hat{n}_{ij \rightarrow i}(u_i)$, $\hat{m}_{ij \rightarrow j}(v_j)$, $\hat{n}_{kl \rightarrow k}(u_k)$ and $m_{l \rightarrow kl}(v_l)$.}
\label{BothFactorGraph_Fig}
\end{figure}

To write the BP equations for the probability measure (\ref{Boltz_XX})
we represent it by a fully connected factor graph, Fig.~\ref{BothFactorGraph_Fig}, where every node
corresponds to a variables node $x_i$, every edge $(ij)$ corresponds to a pair-wise
factor node $e^{g(Y_{ij},x^\top_i x_j /\sqrt{N})}$, and every node is related to a
single site factor node $P_X(x_i)$.  

We introduce the messages $n_{i
  \rightarrow ij}(x_i)$, $\tilde n_{ij \rightarrow i}(x_i)$ respectively as
the $r$-dimensional messages from a variable node to a factor node and from a factor node
to a variable node. The belief propagation equations then are
\bea
\tilde{n}_{ki \rightarrow i}(x_i) &=& \frac{1}{Z_{ki \rightarrow i}}
\int {\rm d}x_k n_{k \rightarrow ki}(x_k) e^{
  g\left(Y_{ki},\frac{x_k^{\top} x_i}{\sqrt{N}} \right) }\, ,
\label{BP_Equation_XX_1}
\\
n_{i \rightarrow ij}(x_i) &=& \frac{P_X(x_i)}{Z_{i \rightarrow ij}} \prod	\limits_{1 \leq k \leq N , k \neq i,j} \tilde{n}_{ki \rightarrow i}(x_i) \,.
\label{BP_Equation_XX_2}
\eea
The factor graph from which these messages are derived is given in Fig.~\ref{BothFactorGraph_Fig}.
The most important assumption made in the BP equations is that the
messages $\tilde n_{ki \to i}(x_i)$ are conditionally on the value $x_i$
independent of each other thus allowing to write the product in
eq.~(\ref{BP_Equation_XX_2}). 

The message in (\ref{BP_Equation_XX_1}) can be expanded as in (\ref{ExpandP_out})  around $w=0$ thanks to the $1/\sqrt{N}$ term. One gets
\be
\tilde n_{ki \rightarrow i}(x_i) = \frac{e^{g(Y_{ki},0)}}{Z_{ki \rightarrow i}} \int {\rm d}x_k n_{k \rightarrow ki}(x_k) \left[ 1 + S_{ki}\frac{x_k^{\top} x_i}{\sqrt{N}}+
\frac{x_i^\top x_k x_k^\top x_i}{2N}
R_{ki} + O\left(\frac{1}{N^{3/2}}\right)\right] \,,
\label{SE_Equation_XX_Expanded}
\ee
where matrices $S_{ij}$ and $R_{ij}$ were defined in
(\ref{Define_S}-\ref{Define_R}). 
One then defines the mean ($r$-dimensional vector) and covariances
matrix (of size $r\times r$) of the message $n_{k
  \rightarrow ki}$ as 
\bea
   {\hat{x}}_{k \rightarrow ki}&=&   \int {\rm d}x_k  n_{k \rightarrow
     ki}(x_k) x_k^\top\, , \\
     \sigma_{x,k\rightarrow ki}     &=&    \int {\rm d}x_k  n_{k
       \rightarrow ki}(x_k) x_k x_k^\top  -  {\hat{x}}_{k \rightarrow
       ki}  {\hat{x}}_{k \rightarrow ki}^\top \, .
\eea
The mean with respect to $n_{k \rightarrow
  ki}(x_k)$  is then taken in (\ref{SE_Equation_XX_Expanded})    and one gets
\begin{multline}
\tilde{n}_{ki \rightarrow i}(x_i) = \frac{1}{Z_{ki \rightarrow i}} \exp\left[g(Y_{ki},0)
+ S_{ki}\frac{{\hat{x}}_{k\to ki}^{\top} x_i}{\sqrt{N}}
-
\frac{x_i^\top {\hat{x}}_{k \rightarrow ki}{\hat{x}}_{k \rightarrow ki}^\top x_i}{2N}S_{ki}^2
\right.
+\\
\left.
\frac{x_i^\top ({\hat{x}}_{k\rightarrow ik} {\hat{x}}_{k\rightarrow ik}^\top + \sigma_{x,k \rightarrow ik}) x_i}{2N}
R_{ki} + O\left(\frac{1}{N^{3/2}}\right) \right] \,.
\label{SE_Equation_XX_Expanded_2}
\end{multline}
Eqs. (\ref{BP_Equation_XX_2}) and (\ref{SE_Equation_XX_Expanded_2}) are combined to get
\begin{equation}
n_{i \rightarrow ij}(x_i) = \frac{P_X(x_i)}{Z_{i \rightarrow ij}} \exp\left( B_{X,i \rightarrow ij}^\top x_i - \frac{x_i^\top A_{X,i \rightarrow ij} x_i}{2}  \right) 
\label{AMP_EQUATION_I-IJ} \,,
\end{equation}
where the $r$-dimensional vector $B_{i \rightarrow ij}$ and the
$r\times r$ matrix $A_{i \rightarrow ij}$ are defined as
\bea
B_{X,i \rightarrow ij} &=& \frac{1}{\sqrt{N}}\sum\limits_{1 \leq k \leq
  N,k \neq  j} S_{ki} {\hat{x}}_{k \rightarrow ki} \, ,
\\
A_{X,i \rightarrow ij} &=& \frac{1}{N}\sum\limits_{1 \leq k \leq N,k
  \neq  j} \left[ S_{ki}^2 {\hat{x}}_{k \rightarrow ki} {\hat{x}}_{k \rightarrow ki}^\top - R_{ki} \left( {\hat{x}}_{k \rightarrow ki}{\hat{x}}_{k \rightarrow ki}^\top + \sigma_{x,k \rightarrow ki}\right)\right] \,.
\eea

The new mean and variance of the message (\ref{AMP_EQUATION_I-IJ})  then needs to be computed. For this we define the function $f^x_{\rm in}$
\begin{equation}
f^x_{\rm in}(A,B) \equiv \frac{\partial}{\partial B} \log\left(\int {\rm d}x P_{X}(x) \exp \left( B^\top x - \frac{x^\top A x}{2} \right) \right)
=
\frac{\int {\rm d}x P_X(x) \exp\left(B^\top x - \frac{x^\top A
      x}{2}\right){x}	}{\int {\rm d}x P_X(x) \exp\left(B^\top x -
    \frac{x^\top A x}{2}\right)}   \label{def_fin}
\end{equation}
as the mean of the normalized density probability
\begin{equation}
{\cal W}(x,A,B)= \frac{1}{Z_x(A,B)} P_X(x) \exp \left(B^\top x - \frac{x^\top A x}{2} \right)
\label{Define_F_in}\,.
\end{equation}
The variance of the message (\ref{AMP_EQUATION_I-IJ})  can be computed
by writing the derivative of $f^x_{\rm in}(A,B)$ with respect to $B$
and getting 
\begin{eqnarray}
\frac{\partial f^x_{\rm in}(A,B)}{\partial B} =
\mathbb{E}_{\cal W(A,B)}(x x^\top) - f^x_{\rm in}(A,B){f^x_{\rm in}}^\top(A,B) \,.
\end{eqnarray}
This expression is the covariance matrix of distribution
(\ref{Define_F_in}). 
Also it is worth nothing that
\begin{equation}
\frac{ \partial \log (Z_x(A,B))}{\partial B} =
f_{\rm in}^x(A,B) \,.
\end{equation}
Adding the time indexes to clarify how these equations are iterated we
get the following {\it relaxed BP} algorithm for the symmetric
low-rank matrix estimation
\bea
B_{X,i \rightarrow ij}^t &=& \frac{1}{\sqrt{N}}\sum\limits_{1 \leq k
  \leq N,k \neq  j} S_{ki}  {\hat{x}}_{k \rightarrow ki}^t \, ,\label{message_B} \\
A_{X,i \rightarrow ij}^t &=& \frac{1}{N}\sum\limits_{1 \leq k \leq N,k
  \neq  j} \left[ S_{ki}^2  {\hat{x}}_{k \rightarrow
    ki}^t{{\hat{x}}_{k \rightarrow ki}^{t,\top}}  - R_{ki} 
  ({\hat{x}}_{k \rightarrow ki}^t{{\hat{x}}_{k \rightarrow
      ki}^{t,\top}} + {\sigma}_{x,k \rightarrow ki}^t) \right]  \, ,\label{message_A}
\\
{\hat{x}}_{i \rightarrow ij}^{t+1} &=& f^x_{\rm in}(A_{i \rightarrow
  ij}^t,B_{i \rightarrow ij}^t) \, ,
\label{Mean_BP_XX}
\\
\sigma_{x,i \rightarrow ij}^{t+1} &=& \frac{\partial f^x_{\rm in}}{\partial B}(A_{i \rightarrow ij}^t,B_{i \rightarrow ij}^t) \,.
\label{Variance_BP_XX}
\eea

\subsection{Low-RAMP: TAPyfication and Onsager terms}
\label{TAP_Equation_Section_XX}
The above relaxed BP algorithm uses $O(N^2)$ messages which can be
memory demanding. But all the messages depend only weakly on the
target node, and hence the algorithm can be reformulated using only
$O(N)$ messages and collecting correcting terms called the Onsager
terms in order to get estimators of the marginals that in the large
size limit are equivalent to the previous ones. We call this
formulation the TAPyfication, because of its analogy to the work of
\cite{ThoulessAnderson77}.
We present the derivation in the case of symmetric low-rank matrix estimation.
We notice that the variables $B_{i \rightarrow ij}$ and $A_{i
  \rightarrow ij}$ depend only weakly on the target node $j$. One can use this fact to close the equations on the marginals of the system.
In order to do that we introduce the variables $B_{X,i}$ and $A_{X,i}$ as
\bea
B_{X,i}^t &=& \frac{1}{\sqrt{N}}\sum\limits_{1 \leq k \leq N} S_{ki}
{\hat{x}}_{k \rightarrow ki}^t\, ,
\label{Equation_B_i}
\\
A_{X,i}^t &=& \frac{1}{N}\sum\limits_{1 \leq k \leq N}\left[ S_{ki}^2  {\hat{x}}_{k \rightarrow ki}^t{{\hat{x}}_{k \rightarrow ki}^{t,\top}}  - R_{ki}  \left( {\hat{x}}_{k \rightarrow ki}^t{{\hat{x}}_{k \rightarrow ki}^{t,\top}} + \sigma_{x,k \rightarrow ki}^t\right)\right]\,. \label{Equation_A_i}
\eea
We also define the variables ${\hat{x}}_i^t$ and $\sigma_{x,i}^t$ as
the estimators of the mean and covariance matrix of  $x_i$, reading
\bea
{\hat{x}}_i^{t+1} &=& f^x_{\rm in}(A_{X,i}^t,B_{X,i}^t) \, ,\label{EquationTAP_UpdateMean}
\\
\sigma_{x,i}^{t+1} &=& \frac{\partial f^x_{\rm in} }{\partial B}(A_{X,i}^t,B_{X,i}^t) \,.
\label{EquationTAP_UpdateVariance}
\eea
In order to close the equations we need to write the $B_{X,i}^t$ and
$A_{X,i}^t$ as a function of the estimators ${\hat{x}}_i^t$ and
$\sigma_{x,i}^t$. 

From the definition of the parameters $A$ and $B$, we have that $\forall
j$, $B_{X,i}^t - B_{X,i \rightarrow ij}^t = \frac{S_{ij}}{\sqrt{N}}
\hat x^t_{j \rightarrow ij} =O\left( \frac{1}{\sqrt{N}} \right)$.
$A_{X,i} - A_{X,i \rightarrow ij} = O\left( \frac{1}{N} \right)$.
One deduces from (\ref{Mean_BP_XX}) and (\ref{Variance_BP_XX}),
(\ref{EquationTAP_UpdateMean}) and (\ref{EquationTAP_UpdateVariance}),
and the Taylor expansion, that in the leading order the difference
between the messages and the node estimators is 
\begin{equation}
{\hat{x}}_{k \rightarrow ki}^t - {\hat{x}}_k^t = f(A_{X,k \rightarrow
  ki}^{t-1},B_{X,k \rightarrow ki}^{t-1})  -
f(A_{X,k}^{t-1},B_{X,k}^{t-1})  =
-\frac{S_{ki}}{\sqrt{N}}\sigma^t_{x,k} {\hat{x}}_{i \to ki}^{t-1} +
O\left(\frac{1}{N}\right) = -\frac{S_{ki}}{\sqrt{N}}\sigma^t_{x,k} {\hat{x}}_{i}^{t-1} + O\left(\frac{1}{N}\right)
\label{Equation_TAP_XX_Mean}\,.
\end{equation}
By plugging (\ref{Equation_TAP_XX_Mean}) into (\ref{Equation_B_i}) one
gets in the leading order
\bea
B_{X,i}^t &=& \frac{1}{\sqrt{N}} \sum\limits_{k = 1}^N S_{ki}
{\hat{x}}_k^t - {\hat{x}}_i^{t-1} \frac{1}{N}\sum\limits_{k = 1}^N
S_{ki}^2 \sigma^t_{x,k} \, , \label{Low-RAMP_B}
\\
A_{X,i}^t &=& \frac{1}{N}\sum\limits_{k = 1}^N \left[ S_{ki}^2
  {\hat{x}}_{k}^t{{\hat{x}}_{k}^{t,\top}}   - R_{ki} \left(
    {\hat{x}}_{k}^t{{\hat{x}}_{k}^{t,\top}} + \sigma_{x,k}^t\right)
\right] \, . \label{Low-RAMP_A}
\eea
These two equations, together with eqs. (\ref{EquationTAP_UpdateMean}-\ref{EquationTAP_UpdateVariance}) give us the low-rank
approximate message passing algorithm (Low-RAMP) with $O
\left(N\right)$ variables.
The second term in the equation (\ref{Low-RAMP_B}) is called the
Onsager reaction term. Notice the iteration index $t-1$ which is
non-intuitive on the first sight and was often misplaced until
recently, see e.g. discussion in \cite{zdeborova2015statistical}. Note
also that there is no Onsager reaction term in
the expression for the covariance $A_{X,i}$, that is because the
individual terms in a sum in $A$ are of order $O(1/N)$ and not
$O(1/\sqrt{N})$. This is a common pattern in AMP-type algorithm, the
Onsager terms appear only in the terms that estimate means, not in the
variances.  

The Low-RAMP algorithm is related in spirit to the AMP algorithm for
linear sparse estimation \cite{DonohoMaleki09}, for instance the function $f_{\rm
  in}$ is the same thresholding function as in the
linear-estimation AMP.  However, the linear estimation AMP is more
involved for a generic output channel and the structure of the two
algorithms are quite different, stemming from the fact that in the
present case all interactions are pairwise whereas for the linear
estimation each interaction involves all the variables, giving rise
to non-trivial terms that do not appear in Low-RAMP. 

The following pseudocode summarizes our implementation of the Low-RAMP algorithm: 

\begin{codebox}
\Procname{$\proc{Low-RAMP symmetric}(S_{ij},H_{ij},r,f_{\rm in}^x,\lambda,\epsilon_{\rm criterium},t_{\rm max},\hat{x}^{\rm init})$}
\li Initialize each $\hat{x}_i \in \mathbb{R}^{r \times 1}$ vector using $\hat{x}^{\rm init}$: $\forall i,\,	\hat{x}_i \gets \hat{x}_i^{\rm init}$.
\li Initialize each $\hat{x}_i^{\rm old} \in \mathbb{R}^r$ vector to zero: $\forall i,\,\hat{x}_i^{\rm old} \gets 0$.
\li Initialize each vector $B_{X,i} \in \mathbb{R}^{r \times 1}$ to zero,
$B_{X,i} \gets 0$.
\li Initialize each N $r \times 1$ vector $B_{X,i}^{\rm old}$ to zero,
$\forall i,\,B_{X,i}^{\rm old} \gets 0$.
\li Initialize to zero each N matrix, $r\times r$ matrix $A_{X,i}$ with, $A_{X,i} \gets 0$.
\li Initialize to zero each N matrix, $r\times r$ matrix $A_{X,i}^{\rm old}$ with, $A_{X,i}^{\rm old} \gets 0$.
\li Initialize to zero each N matrix, $r\times r$ matrix $\sigma_{X,i}$ with, $\sigma_{X,i} \gets 0$.

\li \While ${\rm conv}*\lambda> \epsilon_{\rm criterium}$ and $t< t_{\rm max}$:
\li     \Do $t \gets t+1$;
\li 	$\forall i,\,B_{X,i}^{\rm new} \gets $ Update with equation (\ref{Low-RAMP_B}) or (\ref{Low-RAMP_Bfully}).
\li 	$\forall i,\,A_{X,i}^{\rm new} \gets $ Update with equation (\ref{Low-RAMP_A}) or (\ref{Low-RAMP_Afully}).
	
\li $\forall i,\:\:B_{X,i} \gets \lambda B_{X,i}^{\rm new} + (1-\lambda) B_{X,i}^{\rm old}$,
\li $\forall i,\:\:A_{X,i} \gets \lambda A_{X,i}^{\rm new} + (1-\lambda) A_{X,i}^{\rm old}$,
\li			$\forall i,\:\: {\hat x_i^{\rm old}} \gets \hat{x}_i ,\:	{\hat{x}}_i \gets f_{\rm in}^x(A_X,B_{X,i})$,
\li			$\forall i,\:\:	 \sigma_{X,i} \gets \frac{\partial f_{\rm in}^x}{\partial B}(A_X,B_{X,i})$,

\End
\li 	${\rm conv} \gets \frac{1}{N}\sum\limits \Vert \hat{x}_i - \hat{x}_i^{\rm old} \Vert$.
\End
\li \Return signal components $\bx$.
\end{codebox}

The canonical initialization we use is
\begin{equation}
\forall i \in [1;N],\:\: \hat{x}^{\rm init}_i \gets 10^{-3}{\cal
  N}(0,I_r)\, .
\label{TAP_Equation_init_r}
\end{equation}
The constant $10^{-3}$ here can be changed, but it is a bad idea to
initialize exactly at zero since $\hat{x}_i = 0$ could be exactly a
fixed point of the equations. In order to analyze the algorithm for a
specific problem it is instrumental to initialize in the
solution:
\begin{equation}
\forall i \in [1;N],\:\: \hat{x}^{\rm init}_i \gets {x}^{0}_{i}\, ,
\label{TAP_Equation_init_pl}
\end{equation}
where  $ {x}^{0}_{i}$ is the planted (ground truth) configuration. 

In the above pseudocode the damping factor $\lambda$ is chosen
constant for the whole duration of the algorithm. It is possible to choose $\lambda$ dynamically  in  order to improve
the convergence. Using the fact that the Low-RAMP algorithm finds a stationary fixed point of the Bethe free energy given in (\ref{PlefkaSection}) one can choose the damping factor $\lambda$ so that at each step so that the Bethe-free-energy increases, this is described in \cite{DBLP:journals/corr/VilaSRKZ14}.	
Another way to choose $\lambda$ is by ensuring that $\sum \Vert \hat{x}^t - \hat{x}^{t+1} \Vert$ does not oscillate too much. If one sees too much oscillations one increases the damping and decreases it otherwise.
Further way to improve convergence is related to
randomization of the update scheme as argued for the related compressed
sensing problem in \cite{caltagirone2014convergence}.

\subsection{Simplifications}

The above generic Low-RAMP equations simplify further under certain
conditions. 

\subsubsection{Advantage of self-averaging}
\label{sec:self_av}

We can further simplify these equations by noticing that in all the
expressions where $S_{ij}^2$ appears we can replace it by its mean
without changing the leading order of the quantities. This
follows from the assumption made in the BP equations, that states that
the messages incoming to a node are independent conditionally on the
value of the node. Consequently the sums in
eqs.~(\ref{message_B}-\ref{message_A}) are sum of $O(N)$ independent
variables and can hence in the leading order be replaced by their
means.

This allows us to write the Low-RAMP equations (\ref{Low-RAMP_B}-\ref{Low-RAMP_A}) in an even simpler
form
\bea
B_{X,i}^t &=& \frac{1}{\sqrt{N}} \sum\limits_{k = 1}^N S_{ki}
{\hat{x}}_k^t - \left( \frac{1}{N\widetilde \Delta}\sum\limits_{k =
    1}^N\sigma^t_{x,k}\right) {\hat{x}}_i^{t-1}\, ,\label{Low-RAMP_Bfully}
\\
A_{X}^t &=& \frac{1}{N  \widetilde \Delta}\sum\limits_{k = 1}^N
{\hat{x}}_{k}^t{{\hat{x}}_{k}^{t,\top}}   - \overline R  \frac{1}{N}
\sum\limits_{k = 1}^N \left( {\hat{x}}_{k}^t{{\hat{x}}_{k}^{t,\top}} +
  \sigma_{x,k}^t\right) \, , \label{Low-RAMP_Afully}
\eea
where we defined	
\begin{eqnarray}
\frac{1}{\widetilde \Delta} &\equiv& \frac{2}{N^2} \sum_{1 \le i<j \le
  N}   S_{ij}^2 \, , \label{def_tilde_Delta}
\\
\overline R &\equiv&  \frac{2}{N^2} \sum_{1 \le i<j \le N}   R_{ij}
\, . \label{def_overline_R}
\end{eqnarray}
Whereas $\widetilde \Delta$ is always positive, $\overline R $ can be positive or negative. 
Together with
eqs.~(\ref{EquationTAP_UpdateMean}-\ref{EquationTAP_UpdateVariance})
the above two expressions are a closed set of equations. Note in
particular that in (\ref{Low-RAMP_Afully}) the dependence on index $i$
disappeared in the leading order.

\subsubsection{Bayes-optimal case}
\label{Nish_TAP}

In the Bayes-optimal inference case we derived expression
(\ref{E_2nd_der}). Putting it together with the definition of the matrix
$R_{ij}$ in eq.~(\ref{Define_R}) and realizing that the average over
sites $i,j$ and the average over $P(Y|w)$ act the same way we get that for $\overline R$ as defined in
(\ref{def_overline_R}) we have $\overline R=0$. This property belongs to
the class of properties called the Nishimori conditions. In the Bayes
optimal case the expression for $A_{X}^t$ simplified further into
\be
A_{X}^t = \frac{1}{N\Delta}\sum\limits_{k = 1}^N
{\hat{x}}_{k}^t{{\hat{x}}_{k}^{t,\top}}   \, ,  \label{define_A_Nish}
\ee
where we used the Bayes-optimality once more to realize that $\widetilde \Delta
= \Delta$ as defined in eq.~(\ref{Define_Delta}). The convenient
property of the Bayes-optimality is that the quantity $A_{X}^t$ now
has to be non-negative. 

\subsubsection{Conventional Hamiltonian : SK model}

Another case that is worth specifying is the conventional Hamiltonian
where $g$ is given by \eqref{can_Ham_Example}.
One gets
\be
A_{X}^t =  - \overline R  \frac{1}{N}
\sum\limits_{k = 1}^N \sigma_{x,k}^t \, .
\ee
It is slightly counter-intuitive that this variance-like term is
negative, but it is always used only in the function $f_{\rm in}^x$
defined in (\ref{def_fin}) where it gets multiplied by the $P_X(x_i)$
therefore if $P_X$ is decaying fast enough or has a bounded support,
the corresponding integrals exist and are finite.

This is a convenient point where we can make the link between the
Low-RAMP algorithm and the TAP equations for the SK model. For the
Ising spins (\ref{Ising}) the function $f^x_{\rm in}$ becomes 
\bea
      f_{\rm in}^x(A,B) = \tanh{(\beta h + B)} \, ,  \label{fin_SK}\quad
      \frac{\partial  f_{\rm in}^x(A,B)}{\partial B} = 1 -
      \tanh^2{(\beta h + B)} \, .  \label{fin_SKb}
\eea
Notice the independence on the parameter $A$. 
The conventional Hamiltonian of the SK model corresponds to
$g(Y,w)=\beta Y w$ so that $S=\beta Y$. Which gives us for the update
of the parameter $B$ eq.~(\ref{Low-RAMP_B})
\be
           B_{X,i}^t = \frac{\beta}{\sqrt{N}} \sum\limits_{k = 1}^N Y_{ki}
{\hat{x}}_k^t - {\hat{x}}_i^{t-1} \frac{\beta^2}{N}\sum\limits_{k = 1}^N
Y_{ki}^2 (1-{\hat{x}}_k^{2,t}  ) \, .
\ee
Together with (\ref{fin_SK}) we get the well known TAP equations
\cite{ThoulessAnderson77} 
\be
    {\hat{x}}_i^{t+1} = \tanh{\left[\beta h + \frac{\beta}{\sqrt{N}} \sum\limits_{k = 1}^N Y_{ki}
{\hat{x}}_k^t - {\hat{x}}_i^{t-1} \frac{\beta^2}{N}\sum\limits_{k = 1}^N
Y_{ki}^2 (1-{\hat{x}}_k^{2,t}  ) \right]} \, .
\ee

\subsection{Summary of Low-RAMP for the bipartite low-rank estimation}
\label{TAP_Equation_Section_UV}

The derivation for the bipartite case $UV^\top$  is completely analogous. 
The {\it relaxed BP} equations read 
\bea
B_{U,i \rightarrow ij}^t &=& \frac{1}{\sqrt{N}}\sum\limits_{1 \leq l
  \leq M,l \neq  j} S_{il}  {\hat{v}}_{l \rightarrow il}^t\, ,
\\
A_{U,i \rightarrow ij}^t &=& \frac{1}{N}\sum\limits_{1 \leq l \leq M,l
  \neq  j} \left[ S_{il}^2 {\hat{v}}_{l \rightarrow il}^t{{\hat{v}}_{l
      \rightarrow il}^{t,\top}} - R_{il}\left( {\hat{v}}_{l
      \rightarrow il}^t{{\hat{v}}_{l \rightarrow il}^{t,\top}} +
    \sigma_{v,l \rightarrow il}^t\right)\right]\, ,
\\
{\hat{u}}_{i \rightarrow ij}^{t} &=& f^u_{\rm in}(A_{U,i \rightarrow
  ij}^t,B_{U,i \rightarrow ij}^t)\, ,
\label{Mean_TAP_U}
\\
\sigma_{u,i \rightarrow ij}^{t} &=& \frac{\partial f^u_{\rm
    in}}{\partial B}(A_{U,i \rightarrow ij}^t,B_{U,i \rightarrow
  ij}^t)\, ,
\label{Variance_TAP_U}
\\
B_{V,j \rightarrow ij}^t &=& \frac{1}{\sqrt{N}}\sum\limits_{1 \leq k
  \leq N,k \neq  j} S_{kj}  {\hat{u}}_{k \rightarrow kj}^t\, ,
\\
A_{V,j \rightarrow ij}^t &=& \frac{1}{N}\sum\limits_{1 \leq k \leq N,k
  \neq i} \left[ S_{kj}^2 {\hat{u}}_{k \rightarrow kj}^t{{\hat{u}}_{k
      \rightarrow kj}^{t,\top}} - R_{kj}\left( {\hat{u}}_{k
      \rightarrow kj}^t{{\hat{u}}_{k \rightarrow kj}^{t,\top}} +
    \sigma_{u,k \rightarrow kj}^t\right)\right]\, ,
\\
{\hat{v}}_{j \rightarrow ij}^{t+1} &=& f^v_{\rm in}(A_{V,j \rightarrow
  ij}^t,B_{V,j \rightarrow ij}^t)\, ,
\label{Mean_TAP_V}
\\
\sigma_{v,j \rightarrow ij}^{t+1} &=& \frac{\partial f^v_{\rm in}}{\partial B}(A_{V,j \rightarrow ij}^t,B_{V,j \rightarrow ij}^t)
\label{Variance_TAP_V} \,.
\eea
Note that here we broke the symmetry between $U$ and $V$ by choosing
an order in the update. We first update estimators of $U$ without
increasing the time index and only then estimators of $V$ while
increasing the time index by one. 

The Low-RAMP equations with their Onsager terms for  the bipartite low-rank matrix
estimation read
\bea
B^t_{U,i}&=& \frac{1}{\sqrt{N}}
 \sum\limits_{l=1}^M S_{il} \hat{v}^t_{l} - 
 \left(\frac{1}{N}
 \sum\limits_{l=1}^M S_{il}^2\sigma_{v,l}^{t} \right)
\hat{u}^{t-1}_{i} \, ,\label{TAP_UV_B_U}
\\
A_{U,i}^t &=& \frac{1}{N}\sum\limits_{l=1}^M \left[ S_{il}^2
  \hat{v}_{l}^t{\hat{v}_{l}^{t,\top}}  - R_{il} \left(
    {\hat{v}}_{l}^t{{\hat{v}}_{l}^{t,\top}} +
    \sigma_{v,l}^t\right)\right]\, ,	 \label{TAP_UV_A_U}
\\
\hat{u}^t_{i} &=& f^u_{\rm in}(A^t_{U,i},B^t_{U,i})\, , 
\\
\sigma^t_{u,i} &=& \left(\frac{\partial f^u_{\rm in}}{ \partial B}
\right)(A^t_{U,i},B^t_{U,i})\, ,
\\
B^t_{V,j} &=& \frac{1}{\sqrt{N}} \sum\limits_{k=1}^N S_{kj} \hat{u}_k^t
- 
\left(\frac{1}{N}
 \sum\limits_{k=1}^N S_{kj}^2\sigma_{u,k}^t\right) \hat{v}^{t}_{j}\, , \label{TAP_UV_B_V}
\\
A_{V,j}^t &=& \frac{1}{N}\sum\limits_{k = 1}^N \left[ S_{kj}^2
  \hat{u}_{k}^t{\hat{u}_{k}^{t,\top}} - R_{kj} \left(
    {\hat{u}}_{k}^t{{\hat{u}}_{k}^{t,\top}} +
    \sigma_{u,k}^t\right)\right]\, , \label{TAP_UV_A_V}
\\
\hat{v}^{t+1}_{j} &=& f^v_{\rm in}(A^t_{V,j},B^t_{V,j})\, ,
\\
\sigma^{t+1}_{v,j} &=& \left(\frac{\partial f^v_{\rm in}}{ \partial B} \right)(A^t_{V,j},B^t_{V,j}) \,. \label{TAP_Sigma_V}
\eea

Due to independence between messages assumed in the BP algorithm we
can simplify the Low-RAMP algorithm introducing
(\ref{def_tilde_Delta}) and (\ref{def_overline_R}) also for the
bipartite case. As well as take advantage of the Bayes optimality and
the Nishimori conditions or of the particular form of the conventional
Hamiltonian. 

\begin{codebox}
\Procname{$\proc{LowRAMP-bipartite}(S_{ij},H_{ij},\Delta,r,f_{\rm in}^u,f_{\rm in}^v,\lambda,\epsilon_{\rm criterium},t_{\rm max}, \hat{u}^{\rm init},\hat{v}^{\rm init})$}
\li Initialize each $N$, $r \times 1$ vector $\hat{u}_i$ using $\hat{u}^{\rm init}$, $\forall i,\,	\hat{u}_i \gets \hat{u}_i^{\rm init}$.
\li Initialize each $M$, $r \times 1$ vector $\hat{v}_j$ using $\hat{v}^{\rm init}$, $\forall j,\,	\hat{v}_j \gets \hat{v}_j^{\rm init}$.
\li Initialize each $M$, $r \times 1$ vector $\hat{u}_i^{\rm old}$ to 0, $\forall i,\,	\hat{u}_i^{\rm old} \gets 0$.
\li Initialize each $M$, $r \times 1$ vector $\hat{v}_j^{\rm old}$ to 0, $\forall j,\,	\hat{v}_j^{\rm old} \gets 0$.
\li Initialize to zero each N, $r \times 1$ vector $B_{U,i}$ to zero,
$\forall i,\,B_{U,i} \gets 0$.
\li Initialize to zero each M, $r \times 1$ vector $B_{V,j}$ to zero,
$\forall j,\,B_{V,j} \gets 0$.
\li Initialize to zero each N matrix, $r\times r$ matrix $A_{U,i}$ with, $\forall j,\,A_{U,i} \gets 0$.
\li Initialize to zero each N matrix, $r\times r$ matrix $A_{U,j}^{\rm old}$ with, $\forall i,\,A_{U,i}^{\rm old} \gets 0$.
\li Initialize to zero each M matrix, $r\times r$ matrix $A_{V,i}$ with, $\forall j,\,A_{V,j} \gets 0$.
\li Initialize to zero each N matrix, $r\times r$ matrix $A_{V,j}^{\rm old}$ with, $\forall j,\,A_{V,j}^{\rm old} \gets 0$.
\li Initialize to zero each N matrix, $r\times r$ matrix $\sigma_{U,i}$ with, $\sigma_{U,i} \gets 0$.
\li Initialize to zero each M matrix, $r\times r$ matrix $\sigma_{V,j}$ with, $\sigma_{V,j} \gets 0$.

\li \While ${\rm conv}*\lambda> \epsilon_{\rm criterium}$ and $t< t_{\rm max}$:
\li     \Do $t \gets t+1$;
\li {\bf Update variables $U$}
\li 	$\forall i,\:\:B_{U,i}^{\rm new} \gets $ Update with equation (\ref{TAP_UV_B_U}).
\li 	$\forall i,\:\:A_{U,i}^{\rm new} \gets $ Update with equation (\ref{TAP_UV_A_U}).

\li $\forall i,\:\:B_{U,i} \gets \lambda B_{U,i}^{\rm new} + (1-\lambda) B_{U,i}^{\rm old}$,
\li $\forall i,\:\:A_{U,i} \gets \lambda A_{U,i}^{\rm new} + (1-\lambda) A_{U,i}^{\rm old}$,

\li			$\forall i,\:\: {\hat u_i^{\rm old}} \gets \hat{u}_i ,\:	{\hat{u}}_i \gets f_{\rm in}^u(A_ {U,i},B_{U,i})$,
\li			$\forall i,\:\: \sigma_{U,i} \gets \frac{\partial f_{\rm in}^u}{\partial B}(A_{U,i},B_{U,i})$.
\li {\bf Update variables $V$}
\li 	$\forall j,\:\:sB_{V,j}^{\rm new} \gets $ Update with equation (\ref{TAP_UV_B_V}).
\li 	$\forall j,\:\:A_{V,j}^{\rm new} \gets $ Update with equation (\ref{TAP_UV_A_V}).

\li $\forall j,\:\:B_{V,j} \gets \lambda B_{V,j}^{\rm new} + (1-\lambda) B_{V,j}^{\rm old}$,
\li $\forall j,\,A_{V,j} \gets \lambda A_{V,j}^{\rm new} + (1-\lambda) A_{V,j}^{\rm old}$,

\li			$\forall j,\:\: {\hat v_j^{\rm old}} \gets \hat{u}_i ,\:	{\hat{u}}_i \gets f_{\rm in}^v(A_ {V,j},B_{V,j})$,
\li			$\forall j,\:\: \sigma_{V,j} \gets \frac{\partial f_{\rm in}^v}{\partial B}(A_{V,j},B_{V,j})$.
\li {\bf Compute distance with previous iteration}
\li 	${\rm conv} \gets \frac{1}{N}\sum\limits \Vert \hat{u}_i - \hat{u}_i^{\rm old} \Vert 
+  \frac{1}{M}\sum\limits \Vert \hat{v}_j - \hat{v}_j^{\rm old} \Vert$.
\End
\End
\li \Return signal components $\bx$
\end{codebox}

The initialization and damping factor $\lambda$ are chosen similarly
as for the symmetric case as discussed in section \ref{TAP_Equation_Section_XX}. 

\subsubsection{Example The Hopfield model}

In order to illustrate the form of the Low-RAMP equations on an
example of a bipartite matrix factorization problem we choose the
Hopfield model \cite{hopfield1982neural}, a well known model of
associative memory. It turns out that the Hopfield model can be mapped
to the bipartite low-rank matrix factorisation problem.
This is analyzed with greater care recently in \cite{mezard2016mean}.
The Hopfield model is defined by the following Boltzmann distribution
\begin{eqnarray}
P_{\rm Hopfield}(U) &\sim& \exp \left( \beta \sum\limits_{1 \leq i
                           \leq j \leq N} J_{ij} u_i u_j   \right) \, ,\label{hop_Boltz}
\\
J_{ij} &=& \frac{1}{N}\sum\limits_{1 \leq \mu \leq M} Y_{i \mu} Y_{j
           \mu} \, , \label{Hebbs}
\end{eqnarray}
where the $u_i$ are $\pm 1$ spin variables, and $Y_{i \mu}$ is $N 
\times M$ matrix collecting the $M$ of $N$-dimensional patterns to
remember. The elements of the matrix $Y_{i \mu}$ are most commonly
considered as independent $\pm 1$ variables. The
relation (\ref{Hebbs}) is then known as the Hebbs rule of associative
memory. 

Using a Gaussian integral we can map the Hopfield model to the
following probability distribution 
\begin{eqnarray}
P(U,V|Y) = \frac{1}{Z_{UV}(Y)}\exp\left(\sum\limits_{1 \leq \mu \leq M} \left[ - \frac{\beta v_\mu^2}{2} + v_\mu \sum\limits_{1
  \leq j \leq N} \frac{\beta Y_{j \mu} u_j}{\sqrt{N}}\right] \right)
  \prod\limits_{1 \leq i \leq N}P_U(u_i) \, ,  \label{Hopfield_mapped}
\end{eqnarray}
where $P_U(u) = [\delta(u-1) + \delta(u+1) ] /2$. Indeed, we can check
that integrating expression (\ref{Hopfield_mapped}) we recover the
Hopfield Boltzmann distribution (\ref{hop_Boltz}). The Hopfield model
is hence equivalent to the bipartite matrix factorization with
Gaussian spins on one side, Ising spins on the other side, the
conventional Hamiltonian and randomly quenched disorder $Y_{\mu j}$. 

In order to write the Low-RAMP equations
(\ref{TAP_UV_B_U}-\ref{TAP_Sigma_V}), known as the TAP equations in
the context of the Hopfield model, we remind that the $S$ and $R$
matrices for the conventional Hamiltonian are given by
eq. (\ref{can_Ham_Example}). 
Further, we remind that for Ising spins $u$, eq.~(\ref{Ising}), with no magnetic
field ($\rho = 1/2$) the thresholding function $f^u_{\rm in}(A,B)$ is given by 
\begin{equation}
    f^u_{\rm in}(A,B) = \tanh(B)\,  ,  \quad \quad  \frac{\partial
      f^u_{\rm in}(A,B) }{\partial B} = 1-\tanh^2{B} \, .   \label{fin_SK_h0}
\end{equation}
It follows for instance that $\sigma_{u,i}^t = 1 -
(\hat{u}_i^{t})^2$. 
Finally for Gaussian spins with zero mean and variance $1/\beta$ the
thresholding function $f^v_{\rm in}(A,B)$ is given by 
\begin{equation}
    f^v_{\rm in}(A,B) = \frac{B}{A+\beta}\, ,  \quad \quad  \frac{\partial
      f^v_{\rm in}(A,B) }{\partial B} = \frac{1}{A+\beta} \, . 
\end{equation}
We insert these expressions into equations
(\ref{TAP_UV_B_U}-\ref{TAP_Sigma_V}) to get
\begin{eqnarray}
\hat v_{\mu}^{t+1}&=&   \frac{  \frac{1}{ \beta \sqrt N}
                      \sum\limits_{i = 1}^{N} { Y_{i \mu} {\hat
                      u^t_{i} }- {{(1-q^t) \hat v_\mu^{t}}}} }{C^t}
                      \, ,\label{Hopfield_U}
\\
\hat u_{i}^{t+1} &=& \tanh \left[{ \frac{\beta}{\sqrt{N}
                     }\sum_{\mu=1}^{M} Y_{i \mu}\hat v_{\mu}^{t+1}  - \frac{\alpha \hat{u}_i^{t}}{C^t}  }
                     \right] \, ,\label{amps} 
	\end{eqnarray}
where
\begin{eqnarray}
q^t &=& \frac{1}{N}\sum\limits_{1 \leq i \leq N} (\hat{u}_i^{t})^2\,,
\\
C^t &=& \frac{1}{\beta} - (1 - q^t)\,.
\end{eqnarray}

It is possible to close the equations on the $v_j^t$ by injecting
\eqref{Hopfield_U} in \eqref{amps}, and then noticing that the term
$\beta \sum_{\mu=1}^{M} Y_{i\mu}\hat v_{\mu}^{t} /\sqrt{N}$ can be replaced by $\tanh^{-1}(\hat{u}_{i}^{t}) + \frac{\alpha \hat{u}_i^{t-1}}{C^{t}}$ \eqref{amps}
\begin{eqnarray}
\hat{u}_{i}^{t+1} &=& \tanh \left[ \frac{1}{C^t} \sum_{j=1}^N
                    {J}_{ij} \hat{u}_j^t  - \frac{\alpha u_i^{t}}{C^t}- \frac{1 - q^t}{C^t} \left(\tanh^{-1}(\hat{u}_{i}^{t}) + \frac{\alpha \hat{u}_i^{t-1}}{C^{t}}\right)\right] \,.
\end{eqnarray}
These are the same equations as the TAP equations for the Hopfield model discussed in
\cite{mezard2016mean} eqs.~(48,49,50), except that here they are expressed in term of the magnetization and not the fields $\tanh^{-1}(u_i^t)$.

\subsection{Bethe Free Energy}
\label{PlefkaSection}

We define the free energy of a given probability measure as the 
logarithm of its normalization (in physics it is usually the negative
logarithm, but in this article we adopt the definition without the minus
sign).  Notably for the symmetric vector-spin
glass model (\ref{Boltz_XX}) we define 
\begin{equation}
       \Phi_{XX^\top}(Y) = \log{ ( Z_X(Y))} - \sum\limits_{1 \leq i < j \leq N} g(Y_{ij},0)\, , 
\end{equation}
where $Z_X(Y)$ is the normalization, i.e. the partition function, in
(\ref{Boltz_XX}). We subtract the constant term on the right hand
side for convenience. 
For the bipartite vector-spin glass model (\ref{Boltz_UV}) analogously 
\begin{equation}
       \Phi_{UV^\top}(Y) = \log{ ( Z_{UV}(Y))} - \sum\limits_{1 \leq i < N, 1 \leq  j \leq M} g(Y_{ij},0) \, .
\end{equation}

We remove the constant term $g(Y_{ij},0)$ so that the Free energy $\Phi_{XX^\top}$ and $\Phi_{UV^\top}$ will be $O(N)$ and self averaging in the large $N$ limit.
The exact free energies are intractable to compute, therefore we use
approximations equivalent to those used for derivation of the Low-RAMP
algorithm to derive the so-called {\it Bethe free energy}. Under the
assumption of replica symmetry this free energy is exact in the
leading order in $N$ and with high probability over the ensemble of instances. To
derive the Bethe free energy we use the Plefka expansion, later extended
by Georges and Yedidia
\cite{0305-4470-15-6-035,0305-4470-24-9-024}. 
The derivation is
presented in the appendix \ref{app:Plefka}.


The Bethe free energy for the symmetric vector-spin glass case,  $XX^\top$, is
\begin{multline}
\Phi_{{\rm Bethe},XX^\top}  =
\max_{\{A_{X,i}\},\{B_{X,i}\}}\Phi_{{\rm Bethe},XX^\top}(\{A_{X,i}\},\{B_{X,i}\})  \\
\Phi_{{\rm Bethe},XX^\top}(\{A_{X,i}\},\{B_{X,i}\}) = \sum\limits_{1 \leq i \leq N} \log({Z}_x(A_{X,i},B_{X,i})) - B_{X,i}^\top {\hat{x}}_i
+ \frac{1}{2}{\rm Tr}\left[A_{X,i}	({\hat{x}}_i {\hat{x}}_i^\top + \sigma_{x,i}) \right]
\\
+\frac{1}{2}\sum\limits_{1 \leq i ,j\leq N} \left[\frac{1}{\sqrt{N}}S_{ij} {\hat{x}}_i^\top {\hat{x}}_j + \frac{R_{ij}}{2N} {\rm Tr}\left[({\hat{x}}_i {\hat{x}}_i^\top + \sigma_{x,i}) ({\hat{x}}_j {\hat{x}}_j^\top + \sigma_{x,j}) \right]-\frac{S_{ij}^2}{2N}  {\rm Tr} \left[ {\hat{x}}_i {\hat{x}}_i^\top   {\hat{x}}_j {\hat{x}}_j^\top  \right]
-\frac{1}{N}S_{ij}^2 {\rm Tr} \left[\sigma_{x,i}   \sigma_{x,j}  \right] \right] \,,
\label{TAP_XX_FreeEnergy}
\end{multline}
where $\hat x_i=f_{\rm in}^x(A_{X,i}, B_{X,i})$ and $\sigma_{x,i}=\partial_B f_{\rm in}^x
(A_{X,i}, B_{X,i}) $ are considered as explicit functions of $A_{X,i}$
and $B_{X,i}$, where the function $f_{
\rm in}^x$ depends on the prior probability distribution $P_X$ via
eq.~(\ref{def_fin}).

For the bipartite vector-spin glass case, $U V^\top$, the Bethe free energy is
\begin{multline}
\Phi_{{\rm  Bethe},UV^\top} =
\max_{  \{A_{U,i}\},\{B_{U,i}\},\{A_{V,j}\},\{B_{V,j}\} }\Phi_{{\rm
    Bethe},UV^\top}(\{A_{U,i}\},\{B_{U,i}\},\{A_{V,j}\},\{B_{V,j}\})\,
,\\
\Phi_{{\rm Bethe},UV^\top}(\{A_{U,i}\},\{B_{U,i}\},\{A_{V,j}\},\{B_{V,j}\}) = \sum\limits_{1 \leq i \leq N} \log({Z}_u(A_{U,i},B_{U,i}))  -  {B_{U,i}}^\top {\hat{u}}_i
+ \frac{1}{2}{\rm Tr}\left[A_{U,i}({\hat{u}}_i {\hat{u}}_i^\top + \sigma_{u,i}) \right]
\\
+
\sum\limits_{1 \leq j \leq M} \log({Z}_v(A_{V,j},B_{V,j})) -  {B_{V,j}}^\top {\hat{v}}_j
+ \frac{1}{2}{\rm Tr}\left[A_{V,j}({\hat{v}}_j {\hat{v}}_j^\top + \sigma_{v,j}) \right]
\\
+\sum\limits_{1 \leq i \leq N,1 \leq j\leq M} \left[ \frac{1}{\sqrt{N}}S_{ij} \hat{u}	_i^\top \hat{v}_j + \frac{1}{2N}R_{ij} {\rm Tr}\left[(\hat{u}_i \hat{u}_i^\top + \sigma_{u,i})(\hat{v}_j \hat{v}_j^\top + \sigma_{v,j})\right]- \frac{S_{ij}^2{\rm Tr}\left( \hat{u}_i \hat{u}_i^\top \hat{v}_j \hat{v}_j^\top\right)}{2N}
-\frac{1}{N}S_{ij}^2 {\rm Tr} \left[\sigma_{u,i}  \sigma_{v,j} \right]\right] \,, \label{TAP_UV_FreeEnergy}
\end{multline}
where the
$\hat{u}_i=f_{\rm in}^u(A_{U,i}, B_{U,i})$, $\hat{v}_j=f_{\rm
  in}^v(A_{V,j}, B_{V,j})$, $\sigma_{u,i}= \partial_B f_{\rm in}^u(A_{U,i},
B_{U,i})$, $\sigma_{v,j}=\partial_B f_{\rm in}^v(A_{V,j}, B_{V,j})$,
again seen as a function of variables $A$, and $B$. Note that  fixed points of the Low-RAMP
algorithm are stationary points of the Bethe free energy as can be
checked explicitly by taking the derivatives of the formulas.

The main usage of the free energy is when there exist multiple fixed
points of the Low-RAMP equations then the one that corresponds to the best achievable mean squared error is the one for which the free
energy is the largest. Another way to use the free energy is in order
to help the convergence of the Low-RAMP equations, the adaptive damping
is used 
\cite{DBLP:journals/corr/abs-1301-6295,DBLP:journals/corr/VilaSRKZ14} and
relies on the knowledge of the above expression for the Bethe free energy.

%

To conclude, we recall that Low-RAMP is distributed in Matlab and
Julia at \url{http://krzakala.github.io/LowRAMP/},
in a version that include the use of the Bethe free energy as a guide
to increase convergence, as in \cite{DBLP:journals/corr/VilaSRKZ14}.

\section{State Evolution}
\label{StateEvolutionSE}

An appealing property of the Low-RAMP algorithm is that its
large-system-size behaviour can be analyzed via the so called state
evolution (or single letter characterization in information theory). 
In the statistical physics context the state evolution is the cavity
method \cite{MezardParisi87b} thanks to which one can derive the replica symmetric solution
from the TAP equations, taking properly into account the distribution
of the disorder (random or quenched). Mathematically, at least in the Bayes optimal
setting, the state evolution for the present systems is a rigorous
statement about the asymptotic behaviour of the Low-RAMP algorithm 
\cite{javanmard2013state}. 

Here we present derivation of the state evolution for the symmetric
matrix factorization and state it for the bipartite case. 
The main idea of state evolution is to describe the current state of
the algorithm using a small number of variables -- called order
parameters in physics. We then compute how the order parameters evolve as the number of iterations increases. 

\subsection{Derivation for the symmetric low-rank estimation}

To derive the state evolution we assume that all updates are done in
parallel with no damping (the state evolution does depend on the
update strategy).
A distinction will be made between $r$ the rank assumed in the
posterior distribution and $r_0$ the true rank of the planted solution.
Let us introduce the order parameters that will be of relevance here
\begin{eqnarray}
M_x^t &=& \frac{1}{N}\sum\limits_{1 \leq i \leq N} \hat{x}_i^t {x_i^{0,
\top}} \in \mathbb{R}^{r \times r_0},\quad \label{SE_m_def}\\
Q_x^t &=& \frac{1}{N}\sum\limits_{1 \leq i \leq N} \hat{x}_i^t \hat{x}_i^{t,\top} \in \mathbb{R}^{r \times r}
,\quad \label{SE_q_def}\\
\Sigma_x^t &=& \frac{1}{N}\sum\limits_{1 \leq i \leq N} \sigma_{x,i}^t \in \mathbb{R}^{r \times r}\label{SE_S_def}
\,,
\end{eqnarray}
where $M_x^t$ is a matrix of size $r\times r_0$, while $Q_x^t$ and
$\Sigma_x^t$ are $r \times r$ matrices. The interpretation of these
order parameters is the following
\begin{itemize}
\item$M_x^t$ measures how much the current estimate of the mean is
  correlated with the planted solution $x_i^0$. Physicist would call this the magnetization of the system.
\item$Q_x^t$ is called the self-overlap.
\item$\Sigma_x^t$ is the mean variance of variables.
\end{itemize}

In this section we will not assume the Bayes optimal setting, and
distinguish between the prior $P_X(x_i)$ and the distribution $P_{X_0}(x_i^0)$ from
which the planted configuration $x_i^0$ was drawn. Similarly, we
will assume $g(Y,w)$ in the posterior distribution, but the data
matrix $Y$ was created from the planted configuration via $P_{\rm
  out}(Y|w)$. In general $P_X \neq P_{X_0}$ and $g(Y,w) \neq \log{
  P_{\rm out}(Y|w)}$. 

We do assume, however, that \eqref{E_1st_der} holds for our choice of
$g(Y,w)$ and $P_{\rm out}(Y|w)$ even when $g(Y,w) \neq
\log{ P_{\rm out}(Y|w)}$. This will indeed hold in all examples
presented in this paper.  Using self-averaging arguments such as in
Sec.~\ref{sec:self_av} the averages over the quenched randomness $P_{\rm
  out}(Y|w)$ and over elements $(ij)$ are
interchangeable. Eq.~(\ref{E_1st_der}) then in practice means that, in
order for the state evolution as derived in this section to be valid,
the Fisher score matrix $S$ should have an empirical mean of elements
of order $o\left({1/\sqrt{N}} \right)$. If this assumption is not met, i.e. we have  ${\mathbb E}(S) = a / \sqrt{N}$ with $a \gg 1 $, it means that the matrix $S /\sqrt{N}$ will have an
eigenvalue of order $a$, while the eigenvalues corresponding to the
planted signal will be $O(1)$. This means that for $a\gg 1$ the
eigenvalues corresponding to the signal will be subdominant and this
would require additional terms in the state evolution.

We know that ${\hat{x}}_i^{t+1} = f_{\rm in}(A_{X,i}^t,B_{X,i}^t)$ and
$\sigma_{x,i}^{t+1} =\frac{\partial f_{\rm in}}{\partial
  B}(A_{X,i}^t,B_{X,i}^t)$,
eqs.~(\ref{EquationTAP_UpdateMean}-\ref{EquationTAP_UpdateVariance}). Therefore
to compute the updated order parameters in the large $N$ limit one
needs to compute the probability distribution of $B_{X,i}^t$ and $A_{X,i}^t$
\begin{eqnarray}
P(B_{X,i}^t | {x_i^0} , Q_x^t , M_x^t , \Sigma_x^t)\, ,
\\
P(A_{X,i}^t | {x_i^0} , Q_x^t , M_x^t, \Sigma_x^t) \,.
\end{eqnarray}
Quantities $B_{X,i}^t$ and $A_{X,i}^t$ are defined by
eq. (\ref{Equation_B_i}) and (\ref{Equation_A_i}). Notably, by the
assumptions of belief propagation the terms in the sums on the right
hand side of eqs.~(\ref{Equation_B_i}) and (\ref{Equation_A_i}) are
independent. By the central limit theorem 
$B_{X,i}^t$ and $A_{X,i}^t$ then behave as Gaussian random variables. 

Therefore all one needs to compute is their mean and variance with
respect to the output channel. Using eq. (\ref{Equation_B_i}) we get
\begin{equation}
\mathbb{E}(B_{X,i}^t) = \frac{1}{\sqrt{N}}\sum\limits_{1 \leq k \leq N} \int {\rm d} Y_{ki} P_{\rm out}\left(Y_{ki} \middle| \frac{{x_k^{0,\top}} {x^0_i}}{\sqrt{N}} \right)  \left( \frac{\partial g}{\partial w}\right)_{Y_{ki},0}    {\hat{x}}_{k \rightarrow ki}^t \,.
\end{equation}
Let us now expand $P_{\rm out}$ around $0$
\begin{equation}
\mathbb{E}(B_{X,i}^t) = \frac{1}{\sqrt{N}}\sum\limits_{1 \leq k \leq
  N} \int {\rm d} Y_{ki} P_{\rm out}\left(Y_{ki} \middle| 0 \right)
\left[1 + \frac{{x_k^{0,\top}}   {x_i^0}}{\sqrt{N}} \left(
    \frac{\partial \log P_{\rm out}(Y_{ki}|w)}{\partial w}\right)_{Y_{ki},0} + O\left( \frac{1}{{N}}\right) \right] \left( \frac{\partial g}{\partial w}\right)_{Y_{ki},0}   {\hat{x}}_{k \rightarrow ki}^t \,.
\end{equation}
Using the above stated assumption of validity of eq. (\ref{E_1st_der})
we can simplify into 
\begin{equation}
\mathbb{E}(B_{X,i}^t) = \frac{1}{N \widehat\Delta}\sum\limits_{1 \leq k \leq  N} {\hat{x}}_{k}^t x_k^{0,\top}  {x^0_i} =   \frac{M_x^t}{\widehat\Delta}  x^0_i \,,
\end{equation}
where we used the definition (\ref{SE_m_def}) of the order parameter
$M^t$ and where we defined $\widehat\Delta$ via
\begin{equation}
\frac{1}{\widehat\Delta} \equiv \mathbb{E}_{P_{\rm out}}\left[ \left(
    \frac{\partial g(Y,w)}{\partial w}\right)_{Y,0}  \left( \frac{
     \partial \log(P_{\rm out}(Y|w))}{\partial w}\right)_{Y,0}\right] \,. \label{def_hat_Delta}
\end{equation}

Let us now compute the variance of $B_{X,i}^t$.
Using the assumption of belief propagation that messages incoming
to a variable are independent in the leading order we get that the
covariance of $B_{X,i}^t$ is the sum of all the covariance matrices of the
terms in the sum defining $B_{X,i}^t$.
\begin{equation}
{\rm Cov}(B^t_{X,i}) = \frac{1}{N}\sum\limits_{1 \leq i \leq N} {\rm Cov}(S_{ki}   \hat{x}_{k \rightarrow ki}^t) \,.
\end{equation}
Doing a similar computation as for the mean one gets in the leading
order 
\begin{equation}
{\rm Cov}(B^t_{X,i}) = \frac{1}{N \widetilde\Delta} \sum\limits_{1 \leq i \leq N}   {\hat{x}}^t_{k \rightarrow ki} {\hat{x}}_{k \rightarrow ki}^{t,\top}  = \frac{  Q_x^t  }{\widetilde\Delta} \,,
\end{equation}
where $\widetilde\Delta$ was introduced in (\ref{def_tilde_Delta}) and thanks to
self-averaging it also equals 
\begin{equation}
\frac{1}{\widetilde\Delta} = \mathbb{E}_{P_{\rm out}(Y|w)}\left[ \left( \frac{\partial g}{\partial w}\right)_{Y,0}^2\right] \,.
\end{equation}
Here one did not even have to expand $P_{\rm out}$ to second order,
the first order was enough to get the leading order of the variance.

The distribution of the $A^t_{X,i}$ now needs to be computed. Using
the definition of $A^t_{X,i}$ eq.~(\ref{Equation_A_i}) and the self-averaging of
section \ref{sec:self_av} we obtain directly that 
\begin{equation}
\mathbb{E}(A_{X,i}^t) = \frac{  Q_x^t  }{\widetilde\Delta} - \overline R\left( Q_x^t+\Sigma_x^t \right)\,,
\end{equation}
where $\overline R$ is defined in (\ref{def_overline_R}) and also equals 
\begin{equation}
\overline R = \mathbb{E}_{P_{\rm out}(Y|w)}\left[ \left( \frac{\partial g}{\partial w}\right)_{Y,0}^2 + \left( \frac{\partial^2 g}{\partial w^2}\right)_{Y,0}\right] \,.
\end{equation}
Here things are simpler then for $B^t_{X,i}$ since $A^t_{X,i}$ concentrates around its
mean, its variance is of smaller order. 

Overall, using (\ref{EquationTAP_UpdateMean}) and
(\ref{EquationTAP_UpdateVariance}) one gets for the state evolution
equations 
\begin{eqnarray}
M_x^{t+1} = \mathbb{E}_{x_0,W}\left[ f_{\rm in}^x\left(\frac{  Q_x^t  }{\widetilde\Delta} - \overline R (Q_x^t + \Sigma_x^t)  , \frac{  M_x^t  }{\widehat\Delta} x_0 +	\sqrt{\frac{Q_x^t}{\widetilde\Delta}}   W\right)x_0^\top \right] \,,
\label{SE_Equation_XX_M}
\\
Q_x^{t+1} = \mathbb{E}_{x_0,W} \left[ f_{\rm in}^x\left(\frac{  Q_x^t  }{\widetilde\Delta} - \overline R (Q_x^t + \Sigma_x^t)  , \frac{  M_x^t  }{\widehat\Delta} x_0 + \sqrt{\frac{Q_x^t}{\widetilde\Delta}}   W\right)f_{\rm in}^x(\cdots,\cdots)^\top \right]
\label{SE_Equation_XX_Q} \,,
\\
\Sigma_x^{t+1} = \mathbb{E}_{x_0,W} \left[ \frac{\partial f_{\rm in}^x}{\partial B}\left(\frac{  Q_x^t  }{\widetilde\Delta} - \overline R (Q_x^t + \Sigma_x^t)  , \frac{  M_x  }{\widehat\Delta} x_0 + \sqrt{\frac{Q_x^t}{\widetilde\Delta}}   W\right)\right]
\label{SE_Equation_XX_Sigma} \,,
\end{eqnarray}
where $W$ and $x_0$ are two independent random variables. $W$ is a
Gaussian noise of mean 0 and covariance matrix $I_r$ and $x_0$ is a
random variable of probability distribution $P_{X_0}$. The thresholding
function $f_{\rm in}^x$ is defined in eq.~(\ref{def_fin}).
This also allows us to know that at a fixed point the marginals $x_i$ will be distributed according to 
\begin{equation}
\hat{x_i}^{t=+\infty} = f_{\rm in}^x\left(\frac{
    Q_x}{\widetilde\Delta} - \overline R (Q_x + \Sigma_x)  , \frac{
    M_x  }{\widehat\Delta} x_0 + \sqrt{\frac{Q_x}{\widetilde\Delta}}
  W\right)\, , \label{Distribution_Marginals}
\end{equation}
where $W$ is a Gaussian variable and mean 0 and covariance $I_r$ and $x_0$ is taken with respect to $P_{X_0}$.
Let us state that the large $N$ limit of the ${\rm MSE}_X =\sum\limits_{i = 1 \ldots N}||{\hat{x}}_i -
x^0_i||_2^2 /N$ is computed from the order parameters
as
\be
{\rm MSE}_X = {\rm Tr}\left[ \langle x_0 x_0^{\top}\rangle -2M_x +
  Q_x\right]\, ,
\label{MSE_X}
\ee
where $\langle x_0 x_0^{\top}\rangle = \mathbb{E}_{{x_0}}(x_0
x_0^{\top})$ is the average with respect to the distribution $P_{X_0}$. 

Note that the state evolution equations only depend on the assumed and truth
noise channels through three variables $\widetilde\Delta$,
$\widehat\Delta$ and $\overline{R}$. In the Bayes-optimal case these
equations will simplify even further and the noise channel will be
described through one parameter $\Delta$, the Fisher information, this is derived in section \ref{SectionSimplificationSE}.
This universality with respect to the output channel has been observed
elsewhere in a special case of the present problem
\cite{deshpande2015finding} (see  e.g.  their  remark  2.5) in  the study of detection of a small hidden clique with approximate message passing.

Finally one additional assumption made in this whole section is that no
Replica Symmetry Breaking (RSB) appears. It is known that RSB does appear for
some regimes of parameters out of the equilibrium Bayes-optimal case. We let the
investigation of RSB in the context of low-rank matrix estimation for
future work, in the examples analyzed in this paper we will restrict
ourselves to the Bayes-optimal case where RSB at equilibrium cannot happen \cite{zdeborova2015statistical}.

\subsection{Summary for the bipartite low-rank matrix factorization}

State evolution can also be written similarly for the $UV^\top$ case.
In that case there are six order parameters
\begin{eqnarray}
M_u^t = \frac{1}{N}\sum\limits_{1 \leq i \leq N} u_i^t u^{0,\top}_i, \quad
Q_u^t &=& \frac{1}{N}\sum\limits_{1 \leq i \leq N} u_i^t u_i^{t,\top},\quad
\Sigma_u^t = \frac{1}{N}\sum\limits_{1 \leq i \leq N} \sigma_{u,i}^t\, ,
\\
M_v^t = \frac{1}{M}\sum\limits_{1 \leq j \leq M} v_j^t v_j^{0,\top},\quad
Q_v^t &=& \frac{1}{M}\sum\limits_{1 \leq j \leq M} v_j^t {v_j^{t,\top}},\quad
\Sigma_v^t = \frac{1}{M}\sum\limits_{1 \leq j \leq M} \sigma_{v,j}^t\,.
\end{eqnarray}
These order parameters are updated according to the following state
evolution  equations
\begin{eqnarray}
M_u^{t} &=& \mathbb{E}_{u_0,W} \left[ f_{\rm in}^u\left(\frac{ \alpha  Q_v^t }{\widetilde\Delta} - \alpha \overline R(Q_v^t + \Sigma_v^t)
, \alpha \frac{ M_v^t }{\widehat\Delta}u_0 + \sqrt{\frac{\alpha Q_v^t }{\widetilde\Delta}}W\right)u_0^\top \right]\,,
\label{SE_Equation_UV_M_U}
\\
Q_u^{t} &=& \mathbb{E}_{u_0,W}  \left[ f_{\rm in}^u\left(\alpha\frac{ Q_v^t }{\widetilde\Delta}- \alpha \overline R(Q_v + \Sigma_v)
, \alpha\frac{ M_v^t}{\widehat\Delta}u_0 + \sqrt{\frac{\alpha Q_v^t }{\widetilde\Delta}}W_v\right)f_{\rm in}^u(\cdots,\cdots)^\top \right]\,,
\label{SE_Equation_UV_Q_U}
\\
\Sigma_u^{t} &=& \mathbb{E}_{u_0,W}  \left[ 
\frac{\partial f_{\rm in}^u}{\partial B}\left(\alpha\frac{ Q_v^t }{\widetilde\Delta}- \alpha \overline R(Q_v^t + \Sigma_v^t)
, \alpha\frac{ M_v^t}{\widehat\Delta}u_0 + \sqrt{\frac{\alpha Q_v^t }{\widetilde\Delta}}W\right) \right]\,,
\label{SE_Equation_UV_Sigma_U}
\\
\label{SE_Equation_UV_M_V}
M_v^{t+1} &=& \mathbb{E}_{v_0,W} \left[ f_{\rm in}^v\left(\frac{ Q_u^t }{\widetilde\Delta}- \overline R(Q_u^t + \Sigma_u^t)
, \frac{ M_u^t}{\widehat\Delta}v_0 + \sqrt{\frac{Q_u^t }{\widetilde\Delta}}W\right)v_0^\top \right]\,,
\\
\label{SE_Equation_UV_Q_V}
Q_v^{t+1} &=& \mathbb{E}_{v_0,W} \left[ f_{\rm in}^v\left(\frac{ Q_u^t }{\widetilde\Delta}- \overline R(Q_u^t + \Sigma_u^t)
, \frac{ M_u^t}{\widehat\Delta}v_0 +  \sqrt{\frac{Q_u^t }{\widetilde\Delta}}W\right)f_{\rm in}^v(\cdots,\cdots)^\top\right]\,,
\\
\label{SE_Equation_UV_Q_Sigma}
\Sigma_v^{t+1} &=& \mathbb{E}_{v_0,W} \left[ \frac{\partial f_{\rm in}^v}{\partial B}\left(\frac{ Q_u^t }{\widetilde\Delta}- \overline R(Q_u^t + \Sigma_u^t)
, \frac{ M_u^t}{\widehat\Delta}v_0 +  \sqrt{\frac{Q_u^t }{\widetilde\Delta}}W\right)\right]
\,.
\end{eqnarray}
In these equations $W$, $u_0$ and $v_0$ are independent random variables,
$W$ is $r$ dimensional Gaussian variable of mean $\vec{0}$ and covariance matrix $I_r$,
$u_0$ and $v_0$ are sampled from density probability $P_{U_0}$ and $P_{V_0}$ respectively.

The large size limit of the ${\rm MSE}_U =\sum\limits_{i = 1 \ldots N}||{\hat{u}}_i -
u^0_i||_2^2 /N$ and  ${\rm MSE}_V = \sum\limits_{j = 1 \ldots
  M}||{\hat{v}}_j - v^0_j||_2^2/M $ can be computed from the order
parameters as
\begin{eqnarray}
{\rm MSE}_U &=& {\rm Tr}\left[\langle u_0 u_0^{\top}\rangle -2M_u +
                Q_u\right]\, ,
\label{MSE_U}
\\
{\rm MSE}_V &=& {\rm Tr}\left[\langle v_0 v_0^{\top}\rangle -2M_v + Q_v\right] \,,
\label{MSE_V}
\end{eqnarray}
with $\langle u_0 u_0^{\top}\rangle = \mathbb{E}_{{u_0}}(u_0
u_0^{\top})$ and $\langle v_0 v_0^{\top}\rangle = \mathbb{E}_{{v_0}}(v_0
v_0^{\top})$.

\subsection{Replica Free Energy}

State evolution can also be used to derive large size limit of the Bethe
free energy (\ref{TAP_XX_FreeEnergy}) and (\ref{TAP_UV_FreeEnergy})
defined as 
\begin{eqnarray}
\Phi_{{\rm RS},XX^\top} \equiv \lim\limits_{N \rightarrow + \infty} \frac{1}{N}\left\langle \log({Z}_X(Y)) - \sum\limits_{1 \leq i < j \leq N} g(Y_{ij},0) \right\rangle \,,
\\
\Phi_{{\rm RS},UV^\top} \equiv \lim\limits_{N \rightarrow + \infty} \frac{1}{N}\left\langle \log({Z}_{UV}(Y)) - \sum\limits_{1 \leq i \leq N,1 \leq j \leq M} g(Y_{ij},0) \right\rangle \,,
\end{eqnarray}
where the average is  taken with respect to density probability
(\ref{Boltz_XX}), or (\ref{Boltz_UV}), the $Z_{X}(Y)$ and $Z_{UV}(Y)$
are the corresponding partition functions. 
We subtract the constant $g(Y_{ij},0)$ for convenience in order to get a quantity that is self averaging in the large $N$ limit.

Alternatively to the state evolution, the average free energy can be derived using the
replica method as we summarize in Appendix \ref{appendix:ReplicaComputation}. The resulting
replica free energy for the symmetric  $X X^\top$ case 
is (assuming the replica symmetric ansatz to hold)
\begin{equation}
\Phi_{{\rm RS},XX^\top} = \max \left\{    \phi_{\rm RS}(M_x,Q_x,\Sigma_x),
	\frac{\partial \phi_{\rm RS}}{\partial M_x} =
\frac{\partial \phi_{\rm RS}}{\partial Q_x} =
\frac{\partial \phi_{\rm RS}}{\partial \Sigma_x} = 0 \right\}\, ,
\end{equation}
where 	
\begin{multline}
\phi_{\rm RS}(M_x,Q_x,\Sigma_x) = 
\frac{{\rm Tr}(  Q_x   Q_x^\top)}{ 4\widetilde\Delta}- \frac{{\rm Tr}(  M_x   M_x^\top)}{2 \widehat\Delta} - \frac{\overline R}{2} {\rm Tr}( (Q_x + \Sigma_x) (Q_x + \Sigma_x)^\top )
\\
+ \mathbb{E}_{W,x_0} \left[ \log\left({Z}_x\left(\frac{  Q_x
    }{\widetilde\Delta} - \overline R (Q_x + \Sigma_x)  , \frac{  M_x
    }{\widehat\Delta} x_0 + \sqrt{\frac{Q_x}{\widetilde\Delta}}
    W\right) \right)	\right] \, ,
 \label{BetheFreeEnergy_XX_DE}
\end{multline}
where the function ${Z}_x(A,B)$ is defined as the normalization in
eq.~(\ref{Define_F_in}).

For the bipartite $UV^\top$ case we have analogously for the replica free energy
\begin{equation}
\Phi_{{\rm RS},UV^\top} = \max \left\{    \phi_{\rm
    RS}(M_u,Q_u,\Sigma_u,M_v,Q_v,\Sigma_v), \frac{\partial \phi_{\rm
      RS}}{\partial M_u} = \frac{\partial \phi_{\rm RS}}{\partial Q_u}
  = \frac{\partial \phi_{\rm RS}}{\partial \Sigma_u} =  \frac{\partial
    \phi_{\rm RS}}{\partial M_v} = \frac{\partial \phi_{\rm
      RS}}{\partial Q_v} = \frac{\partial \phi_{\rm RS}}{\partial
    \Sigma_v} = 0 \right\}\, ,
\end{equation}
where
\begin{multline}
\phi_{\rm RS}(M_u,Q_u,\Sigma_u,M_v,Q_v,\Sigma_v) = 
\frac{\alpha {\rm Tr}( Q_v Q_u^\top)}{ 2\widetilde\Delta}- \frac{\alpha {\rm Tr}(M_v M_u^\top)}{\widehat\Delta} - \alpha \overline R {\rm Tr}((Q_v + \Sigma_v)(Q_u + \Sigma_u)^\top )
\\
+ \mathbb{E}_{W,u_0} \left[ \log \left({Z}_u\left( \frac{\alpha Q_v}{\widetilde\Delta} - \alpha \overline R(Q_v + \Sigma_v) , \frac{\alpha M_v}{\widehat\Delta} u_0 + W \sqrt{\frac{\alpha Q_v}{\widetilde\Delta}}\right)	\right) \right]
\\
+ \alpha \mathbb{E}_{W,v_0}
\left[ \log\left({Z}_v\left(\frac{Q_u}{\widetilde\Delta} - \overline R(Q_u +
    \Sigma_u) , \frac{M_u}{\widehat\Delta} v_0 +
    \sqrt{\frac{Q_u}{\widetilde\Delta}} W\right)	 \right) \right]\, ,
 \label{BetheFreeEnergy_UV_DE}
\end{multline}
with the function ${Z}_u(A,B)$ and ${Z}_v(A,B)$ also defined as the normalization in
eq.~(\ref{Define_F_in}). 

It is worth noting that there is a close link between the expression
of the replica free energy and the state evolution equations. Namely
fixed points of the state evolution equations are stationary points of
the replica free energy and vice-versa. Therefore, by looking for a stationary point of these equations one finds back the state evolution equations (\ref{SE_Equation_XX_M}-\ref{SE_Equation_XX_Sigma},\ref{SE_Equation_UV_M_U}-\ref{SE_Equation_UV_Q_Sigma}).

\subsection{Simplification of the SE equations}
\subsubsection{Simplification in the Bayes optimal setting}
\label{SectionSimplificationSE}

The state evolution equations simplify considerably when we restrict
ourselves to the Bayes-optimal setting defined in Sec.~\ref{sec:Bayes_optimal} by
eq.~(\ref{Bayes_optimal}). 

From the definitions of $\widehat\Delta$ in eq.~(\ref{def_hat_Delta}) and
$\widetilde\Delta$ in eq.~(\ref{def_tilde_Delta}) and using the identity (\ref{Bayes_optimal})
defining the Bayes optimal setting we obtain
\begin{equation}
\frac{1}{\widehat\Delta} = \frac{1}{\widetilde\Delta} =\frac{1}{\Delta} = \mathbb{E}_{P_{\rm out}(Y,w=0)}\left[\left( \frac{\partial g}{\partial w} \right)_{Y,w=0}^2\right] \,,
\end{equation}
where $\Delta$ is the Fisher information of the output channel defined
in eq.~(\ref{Define_Delta}). Note for instance that for the Gaussian
input channel (\ref{gauss_out}), $\Delta$ is simply the variance of the
Gaussian noise. The bigger the $\Delta$ the harder the inference
problem becomes. The smaller the $\Delta$ the easier the inference is.


Further consequence of having Bayes optimality (\ref{Bayes_optimal}) is that $\overline R = \mathbb{E}(R_{ij}) = 0$ as proven in equation \eqref{E_2nd_der}.
This simplifies greatly the state evolution equations into 
\begin{eqnarray}
M_x^{t+1} &=& \mathbb{E}_{x_0,W}\left[ f_{\rm in}^x\left(\frac{  Q_x^t  }{\Delta} , \frac{  M_x^t  }{\Delta} x_0 +	\sqrt{\frac{Q_x^t}{\Delta}}   W\right)x_0^\top \right]
\label{SE_Equation_XX_Maux} \,,
\\
Q_x^{t+1} &=& \mathbb{E}_{x_0,W} \left[ f_{\rm in}^x\left(\frac{  Q_x^t  }{\Delta} , \frac{  M_x^t  }{\Delta} x_0 + \sqrt{\frac{Q_x^t}{\Delta}}   W\right)f_{\rm in}^x(\cdots,\cdots)^\top \right]
\label{SE_Equation_XX_Qaux} \,,
\end{eqnarray}
where $x_0$ and $W$ are as before independent random variables, $W$ is
Gaussian of mean 0 and covariance matrix $I_r$, and $x_0$ has probability distribution $P_{X_0}$. 

Another property that arises in the Bayes optimal setting
\eqref{Bayes_optimal} and follows from the Nishimori condition
(\ref{Nish_gen}) and the definition of the order parameters $M_x$, $Q_x$ and $\Sigma_x^t$ in
(\ref{SE_m_def}-\ref{SE_S_def}) is that 
\begin{equation}
Q_x^t = M_x^t  = M_x^\top  \, , \quad \quad Q_x^t + \Sigma_x^t= Q_x + \Sigma_x
=\langle x_0 x_0^{\top}\rangle\, .   \label{SE_Nish_cond}
\end{equation}
Enforcing $Q^t_x=M^t_x$ simplifies the state evolution equations further
so that for the symmetric matrix factorization one gets
\begin{equation}
M_x^{t+1} = \mathbb{E}_{x_0,W}\left[ f_{\rm in}^x\left(\frac{  M_x^t
    }{\Delta} , \frac{  M_x^t  }{\Delta} x_0 +
    \sqrt{\frac{M_x^t}{\Delta}}   W\right)x_0^\top \right] \,, \label{SE_Nish_XX}
\end{equation}
where $x_0$ and $W$ are independent random variables distributed as above. 
For the rest of the article we define the Bayes-optimal state
evolution function $f^{\rm SE}_{P_X}$ for prior $P_X$
\begin{equation}
M_x^{t+1} = f^{\rm SE}_{P_X}\left( \frac{M_x^t}{\Delta} \right)\, , \label{f_SE_defined}
\end{equation}
Let us comment on the output channel universality as discussed in Sec.
\ref{sec:universality}. In the Bayes-optimal setting the channel
universality becomes particularly simple and striking. For an arbitrary output channel $P_{\rm out}(Y|w)$ (for which the
expansion done in section \ref{sec:universality} is meaningful) we have the
following  
\begin{itemize}
\item The Low-RAMP algorithm in the Bayes optimal setting depends on
  the noise channel only through the Fisher score matrix $S$ as
  defined in eq. (\ref{Define_S}). This is specified in section (\ref{Nish_TAP}).
\item The state evolution in the Bayes-optimal setting depends on the
  output channel through the Fisher information of the channel
  $\Delta$ (\ref{Define_Delta}) as described in section
  \ref{SectionSimplificationSE}. As a consequence the minimal
  achievable error, the minimal efficiently achievable error and all other
  quantities that can be obtained from the state evolution depend on
  the output channel only trough the Fisher information $\Delta$. 
\end{itemize}
The replica free energy  (\ref{BetheFreeEnergy_XX_DE}) in the
Bayes-optimal case becomes 
\begin{equation}
\phi_{\rm RS}(M_x) = \mathbb{E}_{W,x_0} \left[ \log\left({Z}_x\left(\frac{  M_x
    }{\Delta}, \frac{ M_x
    }{\Delta} x_0 + \sqrt{\frac{M_x}{\Delta}}
    W\right) \right) \right] - \frac{{\rm Tr}(  M_x   M_x^\top)}{4 \Delta} \, .
 \label{BetheFreeEnergy_XX_DE_Nish}
\end{equation}
This was first derived in \cite{lesieur2015mmse}, and proven for a special case in
\cite{deshpande2016asymptotic}, and later in full generality in \cite{krzakala2016mutual,LelargeMiolane16,miolane2017fundamental}.
We also remind that the order parameter $M_x$ used in the state
evolution is related to the mean-squared error as 
\be
{\rm MSE}_X = {\rm Tr}\left[ \langle x_0 x_0^{\top}\rangle -M_x \right]\, ,
\label{MSE_X_Nish}
\ee

For the bipartite vector spin models, $UV^\top$ case, the state
evolution in the Bayes optimal setting reads
\begin{eqnarray}
M_u^{t} &=& \mathbb{E}_{u_0,W} \left[ f_{\rm in}^u\left(\frac{ \alpha  M_v^t }{\Delta}
, \alpha \frac{ M_v^t }{\Delta}u_0 + \sqrt{\frac{\alpha M_v^t }{\Delta}}W\right)u_0^\top \right]\,,
\label{SE_Equation_UV_M_U_Nish}
\\
\label{SE_Equation_UV_M_V_Nish}
M_v^{t+1} &=& \mathbb{E}_{v_0,W} \left[ f_{\rm in}^v\left(\frac{ Q_u^t }{\Delta}, \frac{ M_u^t}{\Delta}v_0 + \sqrt{\frac{Q_u^t }{\Delta}}W\right)v_0^\top \right]\,.
\end{eqnarray}
The replica free energy (\ref{BetheFreeEnergy_UV_DE}) in the Bayes
optimal setting becomes 
\begin{multline}
\phi_{\rm RS}(M_u,M_v) = \mathbb{E}_{W,u_0} \left[  \log\left( {Z}_u\left( \frac{\alpha M_v}{\Delta}, \frac{\alpha M_v}{\Delta} u_0 + W \sqrt{\frac{\alpha M_v}{\Delta}}\right)	\right)\right]
\\ + \alpha \mathbb{E}_{W,v_0} \left[ \log \left({Z}_v\left(\frac{M_u}{\Delta},
    \frac{M_u}{\Delta} v_0 + \sqrt{\frac{M_u}{\Delta}} W\right) \right)
\right] - \frac{\alpha {\rm Tr}( M_v M_u^\top)}{ 2\Delta}\, ,
 \label{BetheFreeEnergy_UV_DE_Nish}
\end{multline}
where once again $M_u = M_u^\top$ and $M_v = M_v^\top$. The global
maximum of the free energy is asymptotically the equilibrium free
energy, the value of $M_u$  and $M_v$ at this maximum is related to
the MMSE via 
\begin{eqnarray}
                   {\rm MSE}_U &=& {\rm Tr}\left[\langle u_0
                                   u_0^{\top}\rangle  -M_u \right]\, ,
\label{MSE_U_Nish}
\\
{\rm MSE}_V &=& {\rm Tr}\left[\langle v_0
                                   v_0^{\top}\rangle -M_v \right] \,.
\label{MSE_V_Nish}
\end{eqnarray} 
Performance of the Low-RAMP in the limit of large system sizes is
given by the fixed point of the state evolution reached with
initialization where both $M_u$  and $M_v$ are close to zero.

\subsubsection{Simplification for the conventional Hamiltonian and randomly quenched disorder}

Another illustrative example of the state evolution we give in this
section is for the conventional Hamiltonian (\ref{can_Ham}) with randomly quenched
disorder, as this is the case most commonly considered in the existing
physics literature. In that case the model (\ref{Boltz_XX})
corresponds to a generic vectorial spin glass model. 
To take into account that the disorder is not planted, but random, we
plug into the generic state evolution
\begin{eqnarray}
P_{\rm out}(Y,w) &=& \frac{1}{\sqrt{2 \pi J^2}} \exp\left(
                     -\frac{Y^2}{2 J^2}\right) \label{random_channel}
\end{eqnarray}
such that $P_{\rm out}(Y,w)$ does not depend on $w$, meaning that the
disorder $Y$ is chosen independently, there is no planting. 
For the conventional Hamiltonian \eqref{can_Ham} and output channel (\ref{random_channel}) we obtain 
\begin{eqnarray}
\overline R = \frac{1}{\widetilde\Delta} &=& \mathbb{E}\left[ Y^2\right] = J^2\,,
\\
\frac{1}{\widehat\Delta} &=& 0 \, .
\end{eqnarray}

The state evolution \eqref{SE_Equation_XX_M} and
\eqref{SE_Equation_XX_Q} then becomes 
\begin{eqnarray}
M_x^{t+1} &=& 0 \, , \\ 
Q_x^{t+1} &=& \mathbb{E}_{W} \left[ f_{\rm in}^x\left(-J^2\Sigma_x^t, J\sqrt{Q_x^t}   W\right)f_{\rm in}^x(\cdots,\cdots)^\top \right]
\label{SE_Equation_XX_Q_Conventional} \,,
\\
\Sigma_x^{t+1} &=& \mathbb{E}_{W} \left[ \frac{\partial f_{\rm in}^x}{\partial B}\left(-J^2\Sigma_x^t, J\sqrt{Q_x^t}   W\right)\right]
\label{SE_Equation_XX_Sigma_Conventional} \,.
\end{eqnarray}
The free energy \eqref{BetheFreeEnergy_XX_DE} is given by 
\begin{multline}
\phi_{\rm RS}(M_x,Q_x,\Sigma_x) = 
\mathbb{E}_{W,x_0} \left[ \log\left( {Z}_x\left(-J^2  \Sigma_x,J \sqrt{Q_x}
    W\right) \right) \right] + 
\frac{J^2 {\rm Tr}(  Q_x   Q_x^\top)}{ 4} - \frac{J^2}{2} {\rm Tr}( (Q_x + \Sigma_x) (Q_x + \Sigma_x)^\top )
\, .
 \label{BetheFreeEnergy_XX_DE_Conventional}
\end{multline}

Specifically, for the Sherrington-Kirkpatrick model \cite{SherringtonKirkpatrick75}, where the
rank is one and the spins are Ising eq.~\eqref{Ising} with $\rho =
1/2$ the $f_{\rm in}^x(A,B)$ is given by (\ref{fin_SK_h0}), the state
evolution becomes  
\begin{eqnarray}
M_x^{t+1} &=& 0      \, ,    \label{SE_SK_M} \\
Q_x^{t+1} &=& \mathbb{E}_{W} \left[ \tanh\left(J \sqrt{Q_x^t}   W\right)^2 \right]
\label{SE_SK_Q} \,,
\\
\Sigma_x^{t+1} &=& 1- Q_x^t
\label{SE_SK_S} \,.
\end{eqnarray}
where $W$ is a Gaussian random variables of zero mean and unit
variance. 
With a free energy \eqref{BetheFreeEnergy_XX_DE} given by 
\begin{equation}
\phi_{\rm RS}(Q_x) = 
\frac{J^2 (1 - Q_x)^2 - J^2}{4} + \mathbb{E}_{W}\left[\log\left( \cosh\left(J\sqrt{Q_x}   W\right) \right)	\right] \,. \label{SE_RS_FREE_ENERGY}
\end{equation}
The reader will notice that these are just the replica symmetric
equations of the Sherrington-Kirkpatrick solution \cite{SherringtonKirkpatrick75}.

\section{General results about low-rank matrix estimation}

\subsection{Analyzis of the performance of PCA}
\label{sec:PCA_analysis}


In this section we analyze the performance of a maximum likelihood
algorithm by estimating the behavior of the (replica symmetric) state
evolution in the limit where the interactions are given by $\exp(
\beta g(Y,w))$ with $\beta \rightarrow +\infty$, and the prior does
not contain hard constraints and is independent of $\beta$. Note that PCA and
related spectral methods correspond to taking $g(Y,w) = -(Y -
w)^2/2$. The presented method allows us to analyze the property of the
generalized spectral method where $g(Y,w)$ can be taken to be any
function including for instance $g(Y,w) = -(D(Y)-w)^2/2$ which would
correspond to performing PCA on an elementwise function $D$ of the
matrix $Y_{ij}$, this can be for instance the Fisher score matrix $S$. 

\subsubsection{Maximum likelihood for the symmetric $XX^\top$ case}
Maximum likelihood can be seen within the Bayesian approach as
analyzing the following posterior for $\beta \to \infty$
\begin{equation}
P(X|Y) = \frac{1}{Z(Y)} \prod\limits_{1 \leq i \leq N}
\frac{\exp(-\Vert x_i \Vert^2_2/2)}{\sqrt{2 \pi}^r} \prod\limits_{1
  \leq i < j \leq N} \exp\left( \beta g\left(Y, \frac{x_i^\top
      x_j}{\sqrt{N}} \right) \right) \, . \label{Probability_PCA}
\end{equation}
The function $g(Y,w)$ defines the parameters $\widehat{\Delta}$ (\ref{def_hat_Delta}),
$\widetilde{\Delta}$ (\ref{def_tilde_Delta}), and $\overline{R}$ (\ref{def_overline_R}). We want to
analyse the overlap between $\hat{X}$ and $X_0$ in the limit $N
\rightarrow \infty$, $\beta \rightarrow \infty$ since then the
posterior will be dominated by the likelihood terms $g(Y,w)$.
We put
here a prior $P_X$ Gaussian to ensure that $Z(Y) < +\infty$. We could
have chosen any $\beta$-independent prior $P_X(x)$ as long as the support of 
$P_X(x)$ is the whole $\mathbb{R}^r$. As $\beta, N \rightarrow
\infty$ the details of $P_X(x)$ will be washed away. One can write the state evolution equations (\ref{SE_Equation_XX_M}-\ref{SE_Equation_XX_Sigma})
\begin{eqnarray}
M_x^{t+1} &=& \beta\Sigma_x^{t+1} \frac{M_x^t
              \Sigma_0}{\widehat{\Delta}}\, , 
\\
Q_x^{t+1} &=& \beta\Sigma_x^{t+1} \left[\frac{M_x^t \Sigma_0
              {M_x^t}^\top}{\widehat{\Delta}^2} +
              \frac{Q_x^t}{\widetilde{\Delta}} \right]
              \beta\Sigma_x^{t+1}\, ,
\\
\beta \Sigma_x^{t+1} &=& \left[
                         \frac{1}{\beta}+Q_x^t\left(\frac{1}{\widetilde{\Delta}}
                         - \overline{R} \right) - \Sigma_x^t \left(
                         \frac{(\beta - 1)}{ {\widetilde{\Delta} }} +
                         \overline{R} \right) \right]^{-1}\, ,
\end{eqnarray}
where $\Sigma_0 = \langle x_0 x_0^\top \rangle\in
\mathbb{R}^{r \times r}$. We take the limit of $\beta \rightarrow \infty$ to get 
\begin{eqnarray}
M_x^{t+1} &=& \Sigma'^{t+1} \frac{M_x^t \Sigma_0}{\widehat{\Delta}}\, ,
\label{PCA_M_x}
\\
Q_x^{t+1} &=& \Sigma_x'^{t+1} \left[\frac{M_x^t \Sigma_0
              {M_x^t}^\top}{\widehat{\Delta}^2} +
              \frac{Q_x^t}{\widetilde{\Delta}} \right]
              \Sigma_x'^{t+1}\, ,
\label{PCA_Q_x}
\\
\Sigma_x'^{t+1} &=& \left[ Q_x^{t} \left(\frac{1}{\widetilde{\Delta}}
                    - \overline{R} \right) -
                    \frac{\Sigma_x'^t}{\widetilde{\Delta}}
                    \right]^{-1}\, ,
\label{PCA_Sigma_x}
\end{eqnarray}
where $\Sigma'^t = \lim_{\beta \rightarrow +\infty} \beta \Sigma^t$. 

In general, effects of replica symmetry breaking have to be taken into
account in the analysis of maximum likelihood. One exception are the
spectral methods for which we take $g(Y,w) = (D(Y)-w)^2/2$, where $D$
is some element-wise function. In that case the maximum likelihood
reduces to computation of the spectrum of the matrix $D(Y)$. Obtaining
the spectrum is a polynomial problems which is a sign that no replica
symmetry breaking is needed to analyze the performance of the spectral
methods on element-wise functions of the matrix $Y$. 

Following the derivation of the state evolution, eq.~\eqref{Distribution_Marginals}, we get that at the
fixed point of the state evolution the spectral estimator $\hat{x}_i$ is distributed according to
\begin{equation}
\hat{x_i} = \Sigma_x'\left[ \frac{M_x}{\widehat{\Delta}} x^0_{i} +
  \sqrt{\frac{Q_x}{\widetilde{\Delta}}}W_i \right] \, ,
\label{EigenVector_PCA}
\end{equation}
where $x^0_{i}$ is the planted signal or the rank-one perturbation,
and the $W_i$ are independent (in the leading order in $N$) Gaussian variables of mean 0 and covariance matrix $I_r$.

\subsubsection{MSE achieved by the spectral methods}


When one uses spectral methods to solve a low rank matrix estimation
problem one computes the $r$ leading eigenvalues of the corresponding
matrix and then one is left with the problem of what to do with the eigenvectors.
Depending on what problem one tries to solve one can for instance
cluster the eigenvectors using the $k$-means algorithm. A more
systematic way is the following: We know by
\eqref{EigenVector_PCA} that the elements of the eigenvectors $\hat{x}_i$ can be written as a random
variable distributed as 
\begin{equation}
\hat{x}_i = \hat{M} x^0_{i} + \sqrt{\hat{Q}}W_i \, ,
\label{EigenVector_PCA_2}
\end{equation}
with $\hat{M}=\Sigma_x' M_x/\widehat{\Delta}$, and
$\hat{Q}=Q_x (\Sigma_x')^2 / \widetilde{\Delta} $, where $M_x$,  $Q_x$ 
and $\Sigma_x'$ are fixed points of the state evolution equations (\ref{PCA_M_x}-\ref{PCA_Sigma_x}). 
This formula allows us to approach the problem as a low-dimensional
Bayesian denoising problem. Writing 
\begin{equation}
P(\hat{x}_i,x^0_{i}) = P(\hat{x}_i|x^0_{i})P_{X_{0}}(x^0_{i}) =
P_{X_{0}}(x^0_{i}) \frac{1}{\sqrt{{\rm Det} \left(2 \pi \hat{Q}
    \right)}} \exp\left( -\left(\hat{x}_i^\top -
    x_{i}^{0,\top}\hat{M}^\top  \right) \hat{Q}^{-1}\left( \hat{x}_i -
    \hat{M}x^0_{i} \right)/2\right)\, ,
\end{equation}
one gets
\begin{equation}
P(x^0_{i}|\hat{x}_i) =
\frac{P(\hat{x}_i|x^0_{i})P_{X}(x^0_{i})}{P(\hat{x}_i)}
=\frac{1}{P(\hat{x}_i)} P_{X_{0}}(x^0_{i}) \frac{1}{\sqrt{{\rm Det}
    \left(2 \pi \hat{Q} \right)}} \exp\left(-\frac{
    \left(\hat{x}_i^\top - x_{i}^{0,\top}\hat{M}^\top  \right)
    \hat{Q}^{-1}\left( \hat{x}_i - \hat{M}x^0_{i} \right)}{2}\right)
\, .
\end{equation}
By taking the average with respect to the posterior probability
distribution one gets for the spectral estimator 
\begin{equation}
\mathbb{E}_{P(x^0_{i}|\hat{x}_i)} \left[ x^0_{i}\right] = f_{\rm
  in}^x\left( \hat{M}^\top \hat{Q}^{-1} \hat{M},
  \hat{M}^\top\hat{Q}^{-1}\hat{x}_i  \right) \, .\label{Mean_Posterior_PCA}
\end{equation}
By combining \eqref{Mean_Posterior_PCA} with \eqref{EigenVector_PCA_2}
one gets that the mean-squared error achieved by the spectral method is
given by 
\begin{equation}
{\rm MSE}_{\rm PCA} = {\mathbb E}_{x_0,W} \left\{ \left[ x_0 - f_{\rm
        in}^x\left(\hat{M}^\top \hat{Q}^{-1} \hat{M},
        \hat{M}^\top\hat{Q}^{-1} \hat{M}x_0 + \sqrt{\hat{M}^\top
          \hat{Q}^{-1} \hat{M}} W \right) \right]^2 \right\}\label{PCA_optimal_reconstruction}\,,
\end{equation}
where the $W$ are once again Gaussian variables of zero mean and unit
covariance, and the variable $x^{0}$ is distributed according to $P_{X_0}$.
In the figures presented in subsequent sections we evaluate the performance of PCA via \eqref{PCA_optimal_reconstruction}.






\subsection{Zero-mean priors, uniform fixed point and relation to
  spectral thresholds}
\label{sec:uniform}

This section summarizes properties of problems for which the prior
distribution $P_x$ has zero mean. We will see that in those cases a
particularly simple fixed point of both the Low-RAMP and its state
evolution exists. We analyze the stability of this fixed points, and
note that linearization around it leads to a
spectral algorithm on the Fisher score matrix. As a result we observe
equivalence between the corresponding
spectral phase transition and a phase transition beyond which Low-RAMP
performs better than a random guess based on the prior. 

We do stress, however, that the results of this section hold only when
the prior has zero mean, and do not hold for generic priors of
non-zero mean. So that the spectral
phase transitions known in the literature are in general not related to the
physically meaningful phase transitions we observe in the
performance of Low-RAMP or in the information theoretically best
performance.

\subsubsection{Linearization around the uniform fixed point}
\label{TrivialFixedPointSection}

From the definition of the thresholding function (\ref{def_fin}), it
follows that that $\hat{x_i} = 0,\: \forall 1 \leq i \leq n$ is a  fixed point of the self-averaged low-RAMP equations (\ref{EquationTAP_UpdateMean},\ref{EquationTAP_UpdateVariance},\ref{Low-RAMP_Bfully},\ref{Low-RAMP_Afully}) whenever
\begin{equation}
\int {\rm d} x P_X(x) x = \langle x \rangle = 0, \, \quad \, {\rm and} \quad \overline{R}
= 0\, .
\end{equation}
We will call $\hat{x_i} = 0,\: \forall 1 \leq i \leq n$ the {\it uniform
fixed point}. The interpretation of this fixed point is that according
to it there is no information about the planted configuration $X_0$ in
the observed values $Y_{ij}$ and the estimator giving the lowest error
is the one that simply sets every variable to zero. When this is the
stable fixed point with highest free energy then this is indeed the
Bayes-optimal estimator.

In previous work on inference and message passing algorithms \cite{krzakala2013spectral,zdeborova2015statistical} we
learned that when a uniform fixed point of the message passing update
exists it is instrumental to expand around it and investigate the
spectral algorithm to which such a procedure leads. 
We follow this strategy here and expand the Low-RAMP equations
around the uniform fixed point. In the linear order in $\hat x$, the term
$A^t_X$ is negligible, from the definition of the thresholding function (\ref{def_fin}) one gets
\begin{equation}
\hat X^{t+1} = \left( \frac{S}{\sqrt{N}}\hat X^{t} - \hat
  X^{t-1}\frac{\Sigma_x}{\widetilde \Delta} \right) \langle x x^\top \rangle\,.
\end{equation}
Where $\Sigma_x$ is the average value of the variance of the
estimators, as defined in (\ref{SE_S_def}), at the uniform fixed point. In the Bayes optimal setting we remind
from equation (\ref{SE_Nish_cond}) that we moreover have $\Sigma_x =
\langle x x^\top \rangle$.





If we consider this last equation as a fixed point equation for $\hat
X$ we see that  columns of $\hat X$ are related to the eigenvectors of
the Fisher score matrix $S$. Expanding around the uniform fixed point
the Low-RAMP equations thus yields a spectral algorithm that is
essentially PCA applied to the matrix $S$ (not the original dataset $Y$).

For the bipartite model ($UV^\top$ case) the situation is
analogous. The self-averaged Low-RAMP equations have a uniform fixed
point $\hat u_i=0 \, \, \forall i$ , $\hat v_j=0\, \,  \forall j$ when
$\overline{R} = 0$ and when the priors $P_U$ and  $P_V$ have zero mean. 
Expanding around this uniform fixed point gives a linear operator
whose singular vectors are related to the left and right singular
vectors of the Fisher score matrix $S$.

\subsubsection{Example of the spectral decomposition of the Fisher score matrix}

Spectral method always come to mind when thinking about estimation of
low-rank matrices. Analysis of the linearized Low-RAMP suggests that
in cases where we have some guess about the form of the output channel
$P_{\rm out}(Y|w)$ then the optimal spectral algorithm should not be ran on the data matrix $Y_{ij}$
but instead on the Fisher score matrix $S_{ij}$ defined by
(\ref{Define_S}). This was derived in \cite{lesieur2015mmse} and further studied in
\cite{perry2016optimality}. In this section we give an example of a case where the
spectrum of $Y_{ij}$ does not carry any information for some region of
parameter, but the one of
$S_{ij}$ does.  

Consider as an example the output channel to be 
\begin{equation}
P_{\rm out}(Y|w) = \frac{1}{2} \exp\left( -{|Y -w|}\right)\,. 
\end{equation}
This is just an additive exponential noise. The Fisher score matrix $S_{ij}$ for this channel is 
\begin{equation}
S_{ij} = \sign(Y_{ij}) \,.
\end{equation} 
Consider the rank one case when the true signal distribution $P_{X_0}$ is
Gaussian of zero mean and variance $\sigma$. 

\begin{figure}
\includegraphics[scale=1]{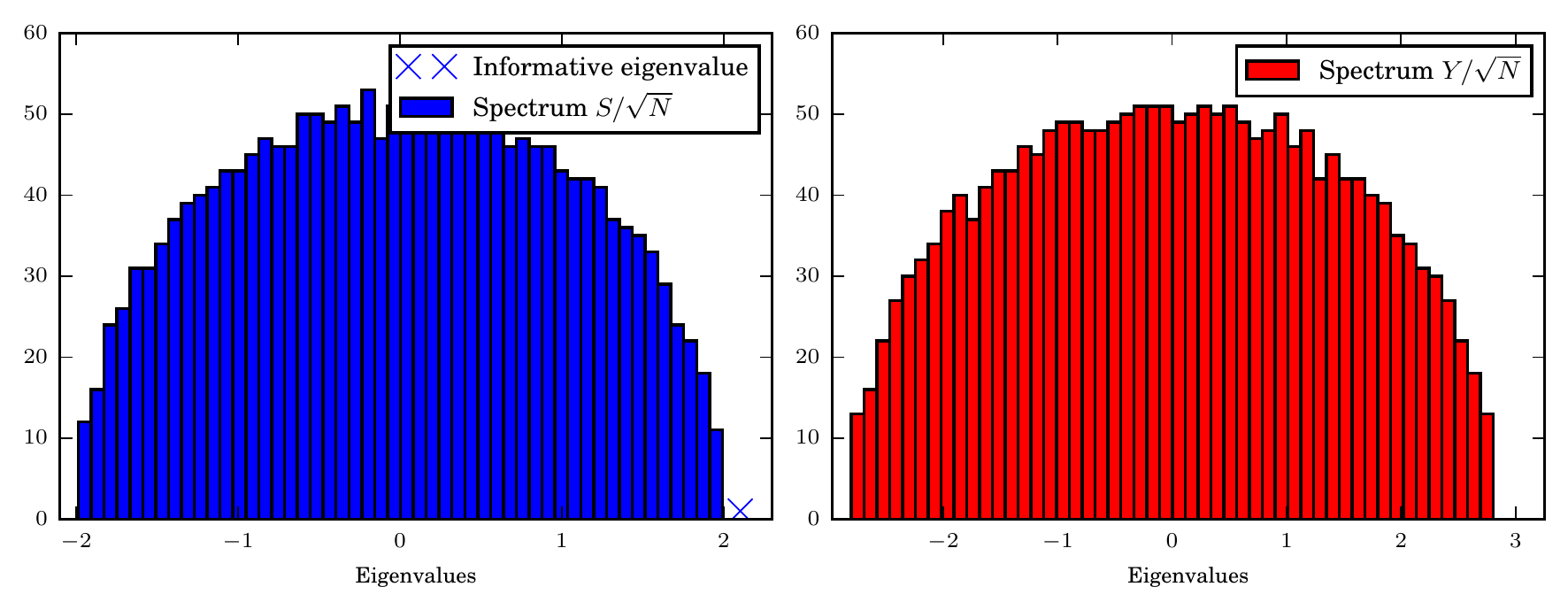}
\caption{Spectrum of the Fisher score matrix $S/\sqrt{N}$ and of the
  data matrix $Y/\sqrt{N}$ for the same
  instance of a problem with exponential output noise in the rank-one
  symmetric $XX^\top$ case. The planted configuration is generated
  from a Gaussian of zero mean and variance $1.4$. We see that an eigenvalue is out of
  the bulk for $S$ but not for $Y$. The data were generated on a system of size $N=2000$.}
\label{Spectrum_S}
\end{figure}

Now let us look at the spectrum of both $Y$ and $S$ in Fig.~\ref{Spectrum_S}.
We plot the spectrum of $S$ and $Y$ for $\sigma  = 1.4$. For this
value of variance we see that an eigenvalue associated with an
eigenvector that carries information about the planted configuration
gets out of the bulk of $S$ but not of $Y$. Even though some information on the signal was encoded into $Y$ in that specific case one had to take the absolute value of $Y$ to be able to recover an informative eigenvalue.

This situation can be quantified using the spectral analysis of
section \ref{sec:PCA_analysis} applied to two different noise channels
$g_1(Y,w) = -\beta(\sign(Y)-w)^2/2$ and 
$g_2(Y,w) = -\beta(Y-w)^2/4$ and taking the limit $\beta \rightarrow \infty$ as in \eqref{Probability_PCA}.
For these noise channels a theoretical analysis of the top eigenvectors of $S$ and $Y$ can be performed using theory presented in section \ref{sec:PCA_analysis}.
Taking square of \eqref{PCA_M_x} and dividing by \eqref{PCA_Q_x} 
one can show that the state evolution equation describing the overlap
of the top positive eigenvectors of $Y$ and $S$ is given by the only stable fixed point of the following update equation
\begin{eqnarray}
\frac{{(m_x^{t+1})}^2}{q_x^{t+1}} = \frac{ \frac{{m_x^{t}}^2}{q_x^{t}}
  \frac{\sigma^2 \widetilde{\Delta}}{\widehat{\Delta}^2}  }{1 +
  \frac{{m_x^{t}}^2}{q_x^{t}} \frac{\sigma \widetilde{\Delta}
  }{\widehat{\Delta}^2} } \, , 
\\
\widehat{\Delta} = \widetilde{\Delta} = 1 : {\rm for}\, g_1(Y,w)\, ,
\\
\widehat{\Delta} = \widetilde{\Delta} = 2 : {\rm for}\, g_2(Y,w)\, ,
\end{eqnarray}
where $\widehat{\Delta}$ and $\widetilde{\Delta}$ are computed when $\beta = 1$.
The trivial fixed point $\frac{{m_x^t}^2}{q_x^t} = 0$ of this equation is unstable as soon as
\begin{equation}
\sigma^2 \geq \frac{  \widehat{\Delta}^2}{\widetilde{\Delta}} \, .
\end{equation}
This analysis tells us that the  top eigenvectors are correlated with the planted solution $x^{0}$ when $\sigma > 1$ for $S$ and $\sigma >  2$ for $Y$.
Therefor for $\sigma=1.4$ the Fisher score matrix has informative
leading eigenvector, while the top eigenvector of $Y$ does not contain any information on the planted solution.

\subsubsection{Stability of the uniform fixed point in Bayes-optimal
  setting}
\label{TrivialFixedPointSection_stab}

In this section we restrict for simplicity to the Bayes optimal
setting defined by
eq. (\ref{Bayes_optimal}).  As we derived in section \ref{StateEvolutionSE} the 
evolution of the Low-RAMP algorithm can be tracked using the state
evolution equations. In the Bayes optimal case we have
$\overline R=0$ therefore the (sufficient) condition for the existence of the
uniform fixed point is to have prior $P_X$ (or both $P_U$ and $P_V$)
of zero mean. The existence of a uniform fixed point of the
Low-RAMP algorithm translates into the existence of a fixed point of
the state evolution with $M^t_x=0$ for the symmetric case (or
$M^t_u=M^t_v=0$ for the bipartite case).

The stability of this fixed point is analyzed by expanding in linear order the state
evolution equation (\ref{SE_Nish_XX}) around the uniform fixed point, taking
into account the definition of the thresholding function (\ref{def_fin}). For the $XX^\top$ case this gives
\begin{equation}
M_x^{t+1} = \frac{\Sigma_{x} M_x^t \Sigma_{x}}{\Delta} \,,
\label{Expansion_SE_XX}
\end{equation}
where  $\Sigma_{x}$ in the Bayes-optimal case is the covariance matrix of
the signal (and prior) distribution as given by
eq. (\ref{SE_Nish_cond}). Calling $\lambda^x_{\rm max}$ the largest
eigenvalue of the covariance of the distribution of
the signal-elements, $\Sigma_{x}$, we obtain a simple criterion for the stability of the uniform fixed point 
\begin{equation}
\left\{
\begin{matrix}
\Delta_c = (\lambda^x_{\rm max})^2 < \Delta \Rightarrow {\rm stable}
\\
\Delta < \Delta_c = (\lambda^x_{\rm max})^2 \Rightarrow {\rm unstable}
\end{matrix}
\right. \label{Trivial_XX_Stability} \,.
\end{equation}
It is useful to specify that for the rank-one case where both
$\Sigma_x$ and $M_x$ are scalars we get $\Sigma_x = \langle 
x^2_0 \rangle$ to be the variance of the prior distribution $P_x$. 
For the rank one, $r=1$, case the stability criterium becomes 
\begin{equation}
\left\{
\begin{matrix}
\Delta_c = \langle 
x^2_0 \rangle^2 < \Delta \Rightarrow {\rm stable}
\\
\Delta < \Delta_c = \langle 
x^2_0 \rangle^2 \Rightarrow {\rm unstable}
\end{matrix}
\right. \label{Trivial_XX_Stability_one} \,.
\end{equation}

Interestingly, the criteria (\ref{Trivial_XX_Stability}) is the same as
the criteria for the spectral phase transition of the Fisher score
matrix~$S$. When the uniform fixed point is not stable the Fisher
score matrix has an eigenvalue going out of the bulk
\cite{baik2005phase,hoyle2004principal}. We stress here that this analysis
is particular to signals of zero mean. If the mean is non-zero the
spectral threshold does not change, but the Bayes optimal performance
gets better and hence superior to PCA. 

The critical value of $\Delta_{c}$ separates two parts of the phase
diagram
\begin{itemize}
\item For $\Delta > \Delta_c$ inference is algorithmically hard or
  impossible. The Low-RAMP algorithm (and sometimes it is
conjectured that all other polynomial algorithms) will not be able to
get a better MSE than corresponding to random guessing from the prior
distribution.
\item For $\Delta < \Delta_c$ inference better than random guessing is
  algorithmically efficiently tractable. The Low-RAMP and also PCA give an MSE strictly better that
random guessing from the prior. 
\end{itemize}

For the bipartite case, $UV^\top$, the stability analysis is a tiny
bit more complicated since there are two order parameters $M_u^{t} $
and $M_v^{t}$. Linearization of the state evolutions leads to 
\begin{eqnarray}
M_u^{t} &=& \alpha \frac{\Sigma_{u} M_v^t \Sigma_{u}}{\Delta} \,,
\label{Expansion_DE_UV_1}
\\
M_v^{t+1} &=& \frac{\Sigma_{v} M_u^t \Sigma_{v}}{\Delta} \,,
\label{Expansion_DE_UV_3}
\end{eqnarray}
where in the Bayes-optimal case the $\Sigma_{u} $ and $\Sigma_{v} $ are
simply the covariances of the prior distribution $P_U$ and $P_V$, i.e.  
$
\Sigma_{u} =  \langle u u^\top \rangle
$
and 
$
\Sigma_{v} = \langle v v^\top \rangle
$.
By replacing (\ref{Expansion_DE_UV_3}) in (\ref{Expansion_DE_UV_1}) one gets
\begin{equation}
M_u^{t+1} = \left(\frac{\sqrt{\alpha} {\Sigma_{u} {\Sigma_{v}} }}{\Delta}
\right) M_u^t \left(\frac{\sqrt{\alpha} {\Sigma_{u} {\Sigma_{v}}}}{\Delta} \right)^\top \,.
\end{equation}
Calling $\lambda_{\rm max}^{uv}$ the largest eigenvalue of the matrix $\Sigma_{u} {\Sigma_{v}} $
gives us the stability criteria of the uniform fixed  point in the
   bipartite case as
\begin{equation}
\left\{
\begin{matrix}
\Delta_c = \sqrt{\alpha} \lambda_{\rm max}^{uv} < \Delta \Rightarrow {\rm Stable}
\\
\Delta < \Delta_c = \sqrt{\alpha} \lambda_{\rm max}^{uv} \Rightarrow {\rm Unstable}
\end{matrix}
\right. \label{Trivial_UV_Stability} \,.
\end{equation}

Also this criteria agrees with the criteria for spectral phase transition for the
Fisher score matrix and also here $\Delta_{c}$ separates two parts of
the phase diagram, one where estimating the signal betten than
randomly sampling from the prior distribution is not posible with
Low-RAMP (and conjecturally with no other polynomial algorithm), and
another where the MSE provided by Low-RAMP or PCA is strictly better
than the one achieved by guessing at random.

\subsection{Multiple stable fixed points: First order phase
  transitions}
\label{sec:1st_order}

The narrative of this paper is to transform the high-dimensional
problem of low-rank matrix factorization to analysis of stable fixed point of
the Low-RAMP algorithms and correspondingly of the state evolution
equations. A case that deserves detailed discussion is when there
exists more than one stable fixed point. The present section is
devoted to this discussion in the Bayes-optimal setting, where the replica symmetric assumption is
fully justified and hence a complete picture can be obtained.

We encounter two types of situations with multiple stable fixed points
\begin{itemize}
\item {\bf Equivalent fixed points due to symmetry:} The less
  interesting type of multiple fixed points arises when there is an
  underlying symmetry in the definition of the problem, then both the
  state evolution and and Low-RAMP equations have multiple fixed
  points equivalent under the symmetry. All these fixed points have the
  same Bethe and replica free energy. One example of such a symmetry is a Gaussian prior of zero mean and
isotropic covariance matrix. Then there is a global rotational
symmetry present. Another example of a symmetry is given by the community detection
problem defined in section \ref{ParagraphCommunityDetection}. 
In the case of $r$ symmetric equally sized groups with connectivity
matrix (\ref{ConnectivityMatrix}) there is a permutational symmetry between communities
so that any fixed of the state evolution or Low-RAMP equations exists in $r!$ versions.
\item {\bf Non-equivalent fixed points}. More interesting case is when
  Low-RAMP and the state evolution equations have multiple stable
  fixed points, not related via any symmetry, having in general
  different free-energies. This is then related to phenomenon that is
  in physics called the {\it first order phase transition}. This is the type that we will discuss in detail in the present section. 
\end{itemize}

In general, it is the fixed point with highest free energy that
provides the true marginals of the posterior distribution (we remind
change of sign in our definition of the free energy w.r.t. the
standard physics definition). From an algorithmic perspective, it is
the fixed point achieved from uninformed initialization (see below)
that gives a way to compute the error achievable by the Low-RAMP
algorithm. A conjecture that appears in a number of papers analyzing
Bayes optimal inference on random instances is that the error reached
by Low-RAMP is the best achievable with a polynomial algorithm. Note,
however, that replica symmetry breaking effects might play a role out
of the equilibrium solution and hence may influence the properties of
the best achievable mean-squared-error. 

\subsubsection{Typical first order phase transition: Algorithmic interpretation}
\label{SubSectionFirstOrderGeneral}

The concept of a first order phase transition is best explained on a
specific example. For the sake of the explanation we will consider the
symmetric $XX^\top$ case in the Bayes optimal setting. For the purpose
of giving a specific example we consider the Gaussian output
channel with variance of the noise $\Delta$, the signal is drawn from
the spiked (i.e. $r=1$) Rademacher-Bernoulli model with fraction of non-zeros being $\rho=0.08$. 

\begin{figure}
\includegraphics[scale=1]{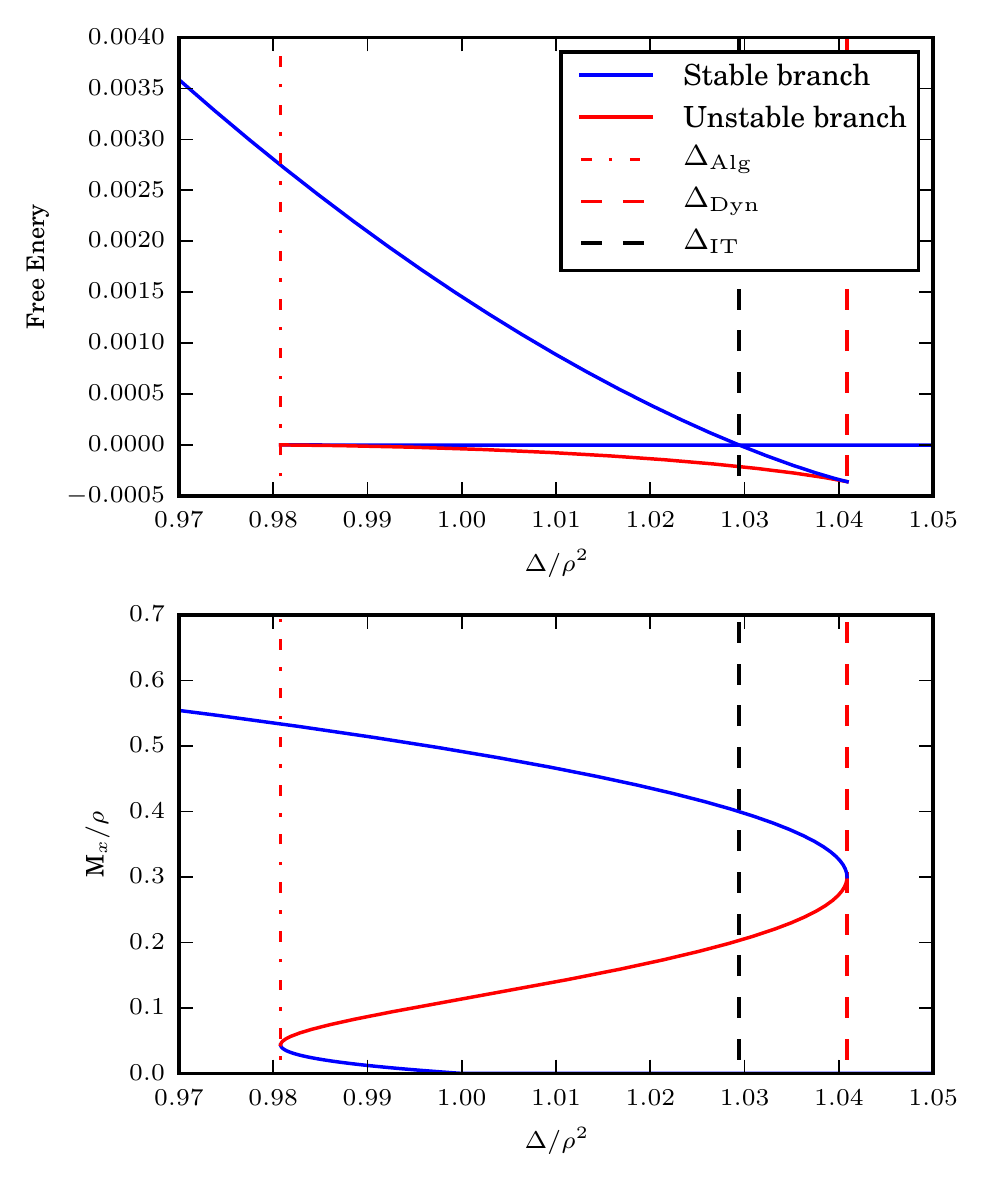}\caption{
In the lower pannel, we plot as a function of $\Delta/\rho^2$ all the fixed points $M_x/\rho$ of state evolution
equations  (\ref{SE_Nish_XX})  for the spiked Rademacher-Bernoulli
model with fraction of non-zero being $\rho=0.08$. Numerically stable
fixed points are in blue, unstable in red. We remind that the order
parameter $M_x$ is related to the mean-squared error as ${\rm MSE}_X = {\rm Tr}\left[ \langle x_0 x_0^{\top}\rangle -M_x \right]$. In the upper pannel, we
plot the replica free energy (\ref{BetheFreeEnergy_XX_DE_Nish})
corresponding to these fixed points, again as a function of
$\Delta/\rho^2$. 
When there are multiple stable fixed points to the SE equations the
one that corresponds to performance of the Bayes optimal estimation is
the one with the largest free energy (we remind that with respect to
the most common physics notation we defined the free energy as the positive
logarithm of the partition function).
The $\Delta$ for which these two branches cross in free-energy is
called the information theoretic phase transition, $\Delta_{\rm IT}$.
The two spinodal transitions  $\Delta_{\rm Alg}$ and $\Delta_{\rm
  Dyn}$ are where the lower MSE an higher MSE stable fixed point disappears. The $\Delta_c/\rho^2=1$ corresponds to the spectral transition at which the uniform fixed point become unstable.
}
\label{MinFreeEnergyFirstOrderExplanation}
\end{figure}

In Fig.~\ref{MinFreeEnergyFirstOrderExplanation} we plot all the fixed
points of the state evolution equation and the corresponding
value of the replica free energy (\ref{BetheFreeEnergy_XX_DE_Nish}) as a function
of $\Delta/\rho^2$. The equations for the state evolution specific to
the spiked Rademacher-Bernoulli model are given by
\eqref{SE_Rademacher_Bernoulli}. For this model the
uniform fixed point exists and is stable down to $\Delta_c = \rho^2$, eq.~(\ref{Trivial_XX_Stability}). The numerically
stable fixed points are drawn in blue, the unstable ones in red.  We
focus on the interval of $\Delta$ where more than one stable
fixed point of the state evolution (\ref{SE_Nish_XX})  exists. 
We use the example of the spiked Rademacher-Bernoulli model for the purpose
of being specific in 
figure~\ref{MinFreeEnergyFirstOrderExplanation}. The definitions and
properties defined in the rest of this section are generic and apply
to all the settings considered in this paper, not only to the spiked Rademacher-Bernoulli model.

Let us define two (possibly equal) stable fixed points of the state evolution as
follows: 
\begin{itemize}
\item $M_{\rm Uninformative}(\Delta)$ is the fixed point of the state
  evolution one reaches when initialising the $M^{t=0} = \epsilon I_r$
  (with $\epsilon$ very small and positive, $I_r$ being the identity matrix). We call this the uninformative initialisation.
\item $M_{\rm Informative}(\Delta)$ is the fixed point of the state
  evolution one reaches when initialising the $M^{t=0} = \langle x
  x^\top \rangle$. This is the informative initialisation where
  we start as if the planted configuration was known.
\end{itemize}
In principle there could be other stable fixed points apart of $M_{\rm
  Uninformative}(\Delta)$ and $M_{\rm Informative}(\Delta)$, but among
all the examples that we analyzed in this paper we have not observed	
any such case. In this paper we hence discuss only the case with at
most two stable fixed points, keeping in mind that if more stable
fixed points exist than the theory does apply straightforwardly as well (the physical
fixed point is still the one of highest free energy), and there could
be several first order phase transitions following each other as
$\Delta$ is changed. 

With two different stable fixed points existing for some values of $
\Delta$ we observe three critical values defined as follows
\begin{itemize}
\item $\Delta_{\rm Alg}$ called the {\it algorithmic spinodal
    transition}, is the value of $\Delta$ at which the fixed
  point $M_{\rm Uninformative}$ stops existing and becomes equal to
  $M_{\rm Informative}$. 
\item $\Delta_{\rm IT}$ called the {\it information theoretic phase
    transition}, is the value of $\Delta$ at which the two fixed
  points $M_{\rm Uninformative} \neq M_{\rm Informative}$  exist
  and have the same replica free energy (\ref{BetheFreeEnergy_XX_DE_Nish}).
\item $\Delta_{\rm Dyn}$ called the {\it dynamic spinodal transition},
  is the value of $\Delta$ at which the fixed point $M_{\rm
    Informative}$ stops existing and becomes equal to $M_{\rm
    Uninformative}$.  
\end{itemize}
In Fig.~\ref{MinFreeEnergyFirstOrderExplanation} these three
transition are marked by vertical dashed lines, the information
theoretic transition in black and the two spinodal transition in red. 

We recall that for the priors of zero mean, where the uniform fixed points
discussed in section \ref{TrivialFixedPointSection} exists, the
stability point of the uniform fixed point $\Delta_c$ (\ref{Trivial_XX_Stability}) is in general unrelated to the
$\Delta_{\rm Alg}$, $\Delta_{\rm IT}$ and $\Delta_{\rm Dyn}$. In the
example presented in Fig.~\ref{MinFreeEnergyFirstOrderExplanation} of spiked Rademacher-Bernoulli model at $\rho=0.08$ we have
$\Delta_{\rm Alg} < \Delta_c < \Delta_{\rm IT}$. In general the
position of $\Delta_c$ with respect to $\Delta_{\rm Alg}$,
$\Delta_{\rm IT}$ and $\Delta_{\rm Dyn}$ can be arbitrary, in the
following sections we will observe several examples. A notable
situation is when $\Delta_c = \Delta_{\rm Alg}$, cases where this
happen are discussed in section \ref{FirstOrderCriteriaSubSection}. It
should be noted that this is the case in the community detection
problem, and since this is a well known and studied example it is sometimes presented in the literature as the generic case. But from the numerous
examples presented in this paper we see that cases where $\Delta_c
\neq \Delta_{\rm Alg}$ are also very common.

Phase transitions are loved and cherished in physics, in the context
of statistical inference the most intriguing properties related to
phase transitions is their implications in terms of average computational
complexity. Notably in the setting of Bayes-optimal low-rank matrix
factorization as studied in this paper we distinguish two different regions
\begin{itemize}
   \item {\bf Phase where Low-RAMP is asymptotically Bayes-optimal:} For $\Delta \ge \Delta_{\rm IT}$ and for $\Delta \le \Delta_{\rm
     Alg}$ the Low-RAMP algorithm in the limit of large system sizes
   gives the information theoretically optimal performance. This is 
   either because the fixed point it reaches is unique (up to symmetries) or because it
   is the one with larger free energy. 
   \item {\bf The hard phase:} For $ \Delta_{\rm Alg} < \Delta < \Delta_{\rm IT}$ the
     estimation error achieved by the Low-RAMP algorithm is strictly
     larger that the lowest information-theoretically achievable
     error. On the other hand, and in line with other statistical physics works on inference
     problems in the Bayes optimal setting, we conjecture that in this
     region no other polynomial algorithm that would achieve better
     error than Low-RAMP exists. This conjecture could be slightly
     modified by the fact that the branch of stable fixed points that
     do not correspond to the MMSE could present aspects of replica
     symmetry breaking which could modify its
     position. Investigation of this is left for future work. 
\end{itemize}

From a mathematically rigorous point of view the results of this paper divide
into three parts: 
\begin{itemize}
\item (a) Those that are rigorous, known from existing literature that
  is not related to statistical physics considerations. An example is
  given by the performance of the spectral methods that is better that
  random guessing for $\Delta < \Delta_c$ \cite{baik2005phase}.
\item (b) A second part regroups the results directly following from
  the analysis of the Bayes-optimal Low-RAMP and the MMSE that are not
  presented rigorously in this paper, but were made rigorous in a
  series of recent works
  \cite{krzakala2016mutual,barbier2016mutual,LelargeMiolane16,miolane2017fundamental}. Most
  of these results are proven and although some cases are still
  missing a rigorous proof, it is safe to assume that it is a question
  of time that researchers will fill the corresponding gaps and weaken
  the corresponding assumptions. For instance the state evolution of
  Low-RAMP \cite{rangan2012iterative,javanmard2013state,DeshpandeM14},
  its Bayes-Optimality in the easy phase (at least for problem where
  the paramagnetic fixed point is not symmetric) and the value of the
  information theoretically optimal MMSE
  \cite{krzakala2016mutual,deshpande2016asymptotic,barbier2016mutual,LelargeMiolane16,miolane2017fundamental}
  are all proven rigorously.
\item (c) The third kind of claims are purely conjectures.  For
  instance, the claim that among all polynomial algorithms the
  performance of Low-RAMP cannot be improved (so that the hard phase
  is indeed hard). Of course proving this in full generality would
  imply that P$\neq$NP and so we cannot expect that such a proof would
  be easy to find. At the same time from a broad perspective of
  understanding average computational complexity this is the most intriguing
  claim and is worth detailed investigation and constant aim to find a
  counter-example.
\end{itemize}

\subsubsection{Note about computation of the first order phase
  transitions}
\label{TechniqueFirstOrderTransitionSection}
In this section we discuss how to solve efficiently the SE equations
in the $XX^\top$ Bayes optimal case for rank one. 
We notice that for a given prior distribution $P_X$ the only way the
symmetric Bayes optimal state evolution, eq.~ (\ref{SE_Nish_XX}),
depends on the noise parameter $\Delta$ is via the ratio
$m^t/\Delta$. One can write the SE equations in the form (\ref{f_SE_defined})
\begin{equation}
m^{t+1} = f_{P_X}^{\rm SE}\(( \frac{m^t}{\Delta}\)) \, .
\label{SE_Simplified_FirstOrder}
\end{equation}
Let us further define fixed points of (\ref{SE_Simplified_FirstOrder})
in a parametric way 
\begin{eqnarray}
\Delta = \frac{f_{P_X}^{\rm SE}(x)}{x}\, ,
\\
m = f_{P_X}^{\rm SE}(x)\, .
\end{eqnarray}
To get a fixed point $(m,\Delta)$ of  (\ref{SE_Simplified_FirstOrder}) we choose a
value of $x$ and compute $(f(x),f(x)/x)$. We observe from the form of
the state evolution (rank-one symmetric Bayes optimal case) that
$f(x)$ is a non-decreasing function of $x$. We further observe that  
$(m,\Delta)$ is a stable fixed points if and only if $\partial_x
\Delta(x) <0$. The two spinodal thresholds $\Delta_{\rm Alg}$ and
$\Delta_{\rm Dyn}$ are defined by loss of existence of corresponding
stable fixed points and they can hence be computed as
\begin{eqnarray}
\Delta_{\rm Alg},\Delta_{\rm Dyn} \in \left\{\Delta(x), x \in \mathbb{R}^+, \frac{\partial \Delta(x)}{\partial x} = 0  \right\}\,.
\label{Delta_Static_x}
\end{eqnarray}

The information theoretic transition $\Delta_{\rm IT}$ relies on
computation of the replica free energy (\ref{BetheFreeEnergy_XX_DE_Nish}). Using the fact that
\begin{equation}
\frac{\partial \phi(m,\Delta)}{\partial m} = \frac{1}{2 \Delta} \left( f_{P_X}^{\rm SE}\left(\frac{m}{\Delta} \right) - m\right)
\label{FreeEnergyEquationDeriv}
\end{equation}
allows us to compute the difference in energy between a fixed point
$m,\Delta$ and the uniform fixed point $m=0$ as  
\begin{equation}
\phi(m(x),\Delta(x)) - \phi(0,\Delta(x)) = \frac{1}{2	}\left[\int_0^x {\rm d}u \,f_{P_X}^{\rm SE}(u)	- \frac{x f^{\rm SE}(x)}{2} \right] \,. \label{FreeEnergyEquationStatic}
\end{equation}
The $x_{\rm IT}$ for which \eqref{FreeEnergyEquationStatic} is zero then gives
the information theoretic phase transition $\Delta_{\rm IT}=f(x_{\rm
  IT})/x_{\rm IT}$. 

\subsubsection{Sufficient criterium for existence of the hard phase}
\label{FirstOrderCriteriaSubSection}

This section is specific to the Bayes-optimal cases when the
prior $P_X$ has a zero mean and hence the uniform fixed point of the
Low-RAMP and the state evolution exists. 
In section \ref{TrivialFixedPointSection_stab} we derived that
the uniform fixed point is stable at $\Delta > \Delta_c$ and
unstable for $\Delta<\Delta_c$. It follows from the theory of
bifurcations that the critical point where a fixed-point changes
stability must be associated with an onset of another close-by fixed
point. In general there are two possibilities
\begin{itemize}
\item {\bf 2nd order bifurcation.} If the fixed point close to the uniform fixed point departs in the
  direction of smaller $\Delta<\Delta_c$, where the uniform fixed point is
  unstable, then this close-by fixed point is stable. This case
  corresponds to Fig.~\ref{MinFreeEnergyFirstOrderExplanation}. Behaviour in the vicinity of the uniform
  fixed point then does not let us distinguish between (a)
  existence of a first order phase transition at lower $\Delta$ (as in
  Fig.~\ref{MinFreeEnergyFirstOrderExplanation}), or (b)
  continuity on the MMSE down to $\Delta=0$ with no algorithmically hard phase
  existing in that case. 
\item {\bf 1st order bifurcation.}  If the fixed point close to the uniform fixed point departs in the
  direction of larger $\Delta>\Delta_c$, where the uniform fixed point is
  stable, then this close-by fixed point is unstable. In that case, the
  fixed point that is stable from $\Delta<\Delta_c$ is not close to
  the uniform fixed point and this case forces existence of a first
  order phase transition with $\Delta_c = 
\Delta_{\rm Alg}$.
\end{itemize}

Expansion of the state evolution (\ref{SE_Nish_XX})  around the
uniform fixed point gives us a closed form criteria to distinguish whether the phase transition happening at
$\Delta_c$ is a 1st or 2nd order bifurcation. 
In case of a 2nd order bifurcation, this expansion allows us to
compute the mean squared error obtained with Low-RAMP close to
$\Delta_c$. 

For specificity we consider the rank-one, $r=1$ case, and expand
eq.~(\ref{SE_Nish_XX}) to 2nd order to get
\begin{equation}
m^{t+1} = m^t \frac{ \langle x_0^2 \rangle^2 }{\Delta} - \frac{(m^t)^2}{\Delta^2}\left(\langle x_0^2 \rangle^3 - \frac{{\langle x_0^3 \rangle}^2}{2}\right) + O\left( (m^t)^3\right) \,,
\label{SecondOrderExpansion}
\end{equation}
where the mean $\langle \cdots \rangle$ of $x_0$ are taken with
respect to $P_{X_0} = P_X$. This is done by expanding $f_{\rm in}^x(A,B)$ to order $4$ in $B$ and $2$ in $A$. All the derivatives
\begin{equation}
\forall i,j \quad \frac{\partial^{i+j} f_{\rm in}^x}{\partial A^i \partial B^j}(A=0,B=0) 
\end{equation}
are linked to moments of the density  probability $P_X$.

The stability criteria $\Delta_c = \langle x_0^2 \rangle^2$ appears once
again as in section \ref{TrivialFixedPointSection_stab}. Below $\Delta < \Delta_c$ the uniform
fixed point $m=0$ is unstable and $m^t$ will converge towards another
fixed point different from $m = 0$. For $\Delta > \Delta_c =
\langle x_0^2 \rangle$ the uniform fixed point $m=0$ is stable. 
Using expansion (\ref{SecondOrderExpansion}) near $\Delta \simeq
\Delta_c=\langle x_0^2 \rangle$ we can write what is the other fixed
point next to $m_{\rm uniform}=0$, we get 
\begin{equation}
m_{\rm close-by} =  \frac{\Delta_c(\Delta_c - \Delta)}{\Delta_c^{3/2} - \frac{ \langle x_0^3
  \rangle^2}{2}} + O\left(\left(\Delta-
      \Delta_c\right)^2 \right) \, .
\label{Expansion_Q_NearTransition}
\end{equation}

By definition of the order parameters in the Bayes-optimal setting we
must have at a fixed point $m \ge 0$, therefore we distinguish two
cases 
\begin{itemize}
        \item If $\langle x_0^3 \rangle^2 < 2\langle x_0^2 \rangle^3$,
          eq.~(\ref{Expansion_Q_NearTransition}) is a stable fixed
          point in the region $\Delta <
          \Delta_c$. Eq.~(\ref{Expansion_Q_NearTransition}) then gives
          the expansion of this fixed point. This situation
          corresponds to Fig.~\ref{MinFreeEnergyFirstOrderExplanation} where the
          Rademacher-Bernoulli prior has zero 3rd moment.  This is the
          2nd order bifurcation at $\Delta_c$. 
        \item $\langle x_0^3 \rangle^2 > 2\langle x_0^2 \rangle^3$
          eq.~(\ref{Expansion_Q_NearTransition}) is an unstable (and
          hence irrelevant) fixed point in the region $\Delta > \Delta_c$. But also in this
          case there must be a stable fixed point for $\Delta <
          \Delta_c$, but this fixed point cannot have a small values
          of $m$. The only way a stable fixed point can appear in this
          case is by a discontinuous (1st order) transition at
          $\Delta_c$. This is the
          1st order bifurcation at $\Delta_c$. 
\end{itemize}



To summarize, we obtained a simple (sufficient) criteria for the
existence of a first order phase transition with $\Delta_{\rm Alg} =
\Delta_c$. Notably there is a 1st order phase transition when
\begin{equation}
\left\{
\begin{matrix}
\langle x_0 \rangle = 0
\\
\langle x_0^3 \rangle^2 > 2\langle x_0^2 \rangle^3
\end{matrix}
\right. \,.  \label{crit_1st}
\end{equation}
The more skewed the prior distribution is the easier the problem
is. Till a point where if the skewness of the signal is bigger than
$\sqrt{2}$ then a first order phenomena will appear in the system. In
this case the MSE achieved by the Low-RAMP algorithm becomes
discontinuously better than ${\rm MSE}_{\rm uniform} = \langle x_0^2 \rangle$

On the other hand when the criteria (\ref{crit_1st}) is not met, then
for $\Delta<\Delta_c$ the MSE achieved by the Low-RAMP algorithm is in
first approximation equal to 
\begin{equation}
{\rm MSE}(\Delta) = \langle x_0^2 \rangle - \frac{\Delta_c(\Delta_c - \Delta)}{\Delta_c^{3/2}- \frac{{\langle x_0^3 \rangle}^2}{2}} + O\left(\left(\Delta- \Delta_c\right)^2 \right) \,.
	\end{equation}
So the MSE obtained by an Low-RAMP algorithm is linear near the
transition in $\Delta - \Delta_c$. 

Interestingly spectral method also give an MSE linear in $\Delta - \Delta_c$.
As derived in section \ref{sec:PCA_analysis} the MSE one achieves using the eigenvectors of matrix $S$ is
\begin{equation}
{\rm MSE}_{\rm Spectral}(\Delta) = \langle x_0^2 \rangle  - \frac{\Delta_c - \Delta}{\sqrt{\Delta_c}} \,.
\end{equation}
From the coefficient of linearity in the error we observe that the
error achieved by PCA is always worse than the error achieved by Low-RAMP.  

Let us remind that uniform fixed points and 1st order phase
transition (at $\Delta_{\rm Alg} = \Delta_c$ or elsewhere) can exist
even if the criteria derived above are not met, there are sufficient,
not necessary conditions. Examples are included in subsequent sections. 

\section{Phase diagrams for Bayes-optimal low-rank matrix estimation}

From now on we restrict our analysis to the Bayes optimal
setting as defined in section \ref{sec:Bayes_optimal}, eq.~\eqref{Bayes_optimal}. The motivation
is to investigate performance of the Bayes-optimal and the Low-RAMP estimators for a
set of benchmark problems.  We investigate
phase diagrams stemming from the state evolution equations and from
the corresponding replica free energies summarized in section \ref{SectionSimplificationSE}.

\subsection{Examples of phase diagram}
\label{sec:Phase_Diagrams}
In this section we present example of phase diagrams for the symmetric
low-rank matrix estimation. 
The first three examples are for rank one, the last two examples are
for general rank.

\subsubsection{Spiked Bernoulli model}

The spiked Bernoulli model is defined by prior (\ref{Bernoulli_X}) with
density $\rho$ of ones, and $1-\rho$ of zeros. This prior has a
positive mean and consequently the state evolution does not have the
uniform fixed point. This is a problem where one tries to recover a submatrix of size $\rho N \times
\rho N$ with mean of elements equal to 1, in a $N \times N$ matrix of
lower mean. Using the Bayes-optimal state evolution (\ref{SE_Equation_XX_M}) one gets,
\begin{eqnarray}
m^t &=& f_{\rm Bernoulli}^{\rm SE}\(( \frac{m^t}{\Delta} \))\, ,
\\
f_{\rm Bernoulli}^{\rm SE}(x) &=& \rho \mathbb{E}_W\left[ \frac{\rho}{\rho +
                      (1-\rho) \exp\left( \frac{-x}{2} + W\sqrt{x}
                      \right)} \right] \, ,
\label{SE_Bernoulli}
\end{eqnarray}
where $W$ is (here and from now on) a Gaussian variables of zero mean
and unit variance. 

\begin{figure}
\includegraphics[scale=0.8]{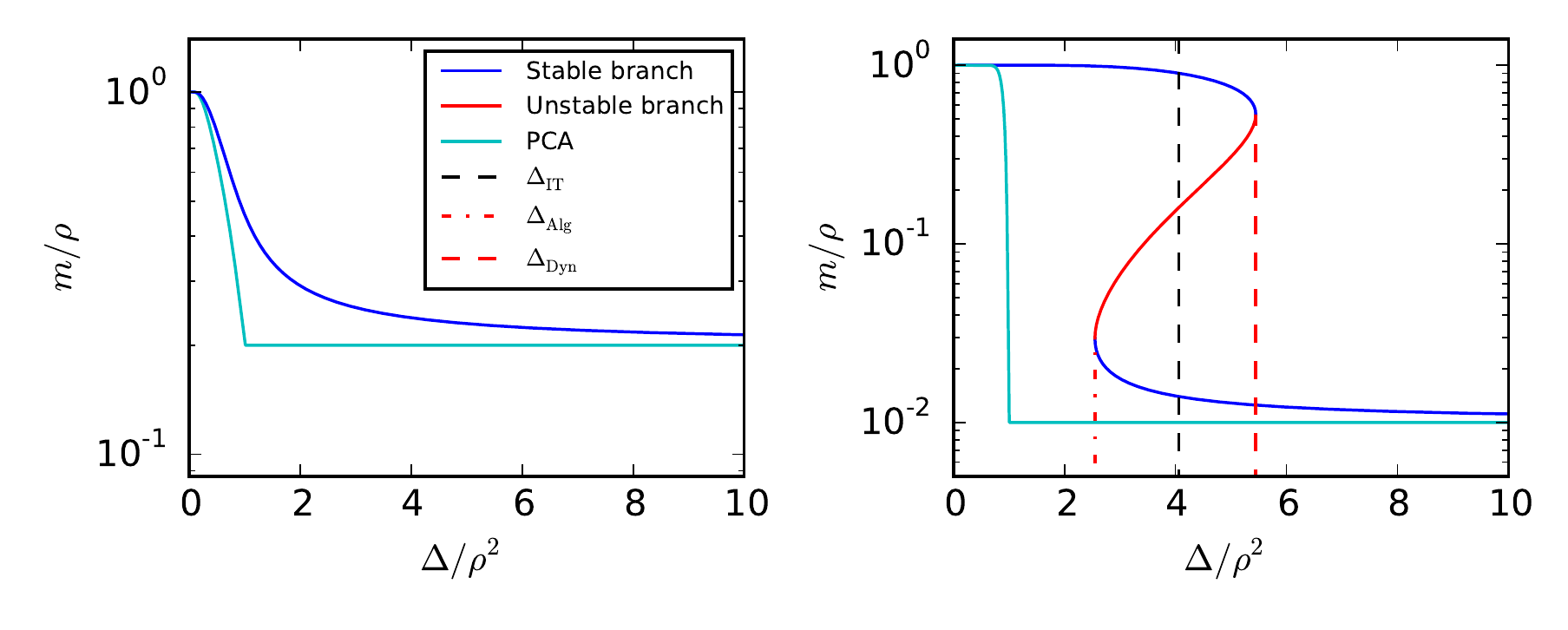}
\caption{Plots of the fixed points of the state evolution for the
  spiked Bernoulli model, MSE being given by ${\rm MSE} = \rho -m$ (\ref{MSE_X_Nish}). The performance of PCA as analyzed in section
  \ref{sec:PCA_analysis} is plotted for comparison. These two  plots are made for $\rho = 0.2$ (left) and  $\rho = 0.01$ (right).}
\label{MSE_Bernoulli}
\end{figure}

\begin{figure}
   \centering
\begin{minipage}{.5\textwidth}
  \centering
  \includegraphics[width=1\linewidth]{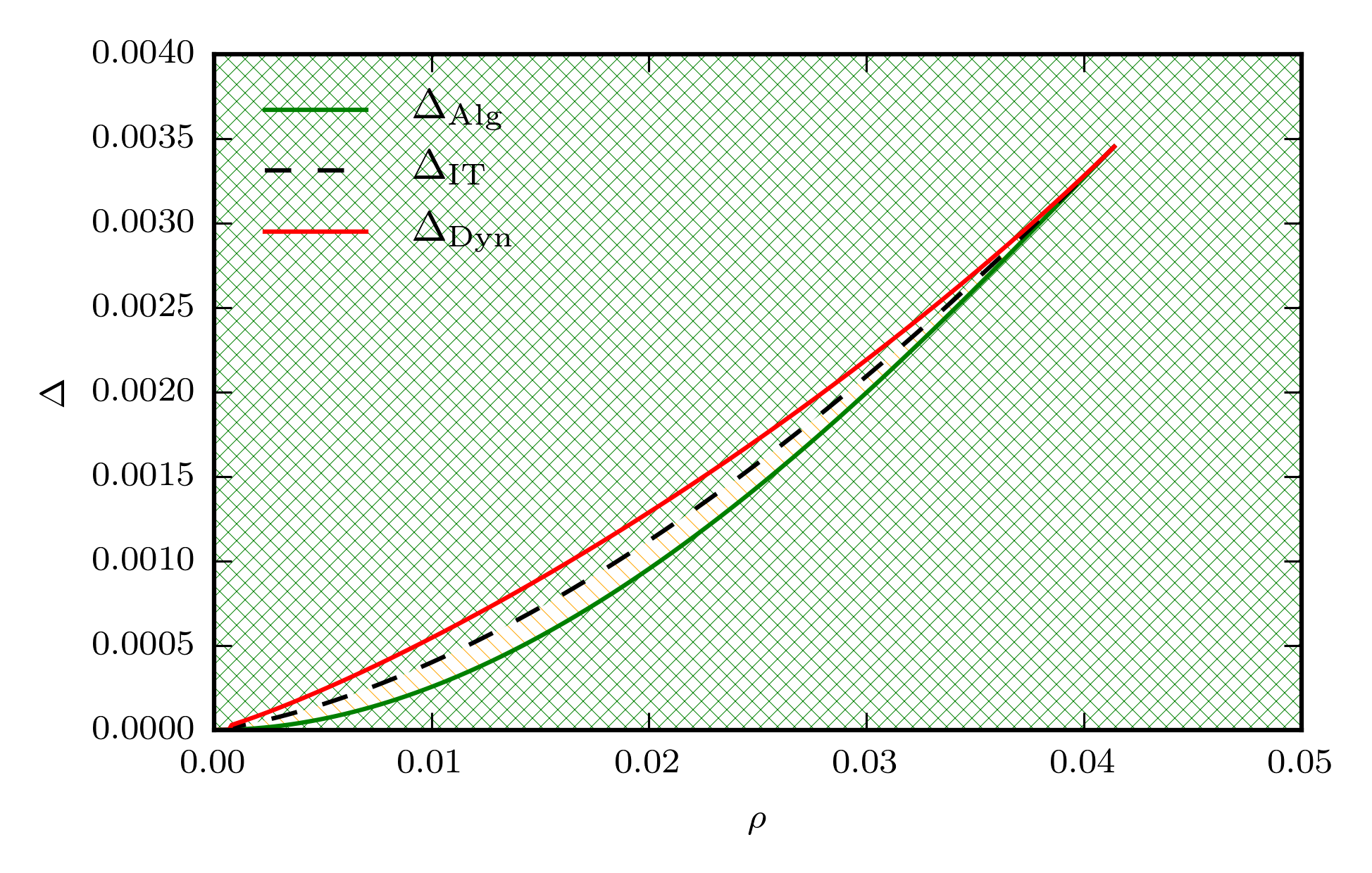}
\end{minipage}%
\begin{minipage}{.5\textwidth}
  \centering
  \includegraphics[width=1\linewidth]{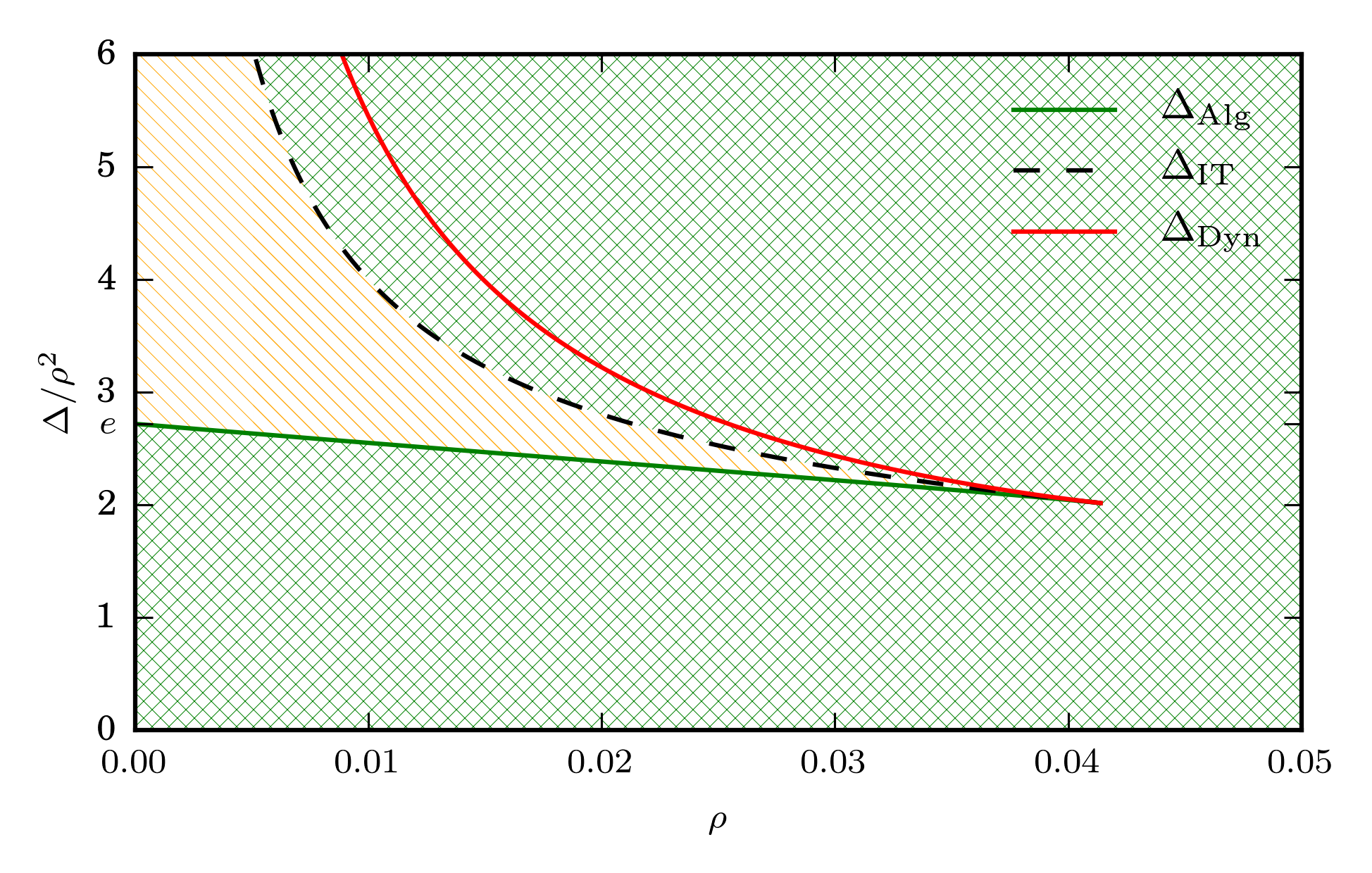}
  \end{minipage}
\caption{The phase diagram of the spiked Bernoulli model
  (\ref{Bernoulli_X}) as a function of the density $\rho$ and effective
  noise $\Delta$ (left) or $\Delta/\rho^2$ (right). There is no phase
  transition in the system for $\rho > 0.04139$ and a 1st order phase
  transition for $\rho < 0.04139$. The lower green curve is the
  algorithmic spinodal $\Delta_{\rm Alg}$ curve, that converges to
  $\Delta_{\rm Alg} =_{\rho \to 0} e \rho^2$. The  dashed  black
  line is the information theoretic threshold~$\Delta_{\rm IT}$.  The
  upper red curve is the dynamical spinodal~$\Delta_{\rm Dyn}$. The
  orange hashed zone is the hard region in which Low-RAMP does not
  reach the Bayes-optimal error. In the rest of the phase diagram
  (green hashed) the Low-RAMP provides in the large size limit the
  Bayes-optimal error. Note that this is exactly the same phase
  diagram as presented in \cite{montanari2015finding} (Fig. 5) for the
problem of finding one dense subgraph, this is thanks to output
channel universality and the fact that large degree sparse graphs have
upon rescaling the same phase diagram as dense graphs.}
\label{SpikedBernoulliPhaseDiagram}
\end{figure}

Depending on the value of $\rho$ there are two kinds of behaviour of
the fixed points of these equations as a function of the effective
noise $\Delta$. We plot the two cases in Fig.~\ref{MSE_Bernoulli}. For
larger values of $\rho$ there is a unique fixed point corresponding to
the MMSE that is asymptotically achieved by the Low-RAMP algorithm,
this is the regime in which the proof of \cite{DeshpandeM14}
applies. For small enough values of $\rho$ we do observe a region of
$\Delta_{\rm Alg}<\Delta< \Delta_{\rm Dyn}$ where there are 3 fixed
points, two stable and one unstable. The replica free energy
associated to the these fixed points crosses at $\Delta_{\rm IT}$ so
that the higher fixed point in the relevant one at $\Delta<\Delta_{\rm
  IT} $, and the lower fixed point at $\Delta>\Delta_{\rm
  IT} $. These phase transitions were defined in section
\ref{SubSectionFirstOrderGeneral}.  In
Fig.~\ref{SpikedBernoulliPhaseDiagram} we plot the phase transitions
$\Delta_{\rm Alg}$, $\Delta_{\rm IT}$ and $\Delta_{\rm Dyn}$ as a
function of $\rho$. The $y$-axes in the left pannel is simply the
effective noise parameter $\Delta$, on the right pannel the same data
are plotted with $\Delta/\rho^2$ on the $y$-axes. We observe that the
$\Delta_{\rm Alg} =_{\rho \to 0} e \rho^2 $, with $e$ being the Euler
number, this is the same asymptotic behaviour as obtained previously
in \cite{montanari2015finding} (Fig. 5). 

\subsubsection{Rademacher-Bernoulli and Gauss-Bernoulli}

Next we analyze the spiked Rademacher-Bernoulli and Gauss-Bernoulli
priors defined by \eqref{P_X_Rademacher}  and
\eqref{P_X_GaussBernoulli}. 
The first thing we notice is that both these distribution have zero
mean and variance $\rho$. This means according to
(\ref{TrivialFixedPointSection}) that there is a uniform fixed point of the SE
equations that is stable for $\Delta > \rho^2$. The skewness of both
these distribution is 0, which means that at $\Delta_c$ there is no
discontinuity of the $\rm MSE$ (\ref{FirstOrderCriteriaSubSection}). The SE equations for these models are
\begin{eqnarray}
m^{t+1} &=& f_{\rm Rademacher-Bernoulli}^{\rm
            SE}\left(\frac{m^t}{\Delta} \right) \, ,\label{SE_Rademacher_Bernoulli_SE}
\\
f_{\rm Rademacher-Bernoulli}^{\rm SE}(x) &=& \rho \mathbb{E}_{W}\left[\tanh\left(x	+ W \sqrt{x}  \right) \frac{\rho }{(1-\rho) \frac{\exp(x/2)}{\cosh\left(x + W \sqrt{x} \right)}  + \rho} \right],
\label{SE_Rademacher_Bernoulli}
\\
m^{t+1} &=& f_{\rm Gaussian-Bernoulli}^{\rm
            SE}\left(\frac{m^t}{\Delta} \right)\, ,
\\
f_{\rm Gauss-Bernoulli}^{\rm SE}(x)/\rho &=& \frac{x}{1 + x} \mathbb{E}_{W} \left[ W^2 \hat{\rho}(x, W \sqrt{x^2 + x})\right],
\label{SE_Gaussian_Bernoulli}
\end{eqnarray}
where $\hat{\rho}$ is
\begin{equation}
{\hat \rho}(a,b^2) = \frac{\rho}{(1-\rho)
  \exp\left(\frac{-b^2}{2(1+a)}\right)(1+a)^\frac{r}{2}+ \rho} \,. 
\end{equation}
Both these models have similar phase diagram.

\begin{figure}[!ht]
\includegraphics[scale=0.8]{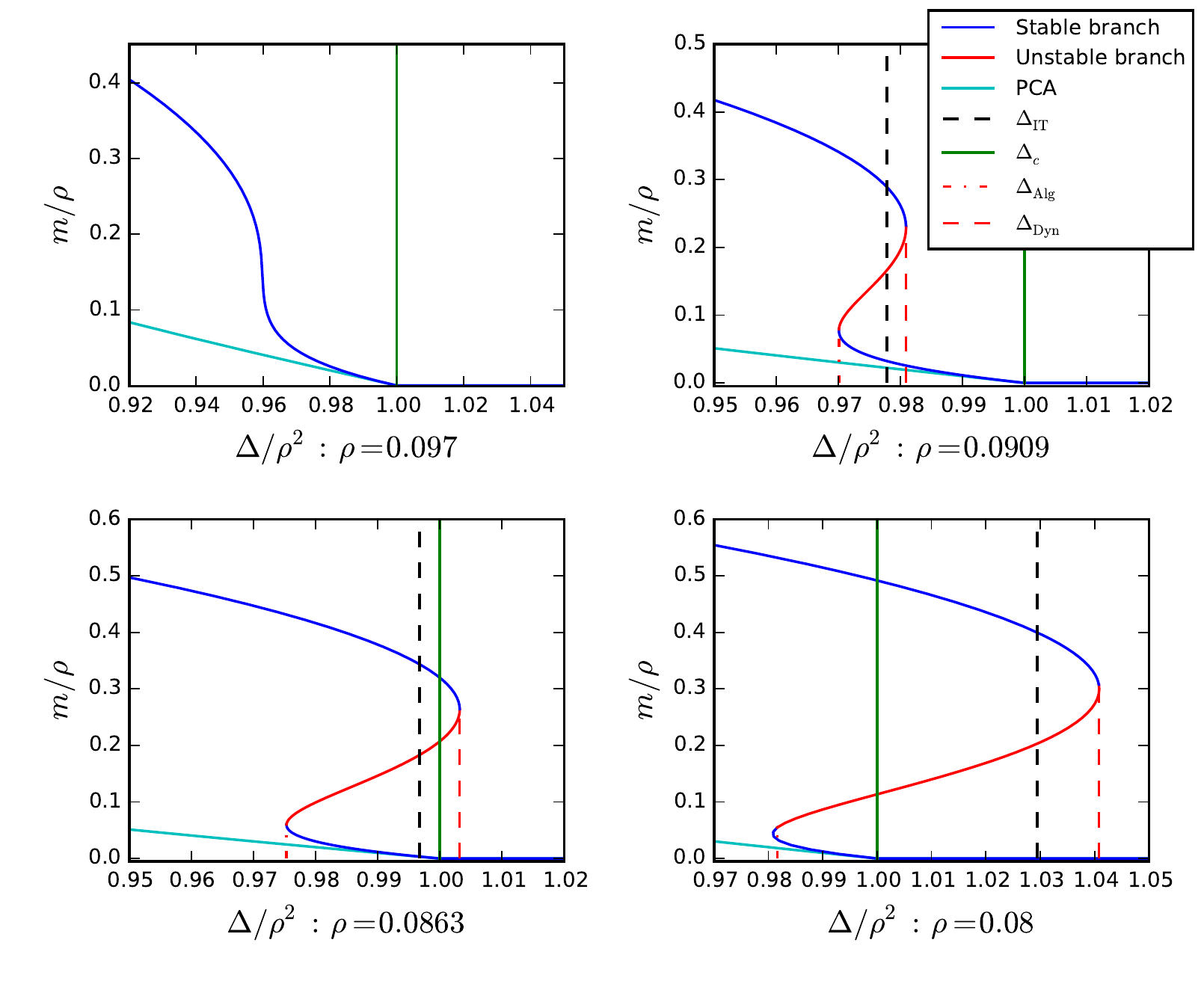}
\caption{We plot all the fixed points of the state evolution equations
  for the spiked Rademacher-Bernoulli model in the Bayes optimal setting as
  a function of ${\Delta/\rho^2}$ for four representative values of
  $\rho$. The blue curves are stable fixed points, red are the
  unstable ones. MSE is given by (\ref{MSE_X_Nish}). The vertical lines mark the two spinodal transition
  (dynamical $\Delta_{\rm Dyn}$ and algorithmic $\Delta_{\rm Alg}$, red dashed), the information theoretic
  transition $\Delta_{\rm IT}$. The vertical green full line mark the stability point
of the uniform fixed point $\Delta_c = \rho^2$. We remind that the error ${\rm MSE} = \rho - m$  achieved by the Bayes optimal
estimator corresponds to the upper branch for $\Delta<
\Delta_{\rm IT}$, and to the lower branch for $\Delta>
\Delta_{\rm IT}$. Error achieved by the Low-RAMP algorithm always
correspond to the lower branch (larger error). Note that in the three
panels where multiple fixed points exists, the only element that
changes is the position of the (spectral) stability threshold
$\Delta_c$ with respect to the other thresholds. 
}
\label{AllFirstOrderCases}
\end{figure}

We first illustrate in Fig.~\ref{AllFirstOrderCases} the different
types of phase transition that we observe for the spiked Rademacher-Bernoulli
model as the density $\rho$ is varied. We plot all the fixed points of equation
(\ref{SE_Rademacher_Bernoulli}) for several values of $\rho$ as a
function of the effective noise $\Delta/\rho^2$.
The four observed case are the following 
\begin{itemize}
\item $\rho = 0.097$ example: For $\rho$ large enough (in the present
  case $\rho>\rho_{\rm tri}= 0.0964$) whatever the $\Delta$ there is only one stable fixed point.
\item For small enough $\rho < \rho_{\rm tri}= 0.0964$ three different
  fixed points exist in a range of $\Delta_{\rm
    Alg}(\rho)<\Delta<\Delta_{\rm Dyn}(\rho)$, where the thresholds $\Delta_{\rm
    Alg}$, and $\Delta_{\rm Dyn}$ are defined by the limits of
  existence of the three fixed points. The information theoretic
  threshold where the free energy corresponding to the two stable
  fixed points crosses is $\Delta_{\rm IT}$. Depending on the values
  of $\rho$ we observed 3 possible scenarios of how $\Delta_c=
\rho^2$ is placed w.r.t. the other thresholds. 
\begin{itemize}
\item $\rho = 0.0909$ example where
  $\Delta_{\rm Dyn}<\Delta_c$. 
\item $\rho = 0.0863$ example where $\Delta_{\rm IT} < \Delta_c <
  \Delta_{\rm Dyn}$.
\item $\rho = 0.08$ examples where $\Delta_{\rm Alg} < \Delta_c <
  \Delta_{\rm IT}$. 
\end{itemize}
\end{itemize}

Finally Fig.~\ref{PhaseDiagram_Gauss_Rademacher_Bernoulli_zoom}
presents the four thresholds $\Delta_{\rm Dyn}$, $\Delta_{\rm IT}$, $\Delta_{\rm Alg}$ and
$\Delta_{c}$ as a function of the density $\rho$.

\begin{figure}
   \centering
\begin{minipage}{.5\textwidth}
  \centering
  \includegraphics[width=1\linewidth]{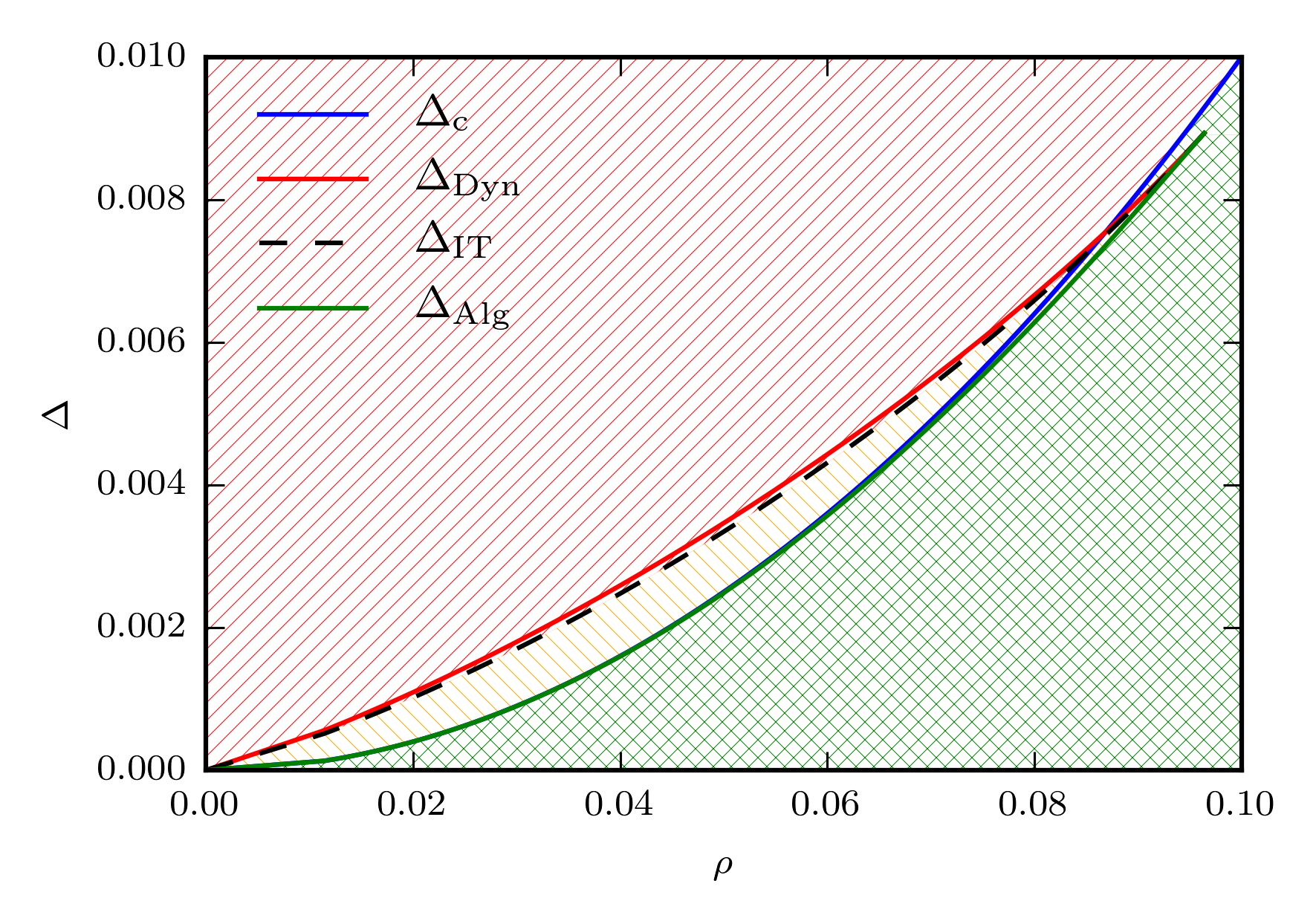}
\end{minipage}%
\begin{minipage}{.5\textwidth}
  \centering
  \includegraphics[width=1\linewidth]{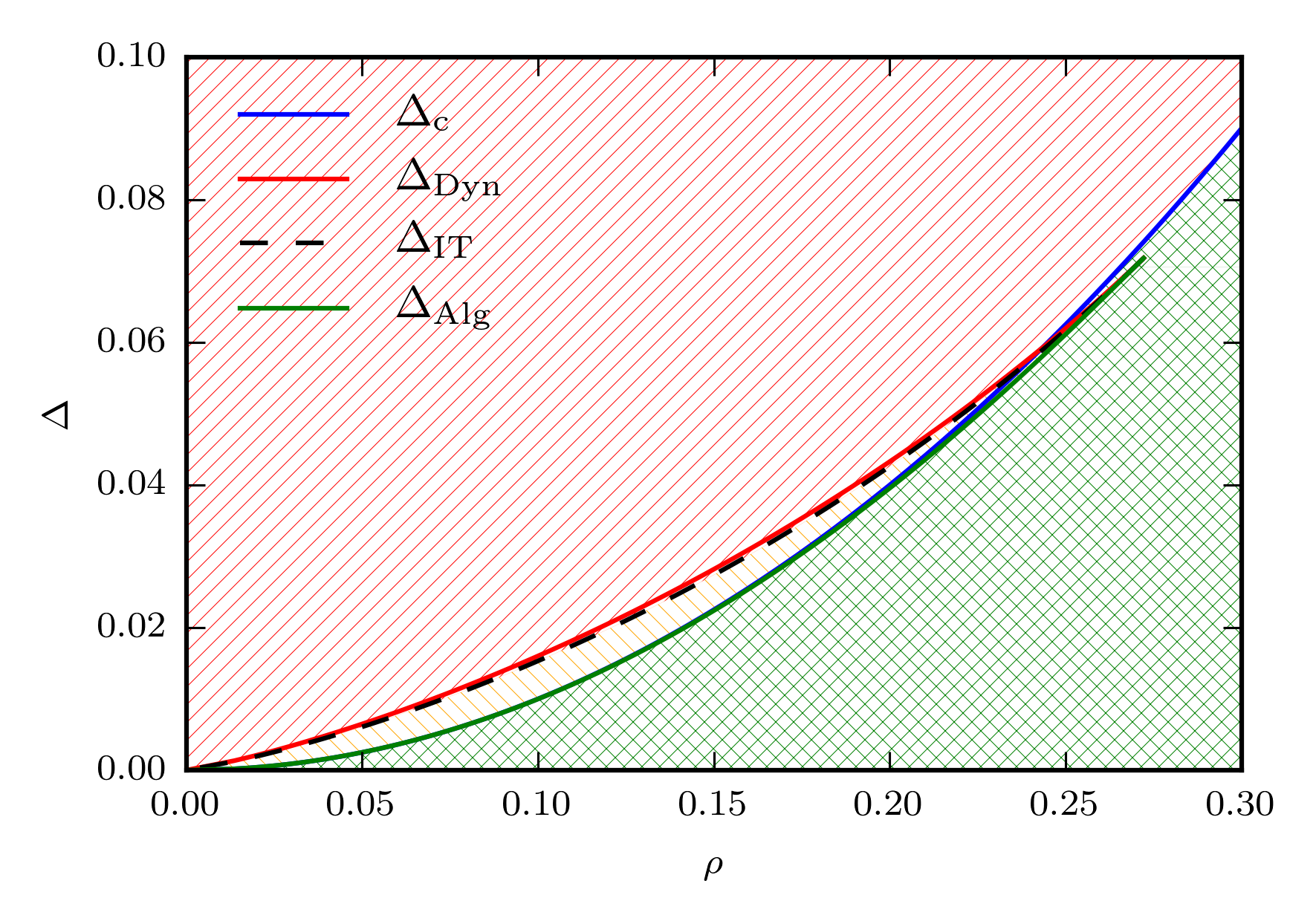}
  \end{minipage}
\caption{Phase diagram of the spiked Rademacher-Bernoulli (left hand
  side pannel) and spiked
  Gauss-Bernoulli (right hand side pannel) models. We plot $\Delta$ as a function of
  $\rho$. The local stability threshold of the uniform fixed point
  $\Delta_c/\rho^2=1$ is in blue. The algorithmic spinodal $\Delta_{\rm Alg}$ (green), the
  dynamical spinodal $\Delta_{\rm Dyn}$ (red) and the information theoretic
  transition $\Delta_{\rm IT}$ (black dashed) all join into a tri-critical point
  located at 
  $(\Delta_{\rm tri}=0.008935, \Delta_{\rm tri}/\rho_{\rm tri}^2 = 0.9612, \rho_{\rm tri}=0.09641)$ for the Rademacher-Bernoulli model (left pannel), and at
  $(\Delta_{\rm tri}=0.07182, \Delta_{\rm tri}/\rho_{\rm tri}^2 = 0.9693, \rho_{\rm tri}=0.2722)$  for the Gauss-Bernoulli model (right pannel). 
The hash materializes the different phases. The easy phase where
the Low-RAMP algorithm is Bayes-optimal and achieves better error than random
guessing is hashed in green crossed lines, the hard phase
where Low-RAMP is suboptimal is hashed in yellow $\setminus \setminus$,
and the impossible phase where even the best achievable error is as bad as
random guessing is hashed in red  $//$.  
}
\label{PhaseDiagram_Gauss_Rademacher_Bernoulli} 
\end{figure}

\begin{figure}
   \centering
\begin{minipage}{.5\textwidth}
  \centering
  \includegraphics[width=1\linewidth]{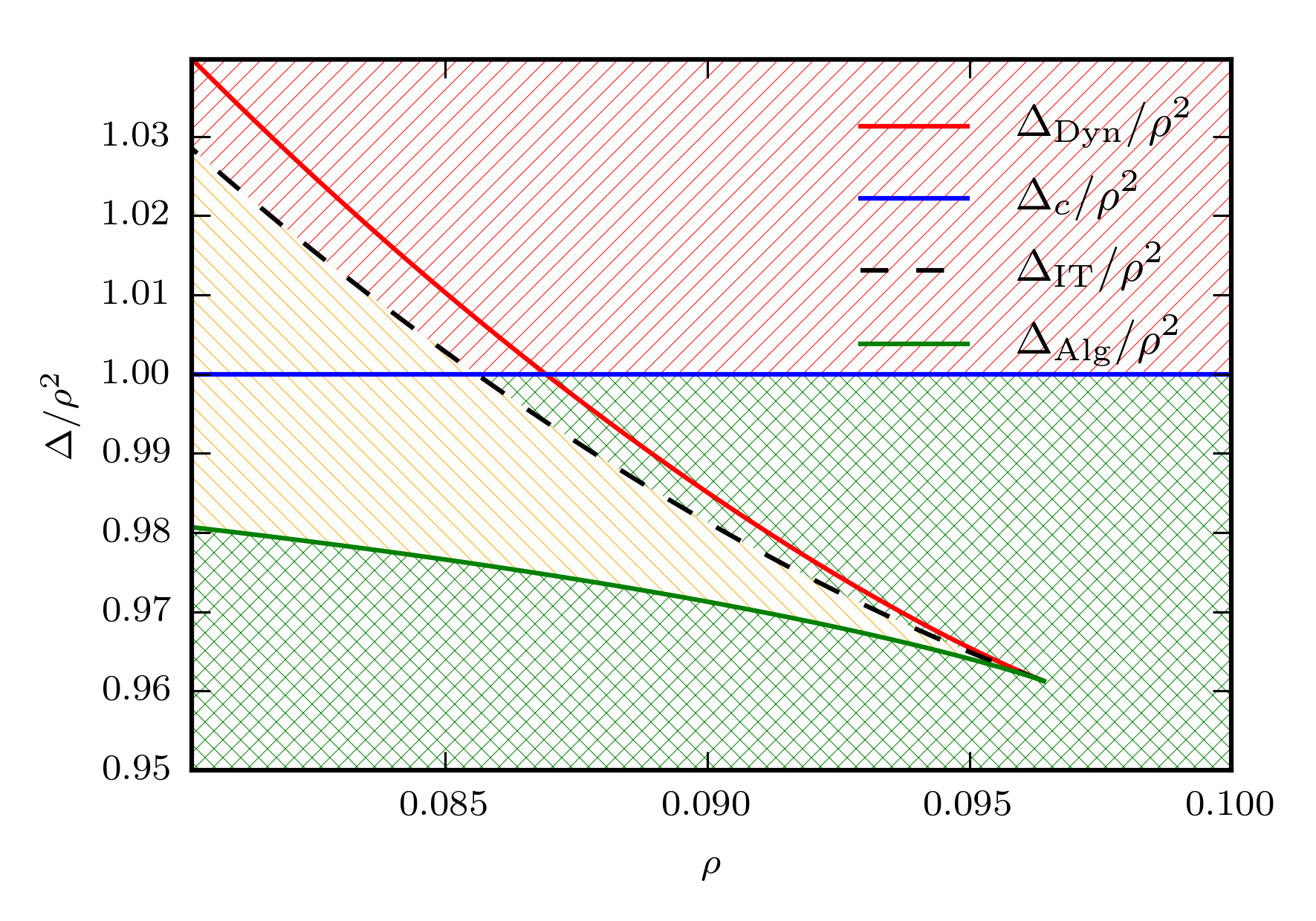}
\end{minipage}%
\begin{minipage}{.5\textwidth}
  \centering
  \includegraphics[width=1\linewidth]{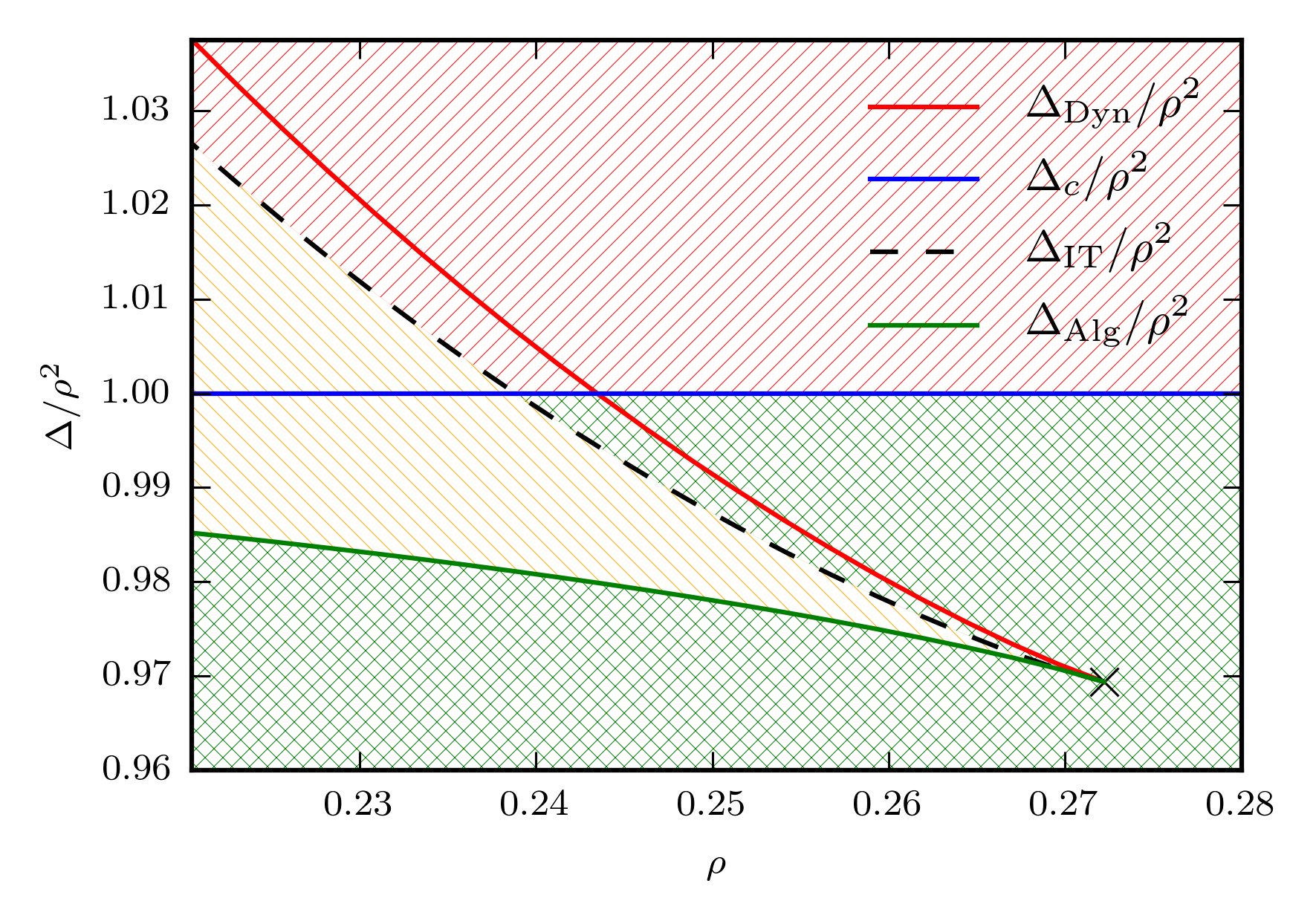}
  \end{minipage}
\caption{The same plot as in
  Fig.~\ref{PhaseDiagram_Gauss_Rademacher_Bernoulli} zoomed into the
  region of the tri-critical point with $y$-axes rescaled by $\rho^2$. 
  The spiked Rademacher-Bernoulli mode on left hand
  side pannel, the spiked Gauss-Bernoulli model on the right hand side.
}
\label{PhaseDiagram_Gauss_Rademacher_Bernoulli_zoom} 
\end{figure}

\begin{figure}[!ht]
\includegraphics[scale=1]{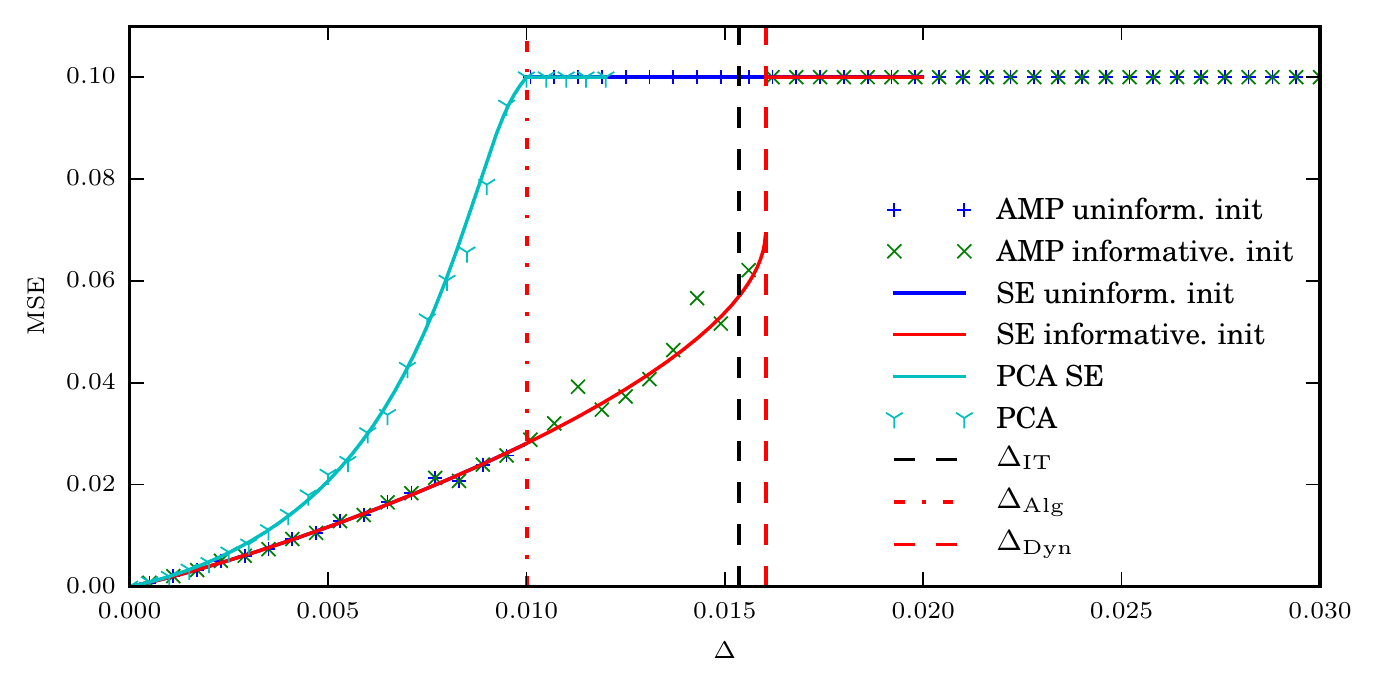}
\caption{Comparison between the state evolution and the fixed point of
  the Low-RAMP algorithm, for the spiked Gauss-Bernoulli model of
  sparse PCA with rank one and density $\rho=0.1$. The phase transitions stemming from state
  evolution are  $\Delta_{\rm Alg} \approx \Delta_{\rm c} = 0.01,
  \Delta_{\rm IT} = 0.0153, \Delta_{\rm Dyn} = 0.0161$. 
The points are the fixed points of the Low-RAMP algorithm run on one
typical instance of the problem of size 
$N = 20000$. Blue pluses is the MSE reached from an uninformative
initialization of the algorithm. Green crosses is the MSE reached from
the informative initialization of the algorithm. 
}
\label{MSE_Sparse_PCA}
\end{figure}

In Fig.~\ref{MSE_Sparse_PCA} we present for completness the comparison
between the fixed points on the state evolution and the fixed points
of the Low-RAMP algorithm for the Gauss-Bernoulli model, with rank
one, Bayes optimal case. The experiment is done on one random instance
of size $N=2\times 10^4$ and we see the agreement is very good, finite
size effect are not very considerable. The data are for the
Gauss-Bernoulli model at $\rho=0.1$, that is in a region where
$\Delta_{\rm Alg}$ is so close to $\Delta_c$ that in this figure the
difference is unnoticeable. 

We also compare to the MSE
reached by the PCA spectral algorithm and from its analysis
eq.~(\ref{PCA_optimal_reconstruction}). We see that whereas both Low-RAMP and PCA
work better than random guesses below $\Delta_c$, the MSE reached by
Low-RAMP is considerably smaller.

\subsubsection{Two balanced groups}
\label{2-groups-balanced}

The next example of phase diagram we present is for community detection with two
balanced (i.e. one group is smaller $\rho N$, but both have the same average degree) groups as defined in eqs.~(\ref{ConnectivityMatrix_Balanced_Intro}-\ref{Proba_Cliques}). This is an example
of a system where the bifurcation at $\Delta_c$ is of a first
order, with $\Delta_c= \Delta_{\rm Alg}$. In this problem the prior
given by eq. (\ref{Proba_Cliques}), with $\langle x_0 \rangle =0$, $\langle (x_0)^2 \rangle=1$.
The output channel is of the stochastic block model type
eq.~(\ref{SBM_out1}-\ref{SBM_out2}),
leading to effective noise parameter 
\begin{equation}
\Delta = \frac{p_{\rm out}(1-p_{\rm out})}{\mu^2} \, , \label{delta_2groups}
\end{equation}
where $\mu$ and $p_{\rm out}$ are the parameters from
(\ref{SBM_out1}-\ref{SBM_out2}).
Therefore the uniform fixed point becomes unstable when $\Delta < 1$.

Using eqs. (\ref{SE_Equation_XX_M}) and (\ref{Proba_Cliques}) we get
for the state evolution for community detection with two balanced
groups 
\begin{equation}
m^{t+1} = f_{\rm Two Balanced}^{\rm SE}\left(\frac{m^t}{\Delta}\right),\forall t ,m^t \in [0;1] \,, \label{DensityEvolution_Cliques}
\end{equation}
where
\begin{equation}
f_{\rm Two Balanced}^{\rm SE}(x) = \int\limits_{-\infty}^{+\infty}\frac{e^{\frac{-u^2}{2}}}{\sqrt{2\pi}}\frac{2\rho(1-\rho)\sinh\left(\frac{x}{2\rho(1-\rho)} + u \sqrt{\frac{x}{\rho(1-\rho)}}\right)}{1+2\rho(1-\rho)\left( \cosh\left(\frac{x}{2\rho(1-\rho)} + u \sqrt{\frac{x}{\rho(1-\rho)}} \right)-1 \right)} \,. \label{F_Rho_Cliques}
\end{equation}

To investigate whether $\Delta_c=1$ is a 1st or 2nd order bifurcation
we compute the second order expansion of the state evolution
equations as we have done in (\ref{SecondOrderExpansion}).
We find that the expansion of the state evolution up to second order
for the two balanced communities is 
\begin{equation}
m^{t+1} = f_{\rho}\left(\frac{m^t}{\Delta}\right) = \frac{m^t}{\Delta}
+ \left(\frac{m^t}{\Delta}\right)^2  \frac{1 - 6\rho(1-\rho)}{2
  \rho(1-\rho)}\, .   \label{frho_2groups}
\end{equation}

\begin{figure}
\includegraphics[scale=0.6]{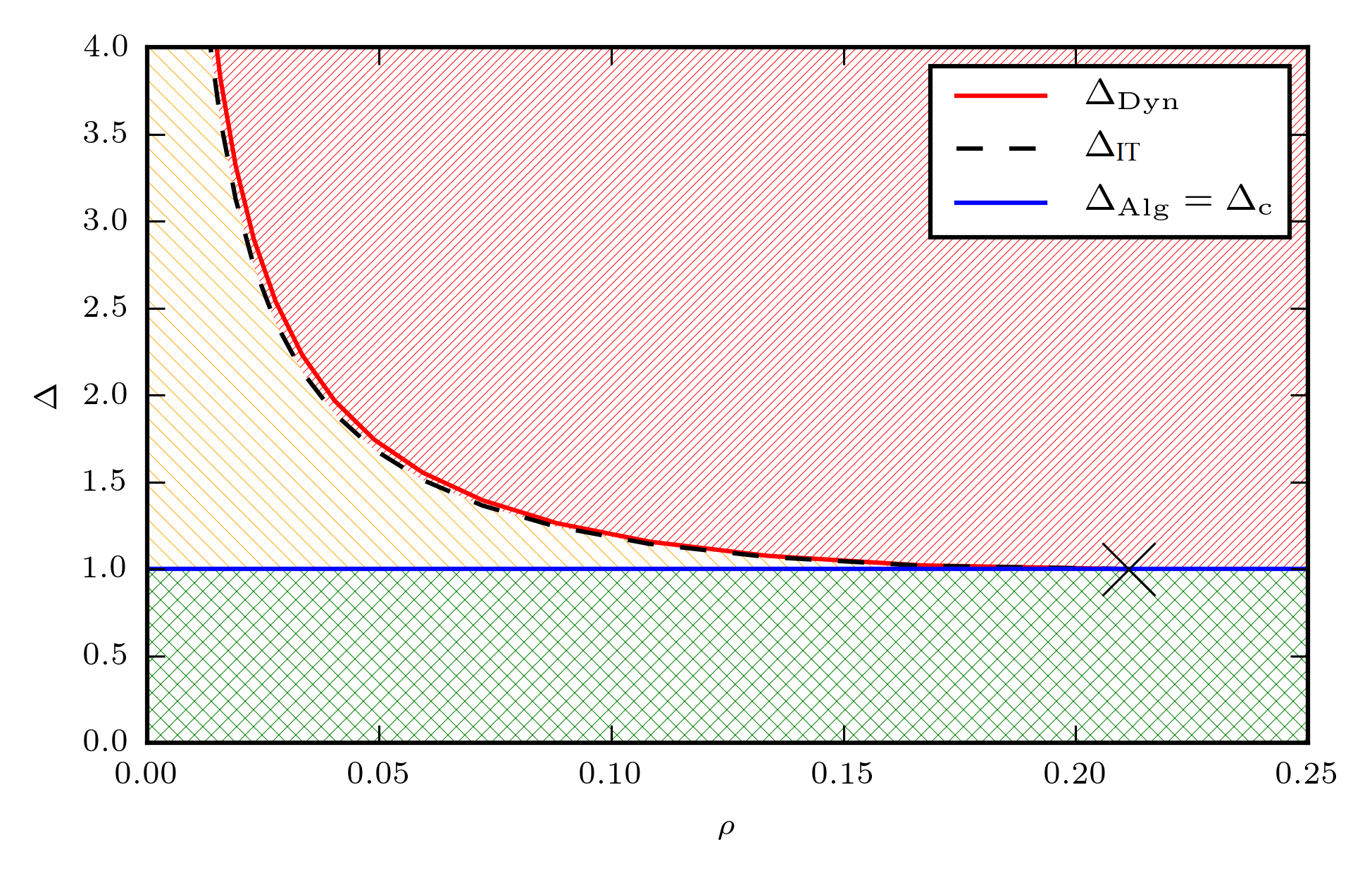}
\caption{We plot here $\Delta_{\rm Alg} = \Delta_c$, $\Delta_{\rm IT}$
  and $\Delta_{\rm Dyn}$ as a function of $\rho$ for the community
  detection with two balanced groups. All the curves merge at $\rho =
  \frac{1}{2} - \frac{1}{\sqrt{12}}$. The hashed regions have the same
meaning as in previous figures, red is the impossible inference phase,
green is easy and yellow is hard inference.}
\label{PhaseDiagram_Balanced}
\end{figure}

In a similar fashion as in section
(\ref{FirstOrderCriteriaSubSection}) it is the sign of the second
order terms that decides between 1st or 2nd order bifurcation at
$\Delta_c=1$. We find that if $\rho(1-\rho) < 1/6$ then the second
order derivative of (\ref{F_Rho_Cliques}) is positive leading to a
jump in MSE when $\Delta$ crosses the value $\Delta=1$ which means
that there will be first order phase transition for all 
\begin{equation}
\rho \in \left[0;\frac{1}{2} - \frac{1}{\sqrt{12}} \approx 0.21 \right] \cup
\left[\frac{1}{2}  + \frac{1}{\sqrt{12}}\approx 0.89 ; 1\right] \,. \label{eq:int}
\end{equation}
It turns out that for the two balanced groups this criteria is both
sufficient and necessary. Out of the interval (\ref{eq:int}) the phase
transition at $\Delta_c$ is of second order, with no discontinuities. 
Defining the phase transitions $\Delta_{\rm Alg}$, $\Delta_{\rm IT}$,
$\Delta_{\rm Dyn}$ as before in section \ref{SubSectionFirstOrderGeneral}, we have $\Delta_{\rm
  Alg}=\Delta_c$ and we plot the three phase transitions for community
detection for two balanced groups in the phase diagram Fig.~\ref{PhaseDiagram_Balanced}.

\subsubsection{Jointly-sparse PCA generic rank}

In this section we discuss analysis of the Gauss-Bernoulli jointly
sparse PCA as defined by the prior distribution
(\ref{P_X_JointedGaussBernoulli}) for a generic rank $r$.  
This just means that each vector $x_i \in {\mathbb R}^r$ is either $\vec{0}$ with
probability $1-\rho$ or is taken from Gaussian density probability of
mean zero and covariance matrix $I_r$ with probability $\rho$. 
The prior distribution (\ref{P_X_JointedGaussBernoulli})  has a zero
mean, therefore the uniform fixed point exist and according to
(\ref{Trivial_XX_Stability}) is stable down to
$\Delta_c = \rho^2$. 

In order to deal with the $r$-dimensional integrals and $r \times r$
dimensional order parameter $M^t$ we notice that there is a rotational
${\rm SO}(r)$ symmetry in the problem, which in the Bayes optimal setting implies
that 
\begin{equation}
M^t = m^t I_r,
\end{equation}
with $m^t$ being a scalar parameter. The problem is hence greatly
simplified, one can then treat the $r$ dimensional integral in (\ref{SE_Nish_XX}). After integration by parts and integration on the sphere one gets
\begin{eqnarray}
\label{JointPCA_DE}
m^{t+1} &=& f_{\rm Joint-GB}^{\rm
            SE}\left(\frac{m^t}{\Delta}\right)\, ,
\\
f_{\rm Joint-GB}^{\rm SE}(x) &=&  \frac{\rho x}{1 + x} \int {\rm d} u \frac{1}{(2
  \pi)^{\frac{r}{2}}} \exp \left(\frac{-u^2}{2}\right) S_r u^{r-1}
\left\{1 + \frac{x u^2 \left[1 - \hat \rho(x,(x^2 +  x)
      u^2)\right]}{r}\right\} \hat \rho(x,(x^2 + x) u^2) \,,
\end{eqnarray}
where $S_r$ is the surface of a unit sphere in $r$ dimensions and where $\hat \rho$ is the posterior probability that a vector is equal to
$\vec{0}$.
\begin{equation}
\label{Estimate_rho}
{\hat \rho}(a,b^2) = \frac{\rho}{(1-\rho) \exp\left(\frac{-b^2}{2(1+a)}\right)(1+a)^\frac{r}{2}+ \rho} \,.
\end{equation}
An expansion of (\ref{JointPCA_DE}) around $m^t=0$ yields
\begin{equation}	
m^{t+1} = \frac{\rho^2 m^t}{\Delta} - \rho^3 \left( \frac{m^t}{\Delta} \right)^2 + O\left( (m^t)^3\right)
\label{ExpansionAround_0_Joint}\,.
\end{equation}
Analogously to the conclusions we reached when studying expansion
(\ref{SecondOrderExpansion}), we conclude that since the second term is
negative there is always a 2nd order bifurcation at $\Delta_c=\rho^2$
with a stable fixed point for $\Delta< \Delta_c$ that stays close to
the uniform fixed point. At the same time this close-by fixed point
typically exists only in a very small interval of $(\Delta_{\rm Alg},
\Delta_c)$, similarly as in Fig.~\ref{MSE_Sparse_PCA}.

\subsubsection{Community detection with symmetric groups}

In this section we discuss the phase diagram of the symmetric communities
detection model as defined in section
\ref{ParagraphCommunityDetection}. The corresponding prior
distribution is \eqref{SBM_P_X} which leads to the function $f_{\rm in}^x$
\begin{equation}
\forall k \in [1:r] ,{f_{\rm in}^x}(A,B)_k = \frac{\exp\left( B_k - A_{kk}/2 \right)}{\sum\limits_{k'=1 \cdots r}\exp\left( B_{k'} - A_{k'k'}/2 \right)}\,.
\label{F_in_Community}
\end{equation}
The corresponding output is given by (\ref{SBM_out1}-\ref{SBM_out2}),
corresponding to the effective noise 
\begin{equation}
\Delta = \frac{p_{\rm out}(1-p_{\rm out})}{\mu^2}  =  \frac{ p_{\rm
    out}(1-p_{\rm out})}{N (p_{\rm in} - p_{\rm out}  )^2} \,.  \label{Delta_SBM}
\end{equation}

Once again we study the phase diagram by analyzing the state evolution
(\ref{SE_Nish_XX}). We can verify that the following form of the order
parameter in invariant under iterations of the Bayes-optimal state evolution 
\begin{equation}
M^t = b^t \frac{I_r}{r} + \frac{(1-b^t)J}{r^2}  \, ,
\end{equation}
where $J$ is the matrix filled with 1.
The order parameter at time $t+1$ will be of the same form with a new $b^{t+1}$.
\begin{itemize}
\item Having $b^t = 0$ is equivalent to having all the variables $x_i$
  saying that they have an equal probability to be in every
  community. This corresponds to initializing the estimators of the
  algorithm to be $\hat x_i^{t=0} = \left(\frac{1}{r} , \cdots ,\frac{1}{r}\right)$ 
\item Having $b^t = 1$ means that the communities have been perfectly
  reconstructed. This corresponds to initializing the algorithm in the
  planted solution.
\end{itemize}
Using (\ref{SE_Nish_XX}) the state evolution equations for $b^t$ can be written.
\begin{equation}
b^{t+1} = {\cal M}_r\left( \frac{b^t}{\Delta} \right) \, , 
\label{SE_b_Community} 
\end{equation}
where
\begin{equation}
{\cal M}_r\left( x \right)= \frac{1}{r - 1}\left[ r \int
\frac{\exp \left( \frac{x}{r} + u_1 \sqrt{\frac{x}{r}} \right)}
{
\exp \left( \frac{x}{r} + u_1 \sqrt{\frac{x}{r}} \right) + \sum\limits_{i = 2}^r \exp \left( u_i \sqrt{\frac{x}{r}} \right)
}
\prod\limits_{i=1}^r{\rm d} u_i \frac{\exp \left(\frac{-u_i^2}{2} \right)}{\sqrt{2\pi}}  -1  \right] \,. \label{Function_Cal_M}
\end{equation}
This can be proven by computing $M^{t+1}_{11}$ using
(\ref{SE_Nish_XX}) and (\ref{F_in_Community}) which yields
\begin{equation}
b^{t+1}\left(\frac{1}{r} - \frac{1}{r^2} \right)= \frac{1}{r}\left[\int
\frac{\exp \left( \frac{x}{r} +u_1 \sqrt{\frac{x}{r}} + u_0\sqrt{\frac{1-b^t}{r^2 \Delta}} \right)}
{
\exp \left( \frac{x}{r} + u_1 \sqrt{\frac{x}{r}}  + u_0\sqrt{\frac{1-b^t}{r^2 \Delta}}  \right) + \sum\limits_{i = 2}^r \exp \left( u_i \sqrt{\frac{x}{r}}  + u_0\sqrt{\frac{1-b^t}{r^2 \Delta}} \right)
}
\prod\limits_{i=0}^r{\rm d} u_i \frac{\exp \left(\frac{-u_i^2}{2} \right)}{\sqrt{2\pi}}\right] \,. \label{Function_Cal_Mb}
\end{equation}
Here we have separated the noise $W$ into two sources $W_{I_r}$ and $W_{J_r}$ (the sum of two independent Gaussian is still a Gaussian) of covariance matrices $\frac{b^t I_r}{r \Delta}$ and $\frac{(1-b^t) J_r}{r^2 \Delta}$.
The first term corresponds to $u_{k}, 1\leq k \leq n$ and the last term to $u_0$.

One observes that eq. (\ref{SE_b_Community}) has always the uniform
fixed point $b^t=0$. This is an example of a non-zero mean prior for
which nevertheless there is a uniform fixed point because other kind
of symmetry is present in the model. Let us expand (\ref{SE_b_Community}) around $b^t=0$ to determine the stability of this fixed point, one gets
\begin{equation}
b^{t+1} = \frac{b^t}{\Delta r^2} + \frac{r-4}{2 \Delta^2 r^4}{b^t}^2 + O\left( {b^t}^3 \right) \label{Expansion_b} \,.
\end{equation}
The uniform fixed point hence becomes unstable for $\Delta > \Delta_c
= \frac{1}{r^2}$.  Translated back into the parameters of the
stochastic block model this gives the easy/hard phase transition at
$|p_{\rm in} - p_{\rm out}| = r \sqrt{p_{\rm out} (1- p_{\rm out})/N}$
well known in the sparse case where $p={\rm const}/N$ from \cite{decelle2011asymptotic}.

In terms of the type of phase transitions there are two cases
\begin{itemize}
\item 2nd order for $r \leq 3$. The second term in
  (\ref{SE_b_Community}) is negative, this means that there will be a
 fixed point close-by to the uniform one  for $\Delta < \Delta_c$. The
 transition is of second  order.
\item 1st order for $r \geq 5$. The second term in
  (\ref{SE_b_Community}) is positive, this means that there will be a
  jump in the order parameter  at the transition $\Delta_c=\Delta_{\rm
  Alg}$. This is the signature a first order phase transition.
\end{itemize}
Rank $r=4$ is a marginal case in which we observed by directly solving the state evolution equations that the transition is continuous. 
We have checked numerically that no first order phase transition
exists for the symmetric community detection problem for $r \in \{2 ,3 ,4 \}$ meaning that first order phenomena appear only for $r \geq 5$.

\begin{figure}[!ht]
   \centering
\begin{minipage}{.5\textwidth}
  \centering
\includegraphics[scale=0.8]{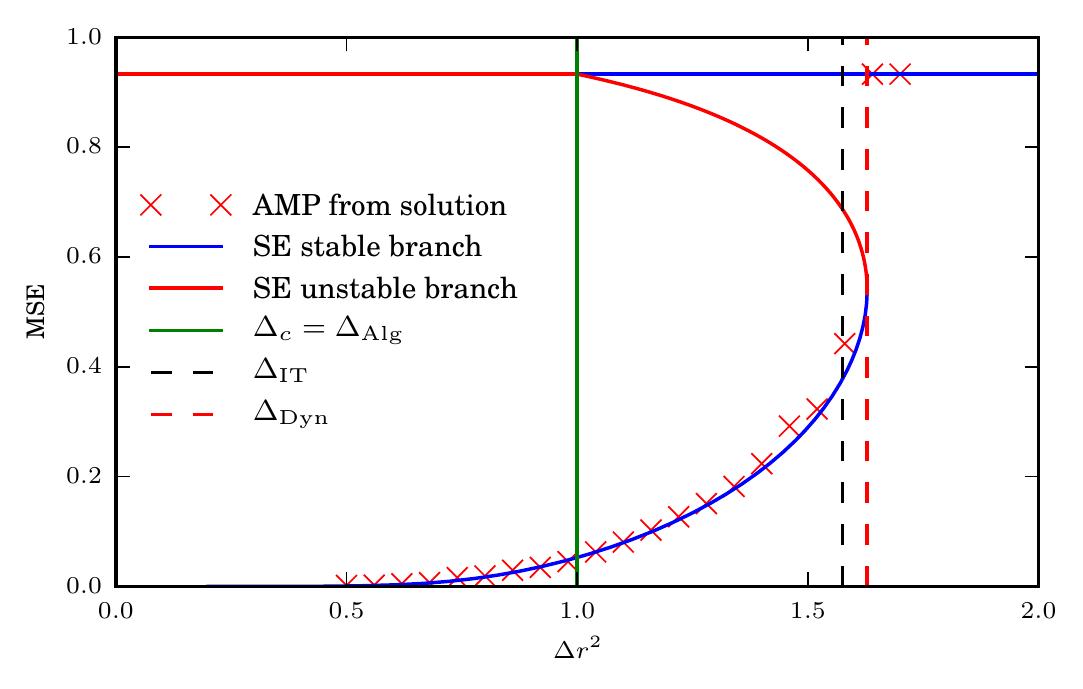}
\end{minipage}%
\begin{minipage}{.5\textwidth}
  \centering
\includegraphics[scale=0.8]{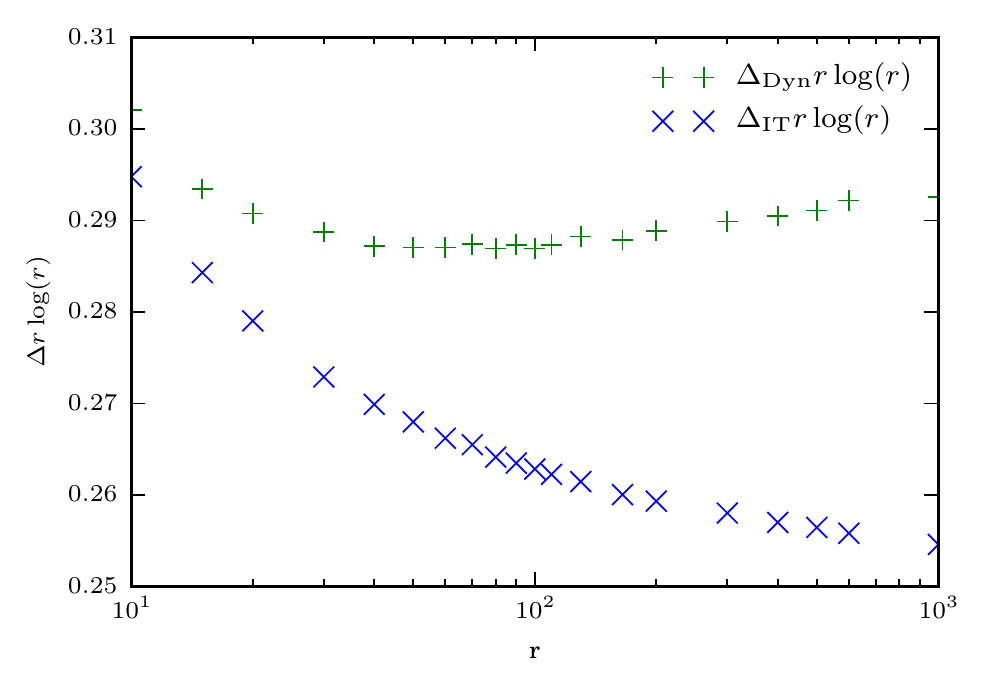}
\end{minipage}
\caption{ Left: We plot MSE deduced from state evolution (lines) and
  from Low-RAMP algorithm (marks) for 
$r=15$ groups and $N = 20 000$ as a function of $\Delta r^2$. The vertical full green
line is $\Delta_c r^2=1$. The
vertical dashed black line is $\Delta_{\rm IT}$ and the full lines
correspond to the MSE obtained from the informative initialization
and have discontinuities at $\Delta_{\rm Dyn}$.
Note that the MSE from does not go to zero at finite positive $\Delta$
instead at small noise one has ${\rm MSE} \sim \exp\left(
    -{\rm const.}/ \Delta \right)$.
Right: We plot  $\Delta r \log{r}$ for the information theoretic
$\Delta_{\rm IT}$ and dynamical spinodal $\Delta_{\rm Dyn}$ phase
transitions obtained from the state evolution using the protocol
described in Appendix~\ref{appendix:AppendixSpinodaleNetwork}. We rescale the $\Delta$ in
this way to compare with the large rank expansion in (\ref{ExpansionSpinodal}) and (\ref{ExpansionStatic}).
	}
\label{Fig_FirstOrderTransition_groups}
\end{figure}

In Fig.~\ref{Fig_FirstOrderTransition_groups} left pannel, we illustrate the first order
phase transition in the state evolution and in the behaviour of the
Low-RAMP algorithms for $r=15$ groups. 

To compute the values of $\Delta_{\rm IT}$ and $\Delta_{\rm Dyn}$, we
write $b^{t+1} = {\cal M}_r(b^t/\Delta)$ and carry a similar analysis
as in section \ref{TechniqueFirstOrderTransitionSection}, as detailed
in appendix \ref{appendix:AppendixSpinodaleNetwork}.
Fig.~\ref{Fig_FirstOrderTransition_groups} (right pannel) summarizes the values in a scaling that anticipated the
large rank expansion done is section \ref{SBM_large_rank}.

\subsection{Large sparsity (small $\rho$) expansions}
\label{sec:small_rho}

In existing literature the sparse PCA was mostly studied for sparsity
levels that are much smaller than a finite fraction of the system size
(as considered in this paper). In order to compare with existing
results we hence devote this section to the study of small density
$\rho$ expansions of the results obtained from the state evolution for
some of the models studied. 

\subsubsection{Spiked Bernoulli, Rademacher-Bernoulli, and
  Gauss-Bernoulli models}

For both the Rademacher-Bernoulli, and
  Gauss-Bernoulli models we have the uniform fixed point stable above
  $\Delta_c = \rho^2$, and in the leading order in $1/\rho$ we have
  $\Delta_{\rm alg} \sim_{\rho \to 0} \Delta_c = \rho^2$. 


The small $\rho$ limit behaviour of the information theoretic
$\Delta_{\rm IT}$ threshold and the dynamical spinodal threshold
$\Delta_{\rm Dyn}$ are given by
\begin{eqnarray}
&&\textbf{Bernoulli\, \,  and \, \, Rademacher-Bernoulli}\nonumber
\\
\Delta_{\rm Dyn}(\rho)   &\sim&_{\rho \rightarrow 0}
                                \frac{-\rho}{2\log(\rho)} \, ,\label{delta_Dyn_RB}
\\
\Delta_{\rm IT}(\rho)   &\sim&_{\rho \rightarrow 0}
                               \frac{-\rho}{4\log(\rho)}\, , \label{delta_IT_RB}
\\
&&\textbf{Gaussian-Bernoulli}\nonumber
\\
\Delta_{\rm Dyn}(\rho)   &\sim&_{\rho \rightarrow 0}
                                \frac{-\rho}{\log(\rho)} \max \left\{
                                \frac{ \frac{ 2\exp\left(
                                \frac{-1}{\beta}\right)}{\sqrt{\pi
                                \beta}} + \erfc \left(
                                \frac{1}{\sqrt{\beta}} \right)}{\beta}
                                , \beta \in \mathbb{R}^+ \right\} \sim
                                0.595 \frac{-\rho}{\log(\rho)}\, ,
\end{eqnarray}
\begin{multline}
\Delta_{\rm IT}(\rho)   \sim_{\rho \rightarrow 0}
                               \frac{-\rho}{\log(\rho)} \max \left\{
                               \frac{\frac{ 2\exp\left(
                                \frac{-1}{\beta}\right)}{\sqrt{\pi
                                \beta}} + \erfc \left(
                               \frac{1}{\sqrt{\beta}} \right)}{\beta}
                               ,
\right.\\ \left.
\int_0^\beta {\rm d}u  \, \frac{2\exp\left(
                                \frac{-1}{u}\right)}{\sqrt{\pi u}} +
                               \left( \frac{1}{\sqrt{u}} \right) =
                               \frac{1}{2}\beta\left[\frac{2\exp\left(
                                \frac{-1}{\beta}\right)}{\sqrt{\pi \beta}} + \erfc \left(
                               \frac{1}{\sqrt{\beta}} \right)  \right]; \beta
                               \in \mathbb{R}^+ \right\} \sim 0.528
                               \frac{-\rho}{\log(\rho)}\, .
\end{multline}
The information theoretic transitions for all these 3 models scale
like $O\left(\frac{-\rho}{\log(\rho)} \right)$ while the algorithmic
transition scales like $O(\rho^2)$. This means that for small $\rho$
there is a large gap between what is information theoretically and
algorithmically achievable. 
Note at this point that the bounds derived in \cite{banks2016information} for sparse PCA have the same
leading order behaviour when $\rho$ is small as (\ref{delta_Dyn_RB})
and (\ref{delta_IT_RB}).

To derive the above small $\rho$ expressions we combine
\eqref{Delta_Static_x} and \eqref{FreeEnergyEquationStatic} and the
following small $\rho$ limit of the state-evolution functions $f^{\rm
  SE}(x)$ 
\begin{eqnarray}
\forall \beta \in \mathbb{R}^+,\lim_{\rho \rightarrow 0} \frac{f^{\rm SE}(-\beta \log(\rho))}{\rho} = 1(\beta > 2),\, \textbf{Bernoulli}
\\
\forall \beta \in \mathbb{R}^+,\lim_{\rho \rightarrow 0} \frac{f^{\rm SE}(-\beta \log(\rho))}{\rho} = 1(\beta > 2),\, \textbf{Rademacher-Bernoulli}
\label{SmallRhoLimit_Rademacher_Bernoulli}
\\
\forall \beta \in \mathbb{R}^+,\lim_{\rho \rightarrow 0} \frac{f^{\rm SE}(-\beta \log(\rho))}{\rho} =\frac{2 \exp\left( \frac{-1}{\beta}\right)}{\sqrt{\beta \pi}} + \erfc\left(\frac{1}{\sqrt{\beta}} \right),\, \textbf{Gaussian-Bernoulli}
\label{SmallRhoLimit_Gauss_Bernoulli}
\end{eqnarray}
Here the functions $f^{\rm SE}$ are the state-evolution update
functions stated in \eqref{SE_Bernoulli},
\eqref{SE_Rademacher_Bernoulli} and \eqref{SE_Gaussian_Bernoulli} for
the different models (when the model is clear from context we will omit the lower index
specifying the model). 
The above is proven by deriving the state evolution equations for each of
these models. The computation is done in appendix \ref{AppendixSmallRho}.

\subsubsection{Two balanced groups, limit of small planted subgraph} 


In this section we analyze the small $\rho$ limit for the two balanced
groups of section \ref{2-groups-balanced}. From the definition of the function
$f_{\rho}$ (\ref{frho_2groups}), and a computation done in appendix \ref{AppendixSmallRho} we get 
\begin{equation}
\lim_{\rho \rightarrow 0} f_{\rho}(
-\beta\rho(1-\rho)\log(\rho(1-\rho))) = {1}(\beta > 2)\, .
\end{equation}
By combining this with \eqref{Delta_Static_x} and \eqref{FreeEnergyEquationStatic} one gets
\begin{eqnarray}
\Delta_{\rm Dyn}(\rho) &\sim&_{\rho \rightarrow 0}
                              \frac{1}{-2\rho(1-\rho)\log(\rho(1-\rho))}\, ,
\\
	\Delta_{\rm IT}(\rho) &\sim&_{\rho \rightarrow 0}
                                     \frac{1}{-4\rho(1-\rho)\log(\rho(1-\rho))}\, .
\end{eqnarray}
The derivations of these limits is done in the appendix
\ref{AppendixSmallRho}.

Note that the limit of small $\rho$ in the two balanced groups
model is closely related to the problem of planted clique. However, in
the planted clique problem the size of the clique to recover is much
smaller that the size of the graph, typically $O(\sqrt{N})$ for efficient
recovery and $O(\log{N})$ for information theoretically possible
recovery. Whereas our equation were derives when $\rho = O(1)$, we can
try to see what would the above scaling imply for $k = \rho N = o(N)$. 
In the planted clique problem the first (smaller) group is fully
connected, this means that the entry $C_{11}$ of the connectivity
matrix $C$ (\ref{ConnectivityMatrix}) is equal to one. Therefore for
$\rho \ll 1$ one has 
\begin{equation}
\mu = (1-p_{\rm out}) \rho \sqrt{N}\, .
\end{equation}
Note that in the canonical definition of the planted clique problem
the average degree of the nodes belonging to the clique are slightly
larger than the average degree of the rest of the graph. The present
case of balanced groups corresponds to a version of the planted clique
problem where next to planting a clique a corresponding number of
edges is added to the rest of the graph to ensure that the average
degree of every node is the same. 

Recalling the definition of $\Delta$ (\ref{delta_2groups}) for the
community detection output channel, and using $\Delta_c=1$ for the
spectral threshold, and $\Delta_{\rm IT} = -1/(4\rho \log \rho)$ for
the information theoretic threshold at small $\rho$ we get 
\begin{eqnarray}
k_c = \sqrt{N}\sqrt{\frac{p_{\rm out}}{1-p_{\rm out}}}\, ,
\\
k_{\rm IT} = \log(N)\frac{4p_{\rm out}}{1-p_{\rm out}}\, .
\end{eqnarray}
We indeed recover the scaling known from the planted clique problem,
see e.g. \cite{deshpande2015finding}.
The $p_{\rm out}$-dependent constant are indeed the tight constants
(for efficient and information theoretically optimal recovery) for the
balanced version of the planted clique where the expected degree of
every node is the same independently of the fact in the node is in the
clique or not.


\subsubsection{Sparse PCA at small density $\rho$}  

We investigate here the small $\rho$ limit
of the bipartite $UV^\top$ spiked Gaussian-Bernoulli, and Rademacher-Bernoulli
model. We remind that $U \in {\mathbb R}^N$, $V \in {\mathbb R}^M$,
while $\alpha=M/N$.  In the model we consider that elements of $U$ are Gaussian of
zero mean and unit variance, 
while $P_V$ is given by (\ref{P_X_Rademacher}) for the Rademacher-Bernoulli model, and by	(\ref{P_X_GaussBernoulli}) for the Gauss-Bernoulli model.
The state evolution equations then read
\begin{eqnarray}
m_u^{t} &=& \frac{\alpha m_v^t}{\Delta + \alpha m_v^t}\, ,
\\
m_v^{t+1} &=& f_{\rm Gauss-Bernoulli}^{\rm SE}\left(
              \frac{m_u^t}{\Delta}\right)\, ,
\\
m_v^{t+1} &=& f_{\rm Rademacher-Bernoulli}^{\rm SE}\left(
              \frac{m_u^t}{\Delta}\right)\, ,
\end{eqnarray}
where $f^{\rm SE}_{\rm Gauss-Bernoulli}$ and $f^{\rm SE}_{\rm Rademacher-Bernoulli}$ are defined in (\ref{SE_Rademacher_Bernoulli}) and (\ref{SE_Gaussian_Bernoulli}) respectively.
By combining these equations one gets
\begin{eqnarray}
m_v^{t+1} &=& f_{\rm Gauss-Bernoulli}^{\rm SE}\left( \frac{\alpha
              m_v^t}{\Delta^2 + \alpha \Delta m_v^t}\right)\, , \label{SE_UV_GAUSS_BERNOULLI}
\\
m_v^{t+1} &=& f_{\rm Rademacher-Bernoulli}^{\rm SE}\left( \frac{\alpha
              m_v^t}{\Delta^2 + \alpha \Delta m_v^t}\right) \, .\label{SE_UV_RADEMACHER_BERNOULLI}
\end{eqnarray}

Because in both these cases $P_U$ and $P_V$ have zero mean the state evolution equations will have the uniform fixed point at $(m_u,m_v) = (0,0)$. This  fixed point becomes unstable when
\begin{equation}
\frac{\rho^2 \alpha}{\Delta^2} > 1, \text{or}, \Delta < \Delta_c =
\rho  \sqrt{\alpha} \, .  \label{Deltac_GB_bip}
\end{equation}
Also in this case the stability transition $\Delta_c$ corresponds to
the spectral transition where one sees informative eigenvalues get out
of the bulk of  the matrix $S$, as is known in the theory of low-rank
perturbations of random matrices \cite{baik2005phase} (these methods
are known not to take advantage of the sparsity). 
For $\rho$ small enough there will be
again a first order phase transitions with $\Delta_{\rm Alg}$,
$\Delta_{\rm IT}$, $\Delta_{\rm Dyn}$ defined as in section \ref{SubSectionFirstOrderGeneral}. 
The asymptotic behaviour of these thresholds is
\begin{eqnarray}
&&\textbf{Rademacher-Bernoulli}\nonumber
\\
\Delta_{\rm Dyn}(\rho)   &\sim&_{\rho \rightarrow 0}
                                \sqrt{\frac{-\rho \alpha}{2 \log(\rho)}
                                }\, , 
\\
\Delta_{\rm IT}(\rho)   &\sim&_{\rho \rightarrow 0}
                               \sqrt{\frac{- \rho \alpha}{4\alpha
                               \log(\rho)}}\, ,   \label{sparse_PCA_IT}
\\
&&\textbf{Gaussian-Bernoulli}\nonumber
\\
\Delta_{\rm Dyn}(\rho)   &\sim&_{\rho \rightarrow 0}
                                \sqrt{\frac{-\rho \alpha}{\log(\rho)}
                                \max \left\{  \frac{\erfc \left(
                                \frac{1}{\sqrt{\beta}} \right)}{\beta}
                                , \beta \in \mathbb{R}^+ \right\} }  \sim 0.771\sqrt{\frac{-\rho \alpha}{\log(\rho)} }  \, ,
\end{eqnarray}
\begin{multline}
\quad \quad \quad\quad\quad\quad\quad  \quad \Delta_{\rm IT}(\rho)   \sim_{\rho \rightarrow 0} \sqrt{
                               \frac{-\rho \alpha}{ \log(\rho)} }  
                                \\ \sqrt{\max
                               \left\{  \frac{ \frac{2\exp\left(
                                \frac{-1}{\beta}\right)}{\sqrt{\pi \beta}}+ \erfc \left(
                               \frac{1}{\sqrt{\beta}} \right)}{\beta}
                               , \int_0^\beta {\rm d}u  \, \frac{2\exp\left(
                                \frac{-1}{u}\right)}{\sqrt{\pi u}} +
                               \left( \frac{1}{\sqrt{u}} \right) =
                               \frac{1}{2}\beta\left[\frac{2\exp\left(
                                \frac{-1}{\beta}\right)}{\sqrt{\pi \beta}} + \erfc \left(
                               \frac{1}{\sqrt{\beta}} \right)  \right]; \beta
                               \in \mathbb{R}^+ \right\} } \\
\sim 0.726 \sqrt{\frac{-\rho \alpha}{\log(\rho)} }\, .
\end{multline}
Reminding that the value of $\Delta_{\rm Alg}$ scales in the same way
as $\Delta_c$ (\ref{Deltac_GB_bip}), we see that a large hard phase
opens as $\rho \to 0$.

To put the above results in relation to existing literature on sparse
PCA \cite{amini2008high,deshpande2014sparse}. The main difference is that the regime considered in
existing literature is that the number of non-zero element in the
matrix $V$, $\rho N = o(1)$. The information theoretic threshold found
in eq.~(\ref{sparse_PCA_IT}) correspond (up to a constant) to
information theoretic bounds found in \cite{amini2008high}. However,
the algorithmic performance that is found in the case of very small
sparsity, by e.g. the covariance thresholding of \cite{deshpande2014sparse}, is not reproduced in our analysis
of linear sparsity $\rho=O(1)$. This suggest that in case when the
sparsity is small but linear in $N$, efficient algorithms that take advantage of
the sparsity might not exist. This regime should be investigated further. 

To derive the small $\rho$ limit we follow a similar strategy as we
did in the symmetric $XX^\top$ case.
The main idea is to find the value of  $\Delta$ for which $m_v = f^{\rm SE}_{P_X}(x)$ is a fixed point of the SE equations (\ref{SE_UV_RADEMACHER_BERNOULLI}) or (\ref{SE_UV_GAUSS_BERNOULLI}).
\begin{equation}
\Delta(x) = \frac{\alpha f^{\rm SE}(x) + \sqrt{\alpha^2 {f^{\rm SE}}^2(x) + 4 \alpha
    \frac{f^{\rm SE}(x)}{x}}}{2}\, ,
\end{equation}
This means that
\begin{equation}
(m_u,m_v) = (x \Delta(x),f^{\rm SE}(x))
\end{equation}
is a fixed point of the state evolution equations for $\Delta =\Delta(x)$.
The free-energy is computed as
\begin{eqnarray}
\phi(m_u = x \Delta(x),m_v = f(x),\Delta = \Delta(x)) - \phi(m_u =
  0,m_v = 0,\Delta = \Delta(x)) = \nonumber \\ = \alpha\int_0^{x} {\rm
  du}f^{\rm SE}(u) + \int_0^{\alpha f^{\rm SE}(x)/\Delta(x)}{\rm d}u
  \frac{u}{1 + u} - xf^{\rm SE}(x)\, .
\end{eqnarray}
Combining this with (\ref{SmallRhoLimit_Gauss_Bernoulli}) and
(\ref{SmallRhoLimit_Rademacher_Bernoulli}) one gets the above asymptotic behaviour.

\subsection{Large rank expansions}
\label{sec:large_r}

Another limit that can be worked out analytically is the large rank
$r$ limit. This section summarizes the results.

\subsubsection{Large-rank limit for jointly-sparse PCA}

\begin{figure}
\includegraphics[scale=1]{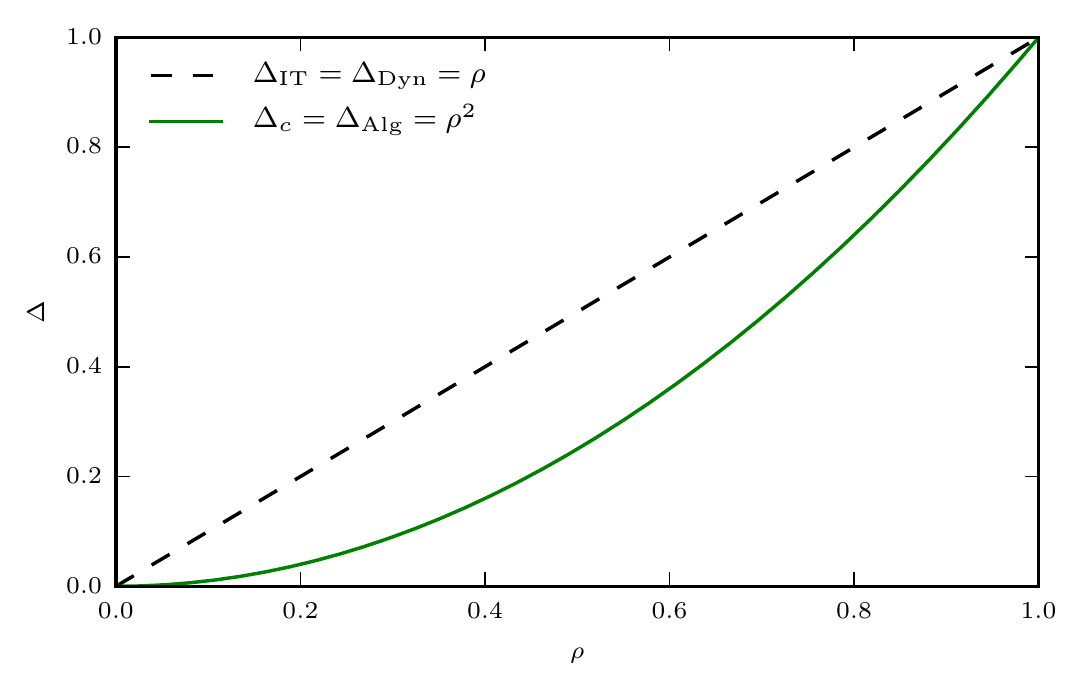}
\caption{Phase diagram of the jointly-sparse PCA model at large rank
  $r$. In the large $r$ limit the algorithmic spinodal merges with
  $\Delta_c$, the dynamical spinodal and the information theoretic one converge toward the line $\Delta = \rho$.}
\label{Infinite_r_phase_diagram}
\end{figure}

We analyze the large $r$ limit of jointly-sparse PCA for which the
state evolution equation is given by (\ref{JointPCA_DE}). We notice
that $u^2$ will have mean $r$ and standard deviation $\sqrt{2r}$. This
essentially means that to get the large $r$ limit of the density
evolution equations one need to replace $u^2$ by $r$ everywhere it
appears. Expanding in large $r$ then gives 
\begin{equation}
m^{t+1} = \frac{\rho m^t}{m^t+\Delta} \left( 1 +\frac{m^t}{\Delta}(1-\hat\rho) \right)\hat\rho \,,
\end{equation}
where 
\begin{equation}
\hat\rho = \lim\limits_{r \rightarrow +\infty} \frac{\rho}{(1-\rho) \exp\left( \frac{r}{2}\left[\frac{m^t}{\Delta} + \log(1+\frac{m^t}{\Delta} \right) \right] + \rho} = \left\{
\begin{matrix}
1 \, \,  {\rm if} \, \,  r {m^t}^2 \rightarrow 0
\\
\rho \, \, {\rm if} \, \, r {m^t}^2 \rightarrow + \infty
\end{matrix}
\right. \,.
\end{equation}
This means that for any $m^t \gg \frac{1}{\sqrt{r}}$ the update equations will be approximately
\begin{equation}
m^{t+1} = \frac{m^t \rho}{m^t + \Delta} + o(1) \,.
\end{equation}
This update equation can be easily analyzed, it only has one stable fixed point located at
\begin{equation}
\max(\rho - \Delta,0) \,.
\end{equation}
Analogous expansion of the replicated free energy leads to the result
that in the large rank limit we have $\Delta_{\rm Dyn} = \Delta_{\rm
  IT} = \rho$ whereas $\Delta_{\rm Alg} = \Delta_c = \rho^2$.  
This is plotted in the large rank phase diagram (\ref{Infinite_r_phase_diagram}).


\subsubsection{Community detection}
\label{SBM_large_rank}

The large rank limit is analyzed for the problem of symmetric
community detection in appendix \ref{appendix:AppendixSpinodaleNetwork}.
The asymptotic behavior of $\Delta_{\rm IT}(r)$ and $\Delta_{\rm Dyn}(r)$ as $r \rightarrow +\infty$ are
\begin{eqnarray}
\Delta_c &=& \frac{1}{r^2} \, ,
\\
\Delta_{\rm Dyn} &=& \frac{1}{2 r \ln(r)} [1+ o_r(1)] \, , \label{ExpansionSpinodal}\\
\Delta_{\rm IT} &=& \frac{1}{4 r \ln(r)} [1+ o_r(1)]\, .\label{ExpansionStatic}
\end{eqnarray}
We see that a large gap opens between $\Delta_c$ and $\Delta_{\rm IT}$ as $r$ grows.
The behavior $\Delta_{\rm IT}$ and $\Delta_{\rm Dyn}$ for moderately
large $r$ is illustrated in figure \ref{Fig_FirstOrderTransition_groups} and we see
that the above limit is reached very slowly.
Using eq. (\ref{Delta_SBM}) this translates into the large $r$ limit
phase transition in terms of parameters of the stochastic block model
as discussed in \cite{lesieur2015mmse}, and proven in \cite{banks2016information}.

\section{Discussion and perspective}
The main results of this paper are summarised in section
\ref{Main_results}.  This concluding section is devoted to open
directions and questions.

In this paper we considered models that lead to fully connected
graphical models. It should be noted that the results and phase
diagrams derived here in this setting are expected to hold as long as
the average degree of the underlying factor graph is large enough
(goes to infinity, after the $N\to \infty$ limit has been taken). 
Such expectation was proved rigorously in some special cases, see
e.g. \cite{deshpande2016asymptotic,caltagirone2016recovering}.

A large part of results reported in this paper concerns the
Bayes-optimal setting for which the mutual information and MMSE was
recently established rigorously
\cite{barbier2016mutual,LelargeMiolane16,miolane2017fundamental}.  

A formidable open question concerns the nature of the hard phase, that
is the region of parameters for which Low-RAMP does not reach the
MMSE. Such regions have been identified in a considerable number of
inference problems and it remains a great challenge to study whether
efficient algorithms that perform better than Low-RAMP exist or to
provide further evidence towards their non-existence. One notable
direction is that whereas the MMSE in Bayes-optimal setting is
correctly described by the replica symmetric solution, the nature of
the metastable branch describing the behaviour of Low-RAMP in the hard
phase may involve effect of replica symmetry breaking, this should be
investigated in future work. 

The results presented in this paper focused on the Bayes-optimal
inference, the theory and methodology presented is not limited to this
case and phase diagrams involving mismatch between the teacher and
student models or parameters could provide further insight into issues such as
model and parameter selection and learning, overfitting, or performance of
minimisation-based algorithms. Effects of replica symmetry breaking
are expected to play a role in this case. This is another interesting direction
for future work.

In the present paper we focused on pair-wise interaction corresponding
to matrix factorizations. In physics as well as in statistics a
generalization to $p$-spin interaction have been studied. In
statistics this type of problems is termed tensor PCA
\cite{richard2014statistical}. Detailed study of tensor PCA is also an interesting
extension. 

In the present paper we treated real-valued vectors spins. The
generalization to spin-values that live on more involved spaces such
as compact groups was recently considered in \cite{perry2016message}
and many interesting questions in this direction remain open. 

Finally in the present work we concentrated on the theoretical
development and phase diagrams for randomly generated data. The
Low-RAMP algorithm (that is distributed in Matlab and Julia at \url{http://krzakala.github.io/LowRAMP/}
) can be successfully applied to real-world data. In that case the set
of theoretical guarantees on its performance is limited, but this is
the case common to many other practically used
approaches. Investigating performance of the Low-RAMP algorithm on
data coming from applications of interest is another formidably rich
and important direction for future work.

\section{Acknowledgement}
FK acknowledges funding from the EU (FP/2007-2013/ERC grant agreement
307087-SPARCS). LZ acknowledges funding for the project SaMURai from
LabEx PALM (Ref : ANR-10-LABX-0039-PALM). We would like to thank Marc
M{\'e}zard, Cristopher Moore, Pierfrancesco Urbani, and Christian
Schmidt for discussions related to this work, and also Stefano Sarao, Fabrizio Antenucci, and Alexis Bozio for proofreading of the manuscript. 

\bibliographystyle{unsrt}
\bibliography{Reference}

\begin{thebibliography}{10}

\bibitem{lesieur2015phase}
Thibault Lesieur, Florent Krzakala, and Lenka Zdeborov{\'a}.
\newblock Phase transitions in sparse {PCA}.
\newblock In {\em IEEE International Symposium on Information Theory
  Proceedings (ISIT)}, pages 1635--1639, 2015.

\bibitem{lesieur2015mmse}
Thibault Lesieur, Florent Krzakala, and Lenka Zdeborov\'a.
\newblock {MMSE} of probabilistic low-rank matrix estimation: Universality with
  respect to the output channel.
\newblock In {\em 53rd Annual Allerton Conference on Communication, Control,
  and Computing (Allerton)}, pages 680--687. IEEE, 2015.

\bibitem{gabrie2015training}
Marylou Gabri{\'e}, Eric~W Tramel, and Florent Krzakala.
\newblock Training restricted {B}oltzmann machine via the
  {T}houless-{A}nderson-{P}almer free energy.
\newblock In {\em Advances in Neural Information Processing Systems}, pages
  640--648, 2015.

\bibitem{tramel2016inferring}
Eric~W Tramel, Andre Manoel, Francesco Caltagirone, Marylou Gabri{\'e}, and
  Florent Krzakala.
\newblock Inferring sparsity: Compressed sensing using generalized restricted
  {B}oltzmann machines.
\newblock In {\em IEEE Information Theory Workshop (ITW)}, pages 265--269.
  IEEE, 2016.

\bibitem{NishimoriBook}
H.~Nishimori.
\newblock {\em Statistical Physics of Spin Glasses and Information Processing:
  An Introduction}.
\newblock Oxford University Press, Oxford, UK, 2001.

\bibitem{zdeborova2015statistical}
Lenka Zdeborov{\'a} and Florent Krzakala.
\newblock Statistical physics of inference: Thresholds and algorithms.
\newblock {\em Advances in Physics}, 65:453--552, 2016.

\bibitem{nishimori2001absence}
Hidetoshi Nishimori and David Sherrington.
\newblock Absence of replica symmetry breaking in a region of the phase diagram
  of the ising spin glass.
\newblock In {\em American Institute of Physics Conference Series}, volume 553,
  pages 67--72, 2001.

\bibitem{eckart1936approximation}
Carl Eckart and Gale Young.
\newblock The approximation of one matrix by another of lower rank.
\newblock {\em Psychometrika}, 1(3):211--218, 1936.

\bibitem{baik2005phase}
Jinho Baik, G{\'e}rard Ben~Arous, and Sandrine P{\'e}ch{\'e}.
\newblock Phase transition of the largest eigenvalue for nonnull complex sample
  covariance matrices.
\newblock {\em Annals of Probability}, pages 1643--1697, 2005.

\bibitem{SherringtonKirkpatrick75}
D.~Sherrington and S.~Kirkpatrick.
\newblock Solvable model of a spin-glass.
\newblock {\em Phys. Rev. Lett.}, 35:1792--1796, 1975.

\bibitem{ThoulessAnderson77}
D.~J. Thouless, P.~W. Anderson, and R.~G. Palmer.
\newblock Solution of 'solvable model of a spin glass'.
\newblock {\em Philosophical Magazine}, 35(3):593--601, 1977.

\bibitem{kosterlitz1976spherical}
JM~Kosterlitz, DJ~Thouless, and Raymund~C Jones.
\newblock Spherical model of a spin-glass.
\newblock {\em Physical Review Letters}, 36(20):1217, 1976.

\bibitem{sommers1981theory}
Hans-Juergen Sommers.
\newblock Theory of a {H}eisenberg spin glass.
\newblock {\em Journal of Magnetism and Magnetic Materials}, 22(3):267--270,
  1981.

\bibitem{gabay1981coexistence}
Marc Gabay and G{\'e}rard Toulouse.
\newblock Coexistence of spin-glass and ferromagnetic orderings.
\newblock {\em Physical Review Letters}, 47(3):201, 1981.

\bibitem{hinton2010practical}
Geoffrey Hinton.
\newblock A practical guide to training restricted {B}oltzmann machines.
\newblock {\em Momentum}, 9(1):926, 2010.

\bibitem{hinton2006fast}
Geoffrey~E Hinton, Simon Osindero, and Yee-Whye Teh.
\newblock A fast learning algorithm for deep belief nets.
\newblock {\em Neural computation}, 18(7):1527--1554, 2006.

\bibitem{lecun2015deep}
Yann LeCun, Yoshua Bengio, and Geoffrey Hinton.
\newblock Deep learning.
\newblock {\em Nature}, 521(7553):436--444, 2015.

\bibitem{kappen1998boltzmann}
Hilbert~J Kappen and FB~Rodriguez.
\newblock Boltzmann machine learning using mean field theory and linear
  response correction.
\newblock {\em Advances in neural information processing systems}, pages
  280--286, 1998.

\bibitem{hopfield1982neural}
John~J Hopfield.
\newblock Neural networks and physical systems with emergent collective
  computational abilities.
\newblock {\em Proceedings of the national academy of sciences},
  79(8):2554--2558, 1982.

\bibitem{mezard2016mean}
Marc M\'ezard.
\newblock Mean-field message-passing equations in the hopfield model and its
  generalizations.
\newblock {\em Phys. Rev. E}, 95:022117, Feb 2017.

\bibitem{hastie2005elements}
Trevor Hastie, Robert Tibshirani, Jerome Friedman, and James Franklin.
\newblock The elements of statistical learning: data mining, inference and
  prediction.
\newblock {\em The Mathematical Intelligencer}, 27(2):83--85, 2005.

\bibitem{wasserman2013all}
Larry Wasserman.
\newblock {\em All of statistics: a concise course in statistical inference}.
\newblock Springer Science \& Business Media, 2013.

\bibitem{lloyd1982least}
Stuart Lloyd.
\newblock Least squares quantization in {PCM}.
\newblock {\em IEEE Transactions on Information Theory}, 28(2):129--137, 1982.

\bibitem{watkin1994optimal}
TLH Watkin and J-P Nadal.
\newblock Optimal unsupervised learning.
\newblock {\em Journal of Physics A: Mathematical and General}, 27(6):1899,
  1994.

\bibitem{barkai1994statistical}
N~Barkai and Haim Sompolinsky.
\newblock Statistical mechanics of the maximum-likelihood density estimation.
\newblock {\em Physical Review E}, 50(3):1766, 1994.

\bibitem{biehl1994statistical}
Michael Biehl and Andreas Mietzner.
\newblock Statistical mechanics of unsupervised structure recognition.
\newblock {\em Journal of Physics A: Mathematical and General}, 27(6):1885,
  1994.

\bibitem{NIPS2013_5074}
Ryosuke Matsushita and Toshiyuki Tanaka.
\newblock Low-rank matrix reconstruction and clustering via approximate message
  passing.
\newblock In C.J.C. Burges, L.~Bottou, M.~Welling, Z.~Ghahramani, and K.Q.
  Weinberger, editors, {\em Advances in Neural Information Processing Systems
  26}, pages 917--925. Curran Associates, Inc., 2013.

\bibitem{parsons2004subspace}
Lance Parsons, Ehtesham Haque, and Huan Liu.
\newblock Subspace clustering for high dimensional data: a review.
\newblock {\em ACM SIGKDD Explorations Newsletter}, 6(1):90--105, 2004.

\bibitem{lesieur2016phase}
Thibault Lesieur, Caterina De~Bacco, Jess Banks, Florent Krzakala, Cris Moore,
  and Lenka Zdeborov{\'a}.
\newblock Phase transitions and optimal algorithms in high-dimensional
  {G}aussian mixture clustering.
\newblock In {\em 54th Annual Allerton Conference on Communication, Control,
  and Computing (Allerton)}, pages 601--608. IEEE, 2016.

\bibitem{johnstone2004sparse}
Iain~M Johnstone and Arthur~Yu Lu.
\newblock Sparse principal components analysis.
\newblock {\em Unpublished manuscript}, 7, 2004.

\bibitem{zou2006sparse}
Hui Zou, Trevor Hastie, and Robert Tibshirani.
\newblock Sparse principal component analysis.
\newblock {\em Journal of computational and graphical statistics},
  15(2):265--286, 2006.

\bibitem{amini2008high}
Arash~A Amini and Martin~J Wainwright.
\newblock High-dimensional analysis of semidefinite relaxations for sparse
  principal components.
\newblock In {\em IEEE International Symposium on Information Theory}, pages
  2454--2458. IEEE, 2008.

\bibitem{deshpande2014sparse}
Yash Deshpande and Andrea Montanari.
\newblock Sparse {PCA} via covariance thresholding.
\newblock In {\em Advances in Neural Information Processing Systems}, pages
  334--342, 2014.

\bibitem{krauthgamer2015semidefinite}
Robert Krauthgamer, Boaz Nadler, Dan Vilenchik, et~al.
\newblock Do semidefinite relaxations solve sparse {PCA} up to the information
  limit?
\newblock {\em The Annals of Statistics}, 43(3):1300--1322, 2015.

\bibitem{berthet2013computational}
Quentin Berthet and Philippe Rigollet.
\newblock Complexity theoretic lower bounds for sparse principal component
  detection.
\newblock In {\em COLT}, pages 1046--1066, 2013.

\bibitem{DeshpandeM14}
Yash Deshpande and Andrea Montanari.
\newblock Information-theoretically optimal sparse {PCA}.
\newblock In {\em IEEE International Symposium on Information Theory (ISIT)},
  pages 2197--2201. IEEE, 2014.

\bibitem{monasson2015estimating}
R{\'e}mi Monasson and Dario Villamaina.
\newblock Estimating the principal components of correlation matrices from all
  their empirical eigenvectors.
\newblock {\em EPL (Europhysics Letters)}, 112(5):50001, 2015.

\bibitem{madeira2004biclustering}
Sara~C Madeira and Arlindo~L Oliveira.
\newblock Biclustering algorithms for biological data analysis: a survey.
\newblock {\em IEEE/ACM Transactions on Computational Biology and
  Bioinformatics (TCBB)}, 1(1):24--45, 2004.

\bibitem{cheng2000biclustering}
Yizong Cheng and George~M Church.
\newblock Biclustering of expression data.
\newblock In {\em Ismb}, volume~8, pages 93--103, 2000.

\bibitem{fortunato2010community}
Santo Fortunato.
\newblock Community detection in graphs.
\newblock {\em Physics Reports}, 486(3):75--174, 2010.

\bibitem{decelle2011inference}
Aurelien Decelle, Florent Krzakala, Cristopher Moore, and Lenka Zdeborov{\'a}.
\newblock Inference and phase transitions in the detection of modules in sparse
  networks.
\newblock {\em Physical Review Letters}, 107(6):065701, 2011.

\bibitem{decelle2011asymptotic}
Aurelien Decelle, Florent Krzakala, Cristopher Moore, and Lenka Zdeborov{\'a}.
\newblock Asymptotic analysis of the stochastic block model for modular
  networks and its algorithmic applications.
\newblock {\em Physical Review E}, 84(6):066106, 2011.

\bibitem{GrossKanter85}
D.~J. Gross, I.~Kanter, and H.~Sompolinsky.
\newblock Mean-field theory of the potts glass.
\newblock {\em Phys. Rev. Lett.}, 55(3):304--307, Jul 1985.

\bibitem{deshpande2016asymptotic}
Yash Deshpande, Emmanuel Abbe, and Andrea Montanari.
\newblock Asymptotic mutual information for the binary stochastic block model.
\newblock In {\em IEEE International Symposium on Information Theory (ISIT)},
  pages 185--189. IEEE, 2016.

\bibitem{barbier2016mutual}
Jean Barbier, Mohamad Dia, Nicolas Macris, Florent Krzakala, Thibault Lesieur,
  and Lenka Zdeborov{'a}.
\newblock Mutual information for symmetric rank-one matrix estimation: A proof
  of the replica formula.
\newblock In {\em Advances In Neural Information Processing Systems}, pages
  424--432, 2016.

\bibitem{caltagirone2016recovering}
Francesco Caltagirone, Marc Lelarge, and L{\'e}o Miolane.
\newblock Recovering asymmetric communities in the stochastic block model.
\newblock {\em arXiv preprint arXiv:1610.03680}, 2016.

\bibitem{montanari2015finding}
Andrea Montanari.
\newblock Finding one community in a sparse graph.
\newblock {\em Journal of Statistical Physics}, 161(2):273--299, 2015.

\bibitem{tubiana2016emergence}
J.~Tubiana and R.~Monasson.
\newblock Emergence of compositional representations in restricted boltzmann
  machines.
\newblock {\em Phys. Rev. Lett.}, 118:138301, Mar 2017.

\bibitem{krzakala2013phase}
Florent Krzakala, Marc M{\'e}zard, and Lenka Zdeborov{\'a}.
\newblock Phase diagram and approximate message passing for blind calibration
  and dictionary learning.
\newblock In {\em IEEE International Symposium on Information Theory
  Proceedings (ISIT)}, pages 659--663. IEEE, 2013.

\bibitem{parker2014bilinear}
Jason~T Parker, Philip Schniter, and Volkan Cevher.
\newblock Bilinear generalized approximate message passing — part {I}:
  Derivation.
\newblock {\em IEEE Transactions on Signal Processing}, 62(22):5839--5853,
  2014.

\bibitem{kabashima2014phase}
Yoshiyuki Kabashima, Florent Krzakala, Marc M{\'e}zard, Ayaka Sakata, and Lenka
  Zdeborov{\'a}.
\newblock Phase transitions and sample complexity in bayes-optimal matrix
  factorization.
\newblock {\em IEEE Transactions on Information Theory}, 62(7):4228--4265,
  2016.

\bibitem{rangan2012iterative}
Sundeep Rangan and Alyson~K Fletcher.
\newblock Iterative estimation of constrained rank-one matrices in noise.
\newblock In {\em IEEE International Symposium on Information Theory
  Proceedings (ISIT)}, pages 1246--1250. IEEE, 2012.

\bibitem{javanmard2013state}
Adel Javanmard and Andrea Montanari.
\newblock State evolution for general approximate message passing algorithms,
  with applications to spatial coupling.
\newblock {\em Information and Inference}, page iat004, 2013.

\bibitem{krzakala2016mutual}
Florent Krzakala, Jiaming Xu, and Lenka Zdeborov{\'a}.
\newblock Mutual information in rank-one matrix estimation.
\newblock In {\em Information Theory Workshop (ITW), 2016 IEEE}, pages 71--75.
  IEEE, 2016.

\bibitem{LelargeMiolane16}
L{\'e}o~Miolane Marc~Lelarge.
\newblock Fundamental limits of symmetric low-rank matrix estimation.
\newblock {\em arXiv:1611.03888 [math.PR]}, 2016.

\bibitem{miolane2017fundamental}
L{\'e}o Miolane.
\newblock Fundamental limits of low-rank matrix estimation.
\newblock {\em arXiv preprint arXiv:1702.00473}, 2017.

\bibitem{MezardParisi87b}
M.~M{\'e}zard, G.~Parisi, and M.~A. Virasoro.
\newblock {\em Spin-Glass Theory and Beyond}, volume~9 of {\em Lecture Notes in
  Physics}.
\newblock World Scientific, Singapore, 1987.

\bibitem{perry2016message}
Amelia Perry, Alexander~S Wein, Afonso~S Bandeira, and Ankur Moitra.
\newblock Message-passing algorithms for synchronization problems over compact
  groups.
\newblock {\em arXiv preprint arXiv:1610.04583}, 2016.

\bibitem{krzakala2014variational}
Florent Krzakala, Andre Manoel, Eric~W Tramel, and Lenka Zdeborov{\'a}.
\newblock Variational free energies for compressed sensing.
\newblock In {\em IEEE International Symposium on Information Theory}, pages
  1499--1503. IEEE, 2014.

\bibitem{DBLP:journals/corr/abs-1301-6295}
Sundeep Rangan, Philip Schniter, Erwin Riegler, Alyson Fletcher, and Volkan
  Cevher.
\newblock Fixed points of generalized approximate message passing with
  arbitrary matrices.
\newblock In {\em IEEE International Symposium on Information Theory
  Proceedings (ISIT)}, pages 664--668. IEEE, 2013.

\bibitem{DBLP:journals/corr/VilaSRKZ14}
Jeremy Vila, Philip Schniter, Sundeep Rangan, Florent Krzakala, and Lenka
  Zdeborov{\'a}.
\newblock Adaptive damping and mean removal for the generalized approximate
  message passing algorithm.
\newblock In {\em IEEE International Conference on Acoustics, Speech and Signal
  Processing (ICASSP)}, pages 2021--2025. IEEE, 2015.

\bibitem{perry2016optimality}
Amelia Perry, Alexander~S Wein, Afonso~S Bandeira, and Ankur Moitra.
\newblock Optimality and sub-optimality of {PCA} for spiked random matrices and
  synchronization.
\newblock {\em arXiv preprint arXiv:1609.05573}, 2016.

\bibitem{banks2016information}
Jess Banks, Cristopher Moore, Roman Vershynin, and Jiaming Xu.
\newblock Information-theoretic bounds and phase transitions in clustering,
  sparse {PCA}, and submatrix localization.
\newblock {\em arXiv preprint arXiv:1607.05222}, 2016.

\bibitem{Yedidia:2003}
Jonathan~S. Yedidia, William~T. Freeman, and Yair Weiss.
\newblock Exploring artificial intelligence in the new millennium.
\newblock chapter Understanding Belief Propagation and Its Generalizations,
  pages 239--269. Morgan Kaufmann Publishers Inc., San Francisco, CA, USA,
  2003.

\bibitem{montanari2006rigorous}
Andrea Montanari and Guilhem Semerjian.
\newblock Rigorous inequalities between length and time scales in glassy
  systems.
\newblock {\em Journal of statistical physics}, 125(1):23--54, 2006.

\bibitem{mezard1986sk}
Marc M{\'e}zard, Giorgio Parisi, and MA~Virasoro.
\newblock Sk model: The replica solution without replicas.
\newblock {\em EPL (Europhysics Letters)}, 1(2):77, 1986.

\bibitem{0305-4470-15-6-035}
T~Plefka.
\newblock Convergence condition of the tap equation for the infinite-ranged
  ising spin glass model.
\newblock {\em Journal of Physics A: Mathematical and General}, 15(6):1971,
  1982.

\bibitem{rangan2010estimation}
Sundeep Rangan.
\newblock Estimation with random linear mixing, belief propagation and
  compressed sensing.
\newblock In {\em 44th Annual Conference on Information Sciences and Systems
  (CISS)}, pages 1--6. IEEE, 2010.

\bibitem{DonohoMaleki09}
David~L. Donoho, Arian Maleki, and Andrea Montanari.
\newblock Message-passing algorithms for compressed sensing.
\newblock {\em Proc. Natl. Acad. Sci.}, 106(45):18914--18919, 2009.

\bibitem{bayati2011dynamics}
Mohsen Bayati and Andrea Montanari.
\newblock The dynamics of message passing on dense graphs, with applications to
  compressed sensing.
\newblock {\em IEEE Transactions on Information Theory}, 57(2):764--785, 2011.

\bibitem{caltagirone2014convergence}
Francesco Caltagirone, Lenka Zdeborov{\'a}, and Florent Krzakala.
\newblock On convergence of approximate message passing.
\newblock In {\em IEEE International Symposium on Information Theory}, pages
  1812--1816. IEEE, 2014.

\bibitem{0305-4470-24-9-024}
A~Georges and J~S Yedidia.
\newblock How to expand around mean-field theory using high-temperature
  expansions.
\newblock {\em Journal of Physics A: Mathematical and General}, 24(9):2173,
  1991.

\bibitem{deshpande2015finding}
Yash Deshpande and Andrea Montanari.
\newblock Finding hidden cliques of size $\sqrt{N/e}$ in nearly linear time.
\newblock {\em Foundations of Computational Mathematics}, pages 1--60, 2015.

\bibitem{krzakala2013spectral}
Florent Krzakala, Cristopher Moore, Elchanan Mossel, Joe Neeman, Allan Sly,
  Lenka Zdeborov{\'a}, and Pan Zhang.
\newblock Spectral redemption in clustering sparse networks.
\newblock {\em Proceedings of the National Academy of Sciences},
  110(52):20935--20940, 2013.

\bibitem{hoyle2004principal}
David~C Hoyle and Magnus Rattray.
\newblock Principal-component-analysis eigenvalue spectra from data with
  symmetry-breaking structure.
\newblock {\em Physical Review E}, 69(2):026124, 2004.

\bibitem{richard2014statistical}
Emile Richard and Andrea Montanari.
\newblock A statistical model for tensor {PCA}.
\newblock In {\em Advances in Neural Information Processing Systems}, pages
  2897--2905, 2014.

\end{thebibliography}

\begin{appendices}

\section{Mean Field equations}
\label{Mean_Field_Appendix}
In order to compare the Low-RAMP algorithm with the commonly used
variational mean field inference we write here the variational mean
field equations in the same notation we used for Low-RAMP. We also
write the mean field free energy.

For the symmetric vector-spin glass the naive mean field equations
read 
\begin{eqnarray}
B_{X,i}^t &=& \sum\limits_{k = 1}^N \frac{1}{\sqrt{N}} S_{ki}
              {\hat{x}}_k^t\, ,
\\
A_{X,i}^t &=& \frac{1}{N}\sum\limits_{k = 1}^N (S_{ki}^2 - R_{ki})
              \left({{\hat{x}}_{k}^t{{\hat{x}}_{k}^{t,\top} +
              \sigma_{x,k}^t}} \right)  \, ,
\\
{\hat{x}}_i^{t+1} &=& f_{\rm in}^x(A_{X,i}^t,B_{X,i}^t)\, ,
\label{Mean_XX}
\\
\sigma_{x,i}^{t+1} &=& \frac{\partial f_{\rm in}^x }{\partial B}(A_{X,i}^t,B_{X,i}^t) \,.
\end{eqnarray}
The mean field free energy for the symmetric $X X^\top$ case reads 
\begin{multline}
F_{X X^\top}^{\rm MF}(\{A_{X,i}\},\{B_{X,i}\}) = \sum\limits_{1 \leq i \leq N} \log({\cal Z}_x(A_{X,i},B_{X,i})) - B_{X,i}^\top {\hat{x}}_i
+ \frac{1}{2}{\rm Tr}\left[A_{X,i}	({\hat{x}}_i {\hat{x}}_i^\top + \sigma_{x,i}) \right]
\\
+\frac{1}{2}\sum\limits_{1 \leq i ,j\leq N}\left[ \frac{1}{\sqrt{N}}S_{ij} {\hat{x}}_i^\top {\hat{x}}_j + \frac{(R_{ij} - S_{ij}^2)}{2N} {\rm Tr}\left[({\hat{x}}_i {\hat{x}}_i^\top + \sigma_{x,i}) ({\hat{x}}_j {\hat{x}}_j^\top + \sigma_{x,j}) \right]
\right] \,.
\end{multline}

For the bipartite $IUV^\top$ case the mean field equations read
\begin{eqnarray}
B^t_{U,i}&=& \frac{1}{\sqrt{N}}
 \sum\limits_{l=1}^M S_{il} \hat{v}^t_{l}\, ,
\\
A_{U}^t &=& \frac{1}{N}\sum\limits_{l=1}^M \left(S_{il}^2 - R_{il}
            \right) \left( \hat{v}_{l}^t{\hat{v}_{l}^{t,\top}} +
            \sigma_{u,l}^t \right)\, ,
\\
\hat{u}^t_{i} &=& f^u_{\rm in}(A^t_{U},B^t_{U,i})\, ,
\\
\sigma^t_{u,i} &=& \left(\frac{\partial f^u_{\rm in}}{ \partial B}
                   \right)(A^t_{U},B^t_{U,i})\, ,
\\
B^t_{V,j} &=& \frac{1}{\sqrt{N}} \sum\limits_{k=1}^N S_{kj}
              \hat{u}_k^t\, ,
\\
A_{V,j}^t &=& \frac{1}{N}\sum\limits_{k = 1}^N (S_{kj}^2 - R_{kj})
              \left(\hat{u}_{k}^t{\hat{u}_{k}^{t,\top}} +
              \sigma_{u,k}^t \right)\, ,
\\
\hat{v}^{t+1}_{j} &=& f^v_{\rm in}(A^t_{V,j},B^t_{V,j})\, ,
\\
\sigma^{t+1}_{v,j} &=& \left(\frac{\partial f^v_{\rm in}}{ \partial B} \right)(A^t_{V,j},B^t_{V,j}) \,.
\end{eqnarray}
The mean field free energy for the bipartite case reads
\begin{multline}
F_{UV^\top}^{\rm MF}(\{A_{U,i}\},\{B_{U,i}\},\{A_{V,j}\},\{B_{V,j}\}) = \sum\limits_{1 \leq i \leq N} \log({\cal Z}_u(A_{U,i},B_{U,i}))  -  {B_{U,i}}^\top {\hat{u}}_i
+ \frac{1}{2}{\rm Tr}\left[A_{U,i}({\hat{u}}_i {\hat{u}}_i^\top + \sigma_{u,i}) \right]
\\
+
\sum\limits_{1 \leq j \leq M} \log({\cal Z}_v(A_{V,j},B_{V,j})) -  {B_{V,j}}^\top {\hat{v}}_j
+ \frac{1}{2}{\rm Tr}\left[A_{V,j}({\hat{v}}_j {\hat{v}}_j^\top + \sigma_{v,j}) \right]
\\
+\sum\limits_{1 \leq i \leq N,1 \leq j\leq M} \left[ \frac{1}{\sqrt{N}}S_{ij} \hat{u}	_i^\top \hat{v}_j + \frac{1}{2N}(R_{ij} - S_{ij}^2) {\rm Tr}\left[(\hat{u}_i \hat{u}_i^\top + \sigma_{u,i})(\hat{v}_j \hat{v}_j^\top + \sigma_{v,j})\right]\right]\,. 
\end{multline}
The difference between the mean field equations and the Low-RAMP
equations from section \ref{TAP_Equation_Section_XX} and \ref{TAP_Equation_Section_UV}  can be seen for both variables $A$ and $B$.

\section{Bethe free energy derived by the Plefka expansion}
\label{app:Plefka}

We present here a way to compute the Bethe free energy for relatively
generic class of Hamiltonians. The main hypothesis for this computation to be  exact in
the large $N$ limit is that the Fisher score matrix $S_{ij}$ can be
well approximated as having elements with weak correlation of order $\frac{1}{N}$ one from another. If this is not true then the following equations do not apply.

As it turns out the mean-field free energy from Appendix \ref{Mean_Field_Appendix} is
just the first order term in the high temperature expansion of the free-energy.
as was derived in \cite{0305-4470-15-6-035,0305-4470-24-9-024} first
by Plefka and then by Georges and Yedidia. The Bethe free energy is
then the 2nd order expansion. 
Suppose one is given a general system with with $N$ classical
variables $x_i$ that weakly interact. The Hamiltonian $H$ defined the
free energy
\begin{equation}
\Phi = \log \left( {\rm Tr} \left[ \exp(\beta H(x_1,\cdots,x_N)) \right] \right) \,.
\end{equation}
We define a new Hamiltonian that fixes the marginal by
\begin{eqnarray}
\beta H_{\rm Field}( \{ \lambda_i \}) = \beta H + \sum\limits_{1 \leq i\leq N}\lambda_i(x_i,\beta)\,, \label{ModifiedHamiltonian}
\\
\Phi_{\rm Field}(\beta,\{ \lambda_i \} )  = \log \left({\rm Tr}\left[
  \exp(\beta H_{\rm Field}(\{ \lambda_i \}) ) \right] \right)\, .
\end{eqnarray}
By taking the Legendre transform one gets
\begin{eqnarray}
\Phi_{\rm Legendre}(\beta,\{ \mu_i(x_i) \} )  &=& \min \left\{
                                                  \Phi_{\rm
                                                  Field}(\beta,\{
                                                  \lambda_i \} ) -
                                                  \sum\limits_{ 1 \leq
                                                  i \leq n} \int {\rm
                                                  d}x \mu_i(x)
                                                  \lambda_i(x)
                                                  \right\}\, , 
\\
\{\Lambda_i(x_i) \}  &=& {\rm argmin} \left\{ \Phi_{\rm
                         Field}(\beta,\{ \lambda_i \} ) -
                         \sum\limits_{ 1 \leq i \leq n} \int {\rm d}x
                         \mu_i(x) \lambda_i(x) \right\}\, , 
\end{eqnarray}
here the minimization is done over the fields $\lambda_i(x_i)$.
The $\mu_i(x_i)$ are marginal density probabilities that one aims to
impose on the system. The $\Lambda_i(x_i)$ are the fields one uses to
fix the marginals equal to $\mu_i(x_i)$. The $\Lambda_i(x_i)$ depend
on the problem $H$, $\beta$ and the $\mu_i(x_i)$.
Because the $\Lambda_i(x_i)$ are defined up to a constant \begin{equation}
\int {\rm d}x \mu_i(x) \Lambda_i(x) = 0
\end{equation}
is imposed.
According to the definition of the Legendre transformation on has
\begin{equation}
\Phi = \max \left\{  \Phi_{\rm Legendre}(\beta,\{ \mu_i(x_i) \} )
\right\}\, .
\end{equation}

To compute $\Phi_{\rm Legendre}(\beta,\{ \mu_i(x_i)$ we resort to a high temperature expansion. This procedure 
is called the Plefka expansion \cite{0305-4470-15-6-035}, it relies on the following high temperature expansion (we stop at order 2)
\begin{equation}
\Phi_{\rm Bethe} = \Phi_{\rm Frustrated}(\beta=0) + \beta \left( \frac{\partial \Phi_{\rm Legendre}}{\partial \beta} \right)(\beta=0) + \frac{\beta^2}{2} \left( \frac{\partial^2 \Phi_{\rm Legendre}}{\partial \beta^2} \right)(\beta=0) \,. \label{Expansion_F_TAP_General_1}
\end{equation}
The expansion was explained in detail in \citep{0305-4470-24-9-024},
here we remind the main steps. 
Let us introduce the following operator
\begin{equation}
U = H - \langle H \rangle - \sum\limits_{1 \leq i\leq N}\frac{\partial \Lambda_i(x_i)}{\partial \beta}\,,
\end{equation}
where the average is taken with respect to probability distribution induced by \eqref{ModifiedHamiltonian}.
One can show that for all observables $O$ one has
\begin{equation}
\forall \beta, \frac{\partial \langle O \rangle}{\partial \beta} = \left\langle \frac{\partial O}{\partial \beta}\right\rangle - \langle OU \rangle \,. \label{O_Derive_beta}
\end{equation}
According to the definition of $F_{\rm Frustrated}$ one can prove that.
\begin{equation}
\forall \beta, \frac{\partial \Phi_{\rm Legendre}}{\partial \beta} = \langle H \rangle \,. \label{FirstOrderFreeEnergy}
\end{equation}
Using (\ref{FirstOrderFreeEnergy}) and (\ref{O_Derive_beta}) we get
\begin{equation}
\forall \beta, \frac{\partial^2 \Phi_{\rm Legendre}}{\partial \beta^2} = \frac{\langle H \rangle}{\partial \beta} = -\langle H U \rangle  \,. 
\end{equation}
Therefore
\begin{equation}
\Phi_{\rm Bethe} = \Phi_{\rm Legendre}(\beta=0) + \beta \langle H \rangle(\beta=0) - \frac{\beta^2}{2}\langle HU \rangle(\beta=0) \,.
\end{equation}
We still need to compute $\frac{ \partial\Lambda_i(x_i)}{\partial \beta}$ at $\beta = 0$. This can be done by computing the derivative of the marginals with respect to $\beta$ and noticing that they have to be zero.
\begin{equation}
\left( \frac{\partial \langle \delta(x_i - \hat{x}_i) \rangle}{\partial \beta}\right) = \left( \frac{\partial \mu_i(\hat{x}_i)}{\partial \beta}\right) = 0 = \langle U \delta(x_i - \hat{x}_i) \rangle =  \left\langle\left( H - <H>  - \sum\limits_{1 \leq i\leq N}\frac{\partial \Lambda_i(x_i)}{\partial \beta} + \int {\rm d}\hat{x}_i \mu_i(\hat{x}_i)\frac{\partial \Lambda_i(\hat{x}	_i)}{\partial \beta} \right)\delta(x_i - \hat{x}_i)  \right\rangle \,.
\end{equation}
From this one deduces
\begin{equation}
\left\langle \frac{\partial \Lambda_i(x_i)}{\partial \beta} \delta(x_i
  - \hat{x}_i)	+ C(i,\beta)\delta(x_i - \hat{x}_i)s \right\rangle =
\langle H \delta(x_i - \hat{x}_i) - \langle H \rangle \delta(x_i -
\hat{x}_i) \rangle \, ,\label{ComputeLambda_1}
\end{equation}
where $C({i,\beta})$ is
\begin{equation}
C(i,\beta) = \int {\rm d}\hat{x}_i \mu_i(\hat{x}_i)\frac{\partial \Lambda_i(\hat{x}_i)}{\partial \beta} \,.
\end{equation}
By definition
\begin{eqnarray}
\langle \delta(x_i - \hat{x}_i)\rangle &=& \mu_i(\hat{x}_i)\, ,
\\
\langle H \delta(x_i - \hat{x}_i)\rangle &=& \mu_i(x_i)\langle H \rangle_{x_i = \hat{x}_i} \,,
\end{eqnarray}
where $\langle H \rangle_{x_i = \hat{x}_i}$ is the average energy
conditioned on the fact that $x_i = \hat{x}_i$. Using (\ref{ComputeLambda_1}) we get
\begin{equation}
\mu_i(\hat{x_i})\left(\frac{\partial \Lambda_i(\hat{x}_i)}{\partial
    \beta} \right) = \mu_i(\hat{x_i})\langle H \rangle_{x_i =
  \hat{x}_i} - \mu_i(\hat{x_i})\langle H \rangle + \mu_i(\hat{x_i})
C(i,\beta)\, , \label{Derivative_Lambda}
\end{equation}
\begin{equation}
\left(\frac{\partial \Lambda_i(\hat{x}_i)}{\partial \beta} \right) = \langle H \rangle_{x_i = \hat{x}_i} - \langle H \rangle +  C(i,\beta) \label{Derivative_Lambda2} \,.
\end{equation}
We can see from (\ref{ModifiedHamiltonian}) that $\Lambda_i(\hat{x}_i,\beta)$ is defined up to a constant we fix that constant by having
\begin{equation}
\int {\rm d}\hat{x}_i \mu_i(\hat{x}_i) \Lambda_i(\hat{x}_i) = 0 \,.
\end{equation}
Therefore we get
\begin{equation}
\forall \beta, \left(\frac{\partial \Lambda_i(\hat{x}_i)}{\partial \beta} \right) = \langle H \rangle_{x_i = \hat{x}_i} - \langle H \rangle\ = \frac{1}{\mu_i(\hat{x}_i)}\langle H \delta(x_i - \hat{x}_i)\rangle - \langle H \rangle\,,\label{Derivative_Lambda3}
\end{equation}
where once again $\langle H \rangle_{x_i = \hat{x}_i}$ is the average of the energy where we have conditioned on the event $x_i = \hat{x_i}$.
Since we do all expansion around  $\beta = 0$ one is able to compute
all the means present in this formula since at $\beta = 0$ the density probability of the system at $\beta=0$ is just
\begin{equation}
P_{\rm Factorised} = \prod\limits_{i=1 \cdots N} \mu_i(x_i)\,.
\end{equation}

By using (\ref{Expansion_F_TAP_General_1}) and (\ref{Derivative_Lambda3}) around $\beta=0$ one gets
\begin{equation}
\Phi_{\rm Bethe} = \sum\limits_{1 \leq i \leq N} S_{\mu_i}
+ \beta \langle H\rangle  + \frac{\beta^2}{2}\left[\langle H^2\rangle - \langle H \rangle^2 - \sum\limits_{1 \leq i \leq N}\int {\rm d}x_i \mu_i(\hat{x}_i)(\langle H \rangle_{x_i = \hat{x}_i} - \langle H \rangle)^2\right] \,, \label{FreeEnergyGeneral}
\end{equation}
where here $S_{\mu_i}$ is here the entropy of the density probability
$\mu_i(x_i)$. In this setting all terms of order $3$ and beyond will
give a zero contribution in the large $N$ limit. This is due to the
fact that the matrix $Y$ has  components largely uncorrelated since
they were taken independently from one another (the signal induces
correlation that are too small to be significant). 
One reason why this approach might fail is that even though order higher than 3 might disappears in the large $N$ limit their sum might be non convergent. This is the sign of a replica symmetry breaking in the system.

From this expression we deduce the Bethe free energy for the $X   X^\top$ case as a function of the variables $A_{X,i}$ and $B_{X,i}$.
The marginals that will be enforced are
\begin{eqnarray}
P_{X,{\rm Marginal}}(x_i) = \frac{1}{{\cal Z}_x(A_{X,i},B_{X,i})}P_X(x_i) \exp\left( B_{X,i}^\top x_i -
  \frac{x_i^\top A_{X,i} x_i}{2}\right) \, , \label{Marginal_X_Min}
\\
P_{U,{\rm Marginal}}(u_i) = \frac{1}{{\cal Z}_u(A_{U,i},B_{U,i})}
  P_U(u_i) \exp\left( B_{U,i}^\top u_i - \frac{u_i^\top A_{U,i}
  u_i}{2}\right) \, ,\label{Marginal_U_Min}
\\
P_{V,{\rm Marginal}}(v_j) = \frac{1}{{\cal Z}_v(A_{V,j},B_{V,j})} P_X(v_j) \exp\left( B_{V,i}^\top v_j - \frac{v_j^\top A_{V,j} v_j}{2}\right) \label{Marginal_V_Min}
\,.
\end{eqnarray}
The Hamiltonians are
\begin{equation}
H_{X  X^\top} = \sum\limits_{i = 1 \cdots N} \log(P_X(x_i)) + \sum\limits_{1 \leq i < j \leq N} \left[ g\left(Y_{ij},\frac{x_i^\top   x_j}{\sqrt{N}} \right)- g(Y_{ij},0) \right]\,, \label{Hamiltonian_XX}
\end{equation}
and
\begin{equation}
H_{UV^\top} = \sum\limits_{i = 1 \cdots N} \log(P_U(u_i)) + \sum\limits_{j = 1 \cdots M} \log(P_V(v_j)) + \sum\limits_{1 \leq i \leq N,1 \leq j \leq M} \left[g\left(Y_{ij},\frac{u_i^\top v_j}{\sqrt{N}} \right)- g(Y_{ij},0) \right]		 \label{Hamiltonian_UV} \,.
\end{equation}
Using (\ref{Marginal_X_Min},\ref{Marginal_U_Min},\ref{Marginal_V_Min})
and (\ref{Hamiltonian_XX},\ref{Hamiltonian_UV}) in
(\ref{FreeEnergyGeneral}) one gets for the symmetric $XX^\top$ case
\begin{multline}
\Phi_{XX^\top}(\{A_{X,i}\},\{B_{X,i}\}) = \sum\limits_{1 \leq i \leq N} \left[\log({\cal Z}_x(A_{X,i},B_{X,i})) - B_{X,i}^\top {\hat{x}}_i
+ \frac{1}{2}{\rm Tr}\left[A_{X,i}	({\hat{x}}_i {\hat{x}}_i^\top + \sigma_{x,i}) \right] \right]
\\
+\frac{1}{2}\sum\limits_{1 \leq i ,j\leq N}\left[ \frac{1}{\sqrt{N}}S_{ij} {\hat{x}}_i^\top {\hat{x}}_j + \frac{R_{ij}}{2N} {\rm Tr}\left[({\hat{x}}_i {\hat{x}}_i^\top + \sigma_{x,i}) ({\hat{x}}_j {\hat{x}}_j^\top + \sigma_{x,j}) \right]-\frac{S_{ij}^2}{2N}  {\rm Tr} \left[ {\hat{x}}_i {\hat{x}}_i^\top   {\hat{x}}_j {\hat{x}}_j^\top  \right]
-\frac{1}{N}S_{ij}^2 {\rm Tr} \left[\sigma_{x,i}   \sigma_{x,j}
\right] \right]\, .
\label{TAP_XX_FreeEnergy_app}
\end{multline}

For the $U V^\top$ case
\begin{multline}
\Phi_{{\rm Bethe},UV^\top}(\{A_{U,i}\},\{B_{U,i}\},\{A_{V,j}\},\{B_{V,j}\}) = \sum\limits_{1 \leq i \leq N} \left[\log({\cal Z}_u(A_{U,i},B_{U,i}))  -  {B_{U,i}}^\top {\hat{u}}_i
+ \frac{1}{2}{\rm Tr}\left[A_{U,i}({\hat{u}}_i {\hat{u}}_i^\top + \sigma_{u,i}) \right] \right]
\\
+
\sum\limits_{1 \leq j \leq M} \left[\log({\cal Z}_v(A_{V,j},B_{V,j})) -  {B_{V,j}}^\top {\hat{v}}_j
+ \frac{1}{2}{\rm Tr}\left[A_{V,j}({\hat{v}}_j {\hat{v}}_j^\top + \sigma_{v,j}) \right] \right]
\\
+\sum\limits_{1 \leq i \leq N,1 \leq j\leq M}\left[ \frac{1}{\sqrt{N}}S_{ij} \hat{u}	_i^\top \hat{v}_j + \frac{1}{2N}R_{ij} {\rm Tr}\left[(\hat{u}_i \hat{u}_i^\top + \sigma_{u,i})(\hat{v}_j \hat{v}_j^\top + \sigma_{v,j})\right]- \frac{S_{ij}^2{\rm Tr}\left( \hat{u}_i \hat{u}_i^\top \hat{v}_j \hat{v}_j^\top\right)}{2N}
-\frac{1}{N}S_{ij}^2 {\rm Tr} \left[\sigma_{u,i}  \sigma_{v,j} \right]\right] \,,
\end{multline}
where the $\hat{x}_i,\hat{u}_i,\hat{v}_j,\sigma_{x,i},\sigma_{u,i},\sigma_{v,j}$ are the mean of the marginals we want to enforce.
These are all function of variables $B_i$, and $A_i$.  One should try
and minimize these function. This is what Low-RAMP does. Fixed points
of the Low-RAMP algorithm are stationary points of the of the Bethe free energy.
%
\section{Replica computation $UV^\top$ case.}
\label{appendix:ReplicaComputation}
In this appendix we present the derivation replica free-energy in the case of the $UV^\top$ case.
In the coming computation the indices $i$ and $k$ will go from $1$ to $N$. And $j$ and $l$ will go from $1$ to $M$.
\be
{\cal Z}(\{ Y_{ij} \}) = \int  \prod_{i} {\rm d} u_i P_U({u}_i)
\prod_{j} {\rm d} v_j P_U({v}_j) \prod_{ij} \exp \left(g\left(
    Y_{ij},\frac{u_i^\top v_j}{\sqrt{N}} \right) - g\left(
    Y_{ij},0\right) \right)\, .
\ee
One would like to compute the average of $ \langle \log({\cal Z}(\{ Y_{ij} \})) \rangle $.
This can be computed using the replica trick
\begin{equation}
\langle \log({\cal Z}(\{ Y_{ij} \})) \rangle_Y = \lim\limits_{n
  \rightarrow 0} \frac{ \langle{\cal Z}^n \rangle - 1 }{n} =
\lim\limits_{n \rightarrow 0} \frac{ \log(\langle{\cal Z}^n
  \rangle)}{n}\, .
\end{equation}
These can be hopefully be computed for any $n \in \mathbb{N}$. We will
then compute this function as $n \rightarrow 0$. We start with evaluating
\begin{equation}
{\cal Z}^n(\{ Y_{ij} \} ) = \int \prod\limits_{a = 1 \cdots n}\prod\limits_{i=1 \cdots N} {\rm d}u^a_i P^a_U(u^a_i)
	     \prod\limits_{j=1 \cdots M} {\rm d}v^a_j P_V^a(v^a_j)
             \exp\prod\limits_{i = 1 \cdots N,j=1\cdots M}
             \exp\left(g\left(Y_{ij},\frac{u_i v_j}{\sqrt{N}} \right)
               - g\left(Y_{ij},0\right)\right)\, .
\end{equation}
Therefore one has
\begin{multline}
\mathbb{E}({\cal Z}^n) = \int \prod\limits_{i = 1 \cdots N,j = 1 \cdots M}  {\rm d}Y_{ij} P_{\rm out}\left(Y_{ij}, w = 0 \right)
\prod\limits_{a = 0 \cdots n} \left( \left[ \prod \limits_{i=1 \cdots N} 	{\rm d}u^a_i P_U^a(u^a_i)\prod\limits_{j=1 \cdots M} {\rm d}v^a_j P_V^a(v^a_j) \right] \right.\\ \left. \left[ \prod\limits_{i = 1 \cdots N,j=1\cdots M} \exp\left(\sum\limits_{a = 0 \cdots n}g^a\left(Y_{ij},\frac{u_i^a v_j^a}{\sqrt{N}} \right) - g^a\left(Y_{ij},0\right)\right) \right] \right)\,,
\end{multline}
where
\begin{itemize}
\item if $a=0$ then $g^a = g^0 = \log(P_{\rm out}(Y,w))$, $P^a_U(u) = P_{U_0}(u)$ and $P^a_V(v) = P_{V_0}(v)$
\item if $a \neq 0$ then $g^a = g$, $P_U^a(u) = P_U(u)$ and $P_V^a(v) = P_V(v)$
\end{itemize}
We expand the function $g^a$ to order 2 and get
\begin{multline}
\mathbb{E}({\cal Z}^n) = \int \prod\limits_{i = 1 \cdots N,j = 1 \cdots M}  {\rm d}Y_{ij} P_{\rm out}\left(Y_{ij}, w = 0 \right)
 \prod\limits_{a = 0 \cdots n}\left(\left[ \prod \limits_{i=1 \cdots N} 	{\rm d}u^a_i P_U^a(u^a_i)\prod\limits_{j=1 \cdots M} {\rm d}v^a_j P_V^a(v^a_j) \right]  \right.\\ \left. \left[\prod\limits_{i = 1 \cdots N,j=1\cdots M} \exp\left(\sum\limits_{a = 0 \cdots n}	\left(\frac{\partial g^a}{\partial w}\right)_{Y_{ij},0}  \frac{u_i^a v_j^a}{\sqrt{N}} + \left(\frac{\partial^2 g^a}{\partial w^2}\right)_{Y_{ij},0}  \frac{(u_i^a v_j^a)^2}{2N} + O\left(\frac{1}{N^1.5} \right)
\right) \right] \right)\,.
\end{multline}
By expanding the exponential to order two one gets
\begin{multline}
\mathbb{E}({\cal Z}^n) = \int \prod\limits_{i = 1 \cdots N,j = 1 \cdots M}  {\rm d}Y_{ij} P_{\rm out}\left(Y_{ij}, w = 0 \right)
 \prod\limits_{a = 0 \cdots n}
 \left[\prod \limits_{i=1 \cdots N} 	{\rm d}u^a_i P_U^a(u^a_i)\prod\limits_{j=1 \cdots M} {\rm d}v^a_j P_V^a(v^a_j) \right] 
\\ \prod\limits_{i = 1 \cdots N,j=1\cdots M} \left[ 1 + \sum\limits_{a = 0 \cdots n}	\left(\frac{\partial g^a}{\partial w}\right)_{Y_{ij},0}  \frac{u_i^a v_j^a}{\sqrt{N}}+ \right.
\\
\left.
\sum\limits_{a = 1 \cdots n}	\left(\frac{\partial g}{\partial w}\right)_{Y_{ij},0}\left(\frac{\partial g^0}{\partial w}\right)_{Y_{ij},0}  \frac{u_i^a v_j^a u_i^0 v_j^0}{N}+
\sum\limits_{1 \leq a < b \leq n}	\left(\frac{\partial g}{\partial w}\right)_{Y_{ij},0}\left(\frac{\partial g}{\partial w}\right)_{Y_{ij},0}  \frac{u_i^a v_j^a u_i^b v_j^b}{N}
+ \right.
\\
\left.
\sum\limits_{a = 0 \cdots n} \left[
\left(\frac{\partial^2 g^a}{\partial w^2}\right)_{Y_{ij},0} + \left(\frac{\partial^2 g^a}{\partial w^2}\right)_{Y_{ij},0}^2 \right] \frac{(u_i^a v_j^a)^2}{2N} + O\left(\frac{1}{N^{1.5}} \right)
\right] \,.
\end{multline}
By averaging with respect to the $Y_{ij}$ one gets
\begin{multline}
\mathbb{E}({\cal Z}^n) = \int
\prod\limits_{a = 0 \cdots n}\left[ \prod \limits_{i=1 \cdots N}
  {\rm d}u^a_i P_U^a(u^a_i)\prod\limits_{j=1 \cdots M} {\rm d}v^a_j
  P_V^a(v^a_j) \right]\\  \prod\limits_{i = 1 \cdots N,j=1\cdots M} \left[1 + 
\sum\limits_{a = 1 \cdots n} \frac{u_i^a v_j^a u_i^0 v_j^0}{N \widehat\Delta}+
\sum\limits_{1 \leq a < b \leq n} \frac{u_i^a v_j^a u_i^b v_j^b}{N \widetilde\Delta}
+ 
\sum\limits_{a = 1 \cdots n} \overline R \frac{(u_i^a v_j^a)^2}{2N} + O\left(\frac{1}{N^{1.5}} \right)
\right] \,.
\\ 
= \int 
\left[ \prod\limits_{a = 0 \cdots n}\prod \limits_{i=1 \cdots N} 	{\rm d}u^a_i P_U^a(u^a_i)\prod\limits_{j=1 \cdots M} {\rm d}v^a_j P_V^a(v^a_j) \right]\\  \exp\left(  \sum\limits_{i = 1 \cdots N,j=1\cdots M}
\sum\limits_{a = 1 \cdots n} \frac{u_i^a v_j^a u_i^0 v_j^0}{N \widehat\Delta}+
\sum\limits_{1 \leq a < b \leq n} \frac{u_i^a v_j^a u_i^b v_j^b}{N \widetilde\Delta}
+ 
\sum\limits_{a = 1 \cdots n} \overline R \frac{(u_i^a v_j^a)^2}{2N} + O\left(\frac{1}{N^{1.5}} \right)
\right) \,.
\end{multline}

We now introduce the order parameters
\begin{eqnarray}
q^u_{ab} = \frac{1}{N}\sum\limits_{i = 1 \cdots N} u_i^a u_i^b\, ,
\\
q^v_{ab} = \frac{1}{M}\sum\limits_{j = 1 \cdots M} v_j^a v_j^b\, .
\end{eqnarray}This leads to 
\begin{multline}
\mathbb{E}({\cal Z}^n) = NM \int \prod_{0 \leq a \leq b \leq n} {\rm d}q^u_{ab} {\rm d}q^v_{ab} \exp\left(\frac{N \alpha}{\widehat\Delta}\sum\limits_{a = 1 \cdots n}q^u_{a0} q^v_{a0} + \frac{N \alpha}{\widetilde\Delta}\sum\limits_{1 \leq  a < b \leq n}q^u_{ab} q^v_{ab}
+ \frac{N \alpha}{2} \sum\limits_{a = \cdots n} \overline R 
q^u_{aa}q^v_{aa}
\right) \hat{\mathit{I}}_u(\{q^u_{ab}\}) \hat{\mathit{I}}_v(\{q^v_{ab}\})
\, ,
\end{multline}
where
\begin{equation}
\hat{\mathit{I}}_u(\{q^u_{ab}\}) = \int
  \left(\prod\limits_{a = 1 \cdots n}\prod \limits_{i=1 \cdots N}
    {\rm d}u^a_i P_V^a(u^a_i)\right)\prod\limits_{0 \leq a \leq b \leq
    n}\delta\left(\sum\limits_{i = 1 \cdots N} u_i^a u_i^b - N
    q^u_{ab}\right) \, ,
\end{equation}
and
\begin{equation}
\hat{\mathit{I}}_v(\{q^v_{ab}\}) = \int \left(\prod\limits_{a = 1 \cdots n}\prod \limits_{j=1 \cdots M} 	{\rm d}v^a_j P_V^a(v^a_j)\right)\prod\limits_{0 \leq a \leq b \leq n}\delta\left(\sum\limits_{i = 1 \cdots N} v_j^a v_j^b - M q^v_{ab}\right)  \,.
\end{equation}
Here $\hat{\mathit{I}}_u(\{q^u_{ab}\})$ and $
\hat{\mathit{I}}_v(\{q^v_{ab}\})$ are the entropy costs one pays in
order for the order parameters to take one specific value.
We can treat this constraint by going to Fourier space and then
rotating the path of integration
\begin{eqnarray}
\mathit{I}_u(\{q^u_{ab}\},\{\hat q^u_{ab}\}) =    \int
  \prod\limits_{a = 0 \cdots n}P_U^a(u_a){\rm d}u^a \exp\left(
    \sum\limits_{0 \leq a \leq b \leq n} \hat{q}^u_{ab}({u^a}
    u^b - q^u_{ab})\right)\, .
\\
\mathit{I}_v(\{q^v_{ab}\},\{\hat q^v_{ab}\}) =    \int
  \prod\limits_{a = 0 \cdots n}P_V^a(v_a){\rm d}v^a \exp\left(
    \sum\limits_{0 \leq a \leq b \leq n} \hat{q}^v_{ab}({v^a}
    v^b - q^v_{ab})\right)\, .
\end{eqnarray}
Therefore
\begin{multline}
\mathbb{E}({\cal Z}^n) = NM \int \prod_{0 \leq a \leq b \leq n} {\rm d}q^u_{ab} {\rm d}\hat q^u_{ab} {\rm d}q^v_{ab} {\rm d}\hat q^v_{ab} \exp\left(\frac{N \alpha}{\widehat\Delta}\sum\limits_{a = 1 \cdots n}q^u_{a0} q^v_{a0} + \frac{N \alpha}{\widetilde\Delta}\sum\limits_{1 \leq  a < b \leq n}q^u_{ab} q^v_{ab}
+ \frac{N \alpha}{2} \sum\limits_{a = \cdots n} \overline R 
q^u_{aa}q^v_{aa}
\right)
\\
 + \mathit{I}_u(\{q^u_{ab}\},\{\hat q^u_{ab}\}))^N + \mathit{I}_v(\{q^v_{ab}\},\{\hat q^v_{ab}\}))^M
\end{multline}
We need to extremize this function with respect to all variables. By
taking the derivative equal to 0 with respect to variables $q^u_{ab}$
and $q^v_{ab}$ we get
\begin{eqnarray}
\hat{q}_{a0}^u &=& \frac{\alpha \hat{q}_{a0}^v}{\widehat\Delta}
                   ,\qquad \hat{q}_{a0}^v =
                   \frac{\hat{q}_{a0}^u}{\widehat\Delta}\, ,
\\
\forall 1 \leq a < b \leq n,\quad \quad \hat{q}_{ab}^u &=& \frac{\alpha
                                               \hat{q}_{ab}^v}{\widetilde\Delta}
                                               ,\qquad \hat{q}_{ab}^v
                                               =
                                               \frac{\hat{q}_{ab}^u}{\widetilde\Delta}\, ,
\\
\hat{q}_{aa}^u &=& \overline R\alpha \hat{q}_{ab}^v,\qquad
                   \hat{q}_{ab}^v = \overline R \hat{q}_{ab}^u\, .
\end{eqnarray}

We now assume the replica symmetric ansatz
\begin{eqnarray}
\forall a,b > 1, \quad q_{ab}^u &=& \delta^a_b (\Sigma_u + Q_u) + (1-
                                    \delta^a_b)Q_u\, ,
\\
\forall a > 1, \quad q_{a0}^u &=& M_u\, ,
\\
\forall a,b > 1, \quad q_{ab}^v &=& \delta^a_b (\Sigma_v + Q_v) + (1-
                                    \delta^a_b)Q_v\, ,
\\
\forall a > 1, \quad q_{a0}^v &=& M_v\, .
\end{eqnarray}
We can then express the free energy.
Let us compute $I_u$ and $I_v$ in that RS ansatz
\begin{equation}
 \mathit{I}_u(\{q^u_{ab}\},\{\hat q^u_{ab}\}) \exp \left(
  \sum\limits_{0 \leq a \leq b \leq n}\hat{q}^u_{ab} q^u_{ab} \right)
=  \int \prod\limits_{a = 0 \cdots n}P_U^a(u_a){\rm d}u^a \exp\left(
  \sum\limits_{0 \leq a \leq b \leq n} \hat{q}^u_{ab}{u^a}
  u^b\right) \, .
\end{equation}
By using Hubbard-Stratonovich identity one gets
\begin{eqnarray}
&\mathit{I}_u(\{q^u_{ab}\},\{\hat q^u_{ab}\}) \exp \left(
  \sum\limits_{0 \leq a \leq b \leq n}\hat{q}^u_{ab} q^u_{ab} \right)
= \\ &= 
\int {\cal D}W P_{U_0}(u_0){\rm d}u^0 \left[ \int {\rm
    d}uP_U(u)\exp\left(  \frac{\alpha M_v}{\widehat\Delta} u u_0 + W
    \sqrt{\frac{\alpha Q_v}{\widetilde\Delta}} u - \left(\frac{\alpha
        Q_v}{\widetilde\Delta} - \alpha \overline R(Q_v + \Sigma_v)
    \right) \frac	{u^2}{2}\right) \right]^n	\, .
\end{eqnarray}
We can now compute the limit as $n \rightarrow 0$.
\begin{multline}
\lim\limits_{n \rightarrow 0}\frac{ \log(\mathit{I}_u(\{q^u_{ab}\},\{\hat q^u_{ab}\}))}{n} = \frac{\alpha Q_v Q_u}{ 2\widetilde\Delta}- \frac{\alpha M_v M_u}{\widehat\Delta} - \alpha \overline R(Q_v + \Sigma_v)(Q_u + \Sigma_u) + \\
\mathbb{E}_{W,u_0} \left[\log\left[ \int {\rm d}u P_U(u)\exp\left(  \frac{\alpha M_v}{\widehat\Delta} u u_0 + W \sqrt{\frac{\alpha Q_v}{\widetilde\Delta}} u - \left(\frac{\alpha Q_v}{\widetilde\Delta} - \alpha \overline R(Q_v + \Sigma_v) \right) \frac	{u^2}{2}\right) \right)	\right]
\\
\lim\limits_{n \rightarrow 0}\frac{\mathit{I}_v(\{q^v_{ab}\},\{\hat q^v_{ab}\})}{N} = \frac{Q_v Q_u}{ 2\widetilde\Delta}- \frac{M_v M_u}{\widehat\Delta} - \overline R(Q_v + \Sigma_v)(Q_u + \Sigma_u) + \\
\mathbb{E}_{W,v_0} \left[\log\left[ \int {\rm d}v P_V(v)\exp\left(
      \frac{M_u}{\widehat\Delta} v v_0 + W
      \sqrt{\frac{Q_u}{\widetilde\Delta}} u -
      \left(\frac{Q_u}{\widetilde\Delta} - \overline R(Q_u + \Sigma_u)
      \right) \frac	{u^2}{2}\right) \right)	\right]\, .
\end{multline}
This finally gives us the replica free energy.
\begin{multline}
\Phi_{\rm RS}(M_u,Q_u,\Sigma_u,M_v,Q_v,\Sigma_v) = 
\frac{\alpha Q_v Q_u}{ 2\widetilde\Delta}- \frac{\alpha M_v M_u}{\widehat\Delta} - \alpha \overline R(Q_v + \Sigma_v)(Q_u + \Sigma_u)
\\
+ \mathbb{E}_{W,u_0} \left[ {\cal Z}_u\left( \frac{\alpha Q_v}{\widetilde\Delta} - \alpha \overline R(Q_v + \Sigma_v) , \frac{\alpha M_v}{\widehat\Delta} u_0 + W \sqrt{\frac{\alpha Q_v}{\widetilde\Delta}}\right)	\right]
\\
+ \alpha \mathbb{E}_{W,v_0} \left[{\cal
    Z}_v\left(\frac{Q_u}{\widetilde\Delta} - \overline R(Q_u +
    \Sigma_u) , \frac{M_u}{\widehat\Delta} v_0 + Ws
    \sqrt{\frac{Q_u}{\widetilde\Delta}}\right)	\right]\, .
\end{multline}
This can also be computed with vectorial notations in both the $X
X^\top$ and $UV^\top$ setting leading to
eqs.~(\ref{BetheFreeEnergy_UV_DE}) and (\ref{BetheFreeEnergy_XX_DE}).

\section{Small $\rho$ expansion}
\label{AppendixSmallRho}

In this appendix we give the small-$\rho$ limits of the state
evolution update functions for the Rademacher-Bernoulli, Gauss-Bernoulli, 2 balanced groups and Bernoulli models.
\paragraph*{\bf Rademacher-Bernoulli model.}
We want to compute $\forall \beta > 0$ the limit $f_{\rm Rademacher-Bernoulli}(-\beta \log(\rho) )/\rho$ \eqref{SE_Rademacher_Bernoulli} as $\rho \rightarrow 0$, 
\begin{eqnarray}
\frac{f_{\rm Rademacher-Bernoulli}^{\rm SE}(-\beta \log(\rho))}{\rho} &=& \mathbb{E}_{W}\left[ \frac{\rho \tanh\left(-\beta \log(\rho)	+ W \sqrt{-\beta \log(\rho)}  \right) }{(1-\rho) \frac{\exp(-\beta \log(\rho)/2)}{\cosh\left(-\beta \log(\rho)	+ W \sqrt{-\beta \log(\rho)} \right)}  + \rho} \right].
\end{eqnarray}
Taking the small $\rho$ limit here we get
\begin{equation}
\lim\limits_{\rho = 0} \frac{f_{\rm Rademacher-Bernoulli}^{\rm SE}(-\beta
  \log(\rho))}{\rho} = \lim\limits_{\rho = 0} \frac{1}{\rho^{\beta
    /2-1} + 1} = 1(\beta > 2) \, , 
\end{equation}
where we used the fact that the noise $W \sqrt{-\beta \log(\rho)}$ is negligible compared to $\log(\rho)$ when $\rho \rightarrow 0$.

\paragraph*{\bf Gauss-Bernoulli model.}
We want to compute $\forall \beta > 0$ the limit $f_{\rm Gauss-Bernoulli}(-\beta \log(\rho) )/\rho$ \eqref{SE_Gaussian_Bernoulli} as $\rho \rightarrow 0$ and $r=1$.

\begin{eqnarray}
f_{\rm Gauss-Bernoulli}^{\rm SE}(x)/\rho &= \mathbb{E}_{W,x_0} \left[ f_{\rm in}^{\rm Gauss-Bernoulli} \left(x,x x_0 + \sqrt{x}W \right) x_0\right],
\\
 & =  \mathbb{E}_{W} \left[ \int \frac{\exp\left(-x_0^2/2 \right)}{\sqrt{2 \pi}} x_0 f_{\rm in}^{\rm Gauss-Bernoulli} \left(x,x x_0 + \sqrt{x}W \right) \right],
\\
&{\rm By}  \, {\rm intregrating \, by \, part \, one \, gets}& \nonumber
\\
&= x \mathbb{E}_{W} \left[ \int \frac{\exp\left(-x_0^2/2 \right)}{\sqrt{2 \pi}} \frac{\partial f_{\rm in}^{\rm Gauss-Bernoulli} }{\partial B}\left(x,x x_0 + \sqrt{x}W \right) \right],
\end{eqnarray}
Here $x x_0 + \sqrt{x} W$ is a Gaussian random variable of mean 0 and variance $x + x^2$. Therefore one has
\begin{equation}
f_{\rm Gauss-Bernoulli}^{\rm SE}(x)/\rho = x \mathbb{E}_{W} \left[ \frac{\partial f_{\rm in}^{\rm Gauss-Bernoulli} }{\partial B}\left(x,\sqrt{x^2 + x}W \right) \right],
\end{equation}
By making another integration by part one gets
\begin{equation}
f_{\rm Gauss-Bernoulli}^{\rm SE}(x)/\rho = \frac{x}{\sqrt{x^2 + x}}
\mathbb{E}_{W} \left[ f_{\rm in}^{\rm Gauss-Bernoulli}\left(x,\sqrt{x^2 + x}W \right) W\right] = \frac{x}{1 + x} \mathbb{E}_{W} \left[ W^2 \hat{\rho}(x, W \sqrt{x^2 + xs}\right],
\end{equation}
where
\begin{equation}
{\hat \rho}(x,W^2 (x^2 + x)) = \frac{\rho}{(1-\rho)
  \exp\left(\frac{-x W^2}{2}\right)\sqrt{1+x}+ \rho} \,. \label{Effective_Rho}
\end{equation}
By writing $x = - \beta \log(\rho)$
Depending on the value of $W$, $\hat{\rho}$ will either go to zero or
one as $\rho \rightarrow 0$. This will appear in the limit of $\hat{\rho}$.
By taking $x = - \beta \log(\rho)$ one gets
\begin{equation}
\lim\limits_{\rho \rightarrow 0} \hat{\rho}(x, W^2 (x^2 + x))
=\lim\limits_{\rho \rightarrow 0}\frac{1}{(1-\rho) \rho^{\frac{\beta
      W^2 }{2}-1}\sqrt{1-\beta \log(\rho)}+ 1} = 1(\beta W^2 > 2). 
\end{equation}
From this we deduce that 
\begin{equation}
\lim\limits_{\rho \rightarrow 0} \frac{f_{\rm Gauss-Bernoulli}^{\rm SE}(-\beta
  \log(\rho)	)}{\rho} = \lim\limits_{\rho \rightarrow 0} \frac{- \beta 
  \log(\rho)}{1 - \beta \log(\rho)} \mathbb{E}_W\left[ W^2 1\left( |W| >
    \sqrt{\frac{2}{\beta}} \right)\right] = \frac{2 \exp(-1/\beta)}{\sqrt{\beta \pi}} +\erfc \((
\frac{1}{\sqrt{\beta}}\)) \, .
\end{equation}
Here we used again the fact that the noise $W \sqrt{-\beta \log(\rho)}$ is negligible compared to $\log(\rho)$ when $\rho \rightarrow 0$.

\paragraph*{\bf 2 balanced groups.}
We want to compute $\forall \beta > 0$ the limit $f_{\rm Balanced}(-\beta \rho(1-\rho)\log(\rho(1-\rho)) )$ \eqref{F_Rho_Cliques} as $\rho \rightarrow 0$.
\begin{eqnarray}
f_{\rm Balanced}^{\rm SE}(- \beta \rho(1-\rho) \log( \rho(1-\rho) )) &=& \mathbb{E}_W\left[
\frac{2\rho(1-\rho)\sinh\left(\frac{-\beta \log( \rho(1-\rho))}{2} + W \sqrt{ -\beta \log(\rho (1-\rho)) }\right)}{1+2\rho(1-\rho)\left( \cosh\left(\frac{-\beta \log( \rho(1-\rho))}{2} + W \sqrt{ -\beta \log(\rho (1-\rho)) } \right)-1 \right)} \right]\,, \nonumber
\\
\lim\limits_{\rho \rightarrow 0	} f_{\rm Balanced}^{\rm SE}(- \beta \rho(1-\rho) \log( \rho(1-\rho) )) &=& \lim\limits_{\rho \rightarrow 0} 
\frac{2\left[\rho(1-\rho)\right]^{1 - \beta/2}}{1+2\left[\rho(1-\rho)\right]^{1-\beta/2}} = 1(\beta > 2)\,. 
\end{eqnarray}
Here we used the fact that the noise $W \sqrt{-\beta \log(\rho(1-\rho))}$ is negligible compared to $\log(\rho(1-\rho))$ when $\rho \rightarrow 0$.

\paragraph*{\bf Spiked Bernoulli model}
The state evolution update function in this model is given by \eqref{SE_Bernoulli}.
Once again we set $x = - \beta \log(\rho)$ to get
\begin{eqnarray}
\frac{f_{\rm Bernoulli}^{\rm SE}(- \beta \log(\rho))}{\rho} &=& \mathbb{E}_W\left[
                                                    \frac{1}{1 +
                                                    (1-\rho)
                                                    \exp\left(
                                                    (\beta/2 -
                                                    1)\log(\rho) +
                                                    W\sqrt{- \beta
                                                    \log(\rho)}
                                                    \right)} \right]\, ,
\\
\lim\limits_{\rho  \rightarrow 0} \frac{f_{\rm Bernoulli}^{\rm SE}(- \beta
  \log(\rho))}{\rho} &=& \lim\limits_{\rho  \rightarrow 0} \frac{1}{1
                         + (1-\rho)\rho^{\beta/2-1}} = 1(\beta > 2)\, .
\end{eqnarray}
Here we used the fact that the noise $W \sqrt{-\beta \log(\rho)}$ is negligible compared to $\log(\rho)$ when $\rho \rightarrow 0$.


To compute the limiting behaviour of the phase transitions we analyzed
the functions $f_{\rm Rademacher-Bernoulli}$,$f_{\rm Gauss-Bernoulli}$
or $f_{\rm Bernoulli}$ (that we will call $f_{\rho}$) as follows. 
We remind 
\begin{equation}
\Delta_{\rm Dyn}(\rho) = \max_{x \in \mathbb{R}^+}
\frac{f_{\rho}(x)}{x} = \frac{\rho}{-\log(\rho)} \max_{\beta \in
  \mathbb{R}^+} \frac{f_{\rho}(-\beta \log(\rho))}{\rho \beta}\, .
\end{equation}
In the small $\rho$ limit
\begin{equation}
\Delta_{\rm Dyn}(\rho) \sim_{\rho \rightarrow 0}
\frac{\rho}{-\log(\rho)} \max_{\beta \in \mathbb{R}^+} \frac{1(\beta >
  2)}{\beta} = \frac{\rho}{-2 \log(\rho)} \, .
\end{equation}
The information-theoretic phase transition in the small $\rho$ limit
is computed as follows
\begin{eqnarray}
\Delta_{\rm IT}(\rho) &=& \max\limits_{x \in \mathbb{R}^+}\left\{
  \Delta(x) ,\int\limits_{0}^{x}{\rm d}m   f_{\rho}\left( u\right)
  =\frac{xf_{\rho(x)}}{2}  \right\} \, ,
\\
\Delta_{\rm IT}(\rho) &\sim_{\rho \rightarrow 0}&
                            \max\limits_{\beta \in
                            \mathbb{R}^+}\left\{ \Delta(-\beta
                            \log(\rho)) ,\int\limits_{0}^{	-\beta
                            \log(\rho)}{\rm d}u   f_{\rho}\left(
                            u\right) =\frac{-\beta \log(\rho)
                            f_{\rho}(-	\beta \log(\rho))}{2}
                            \right\}\, , 
\\
\Delta_{\rm IT}(\rho) &\sim_{\rho \rightarrow 0}&
                            \frac{\rho}{-\log(\rho)}\max\limits_{\beta
                            \in \mathbb{R}^+}\left\{ \frac{1(\beta >
                            2)}{\beta} ,\int\limits_{0}^{\beta}{\rm
                            d}u   1(u  > 2) =\frac{\beta 1(\beta  >
                            2)}{2}  \right\}\, ,
\\
\Delta_{\rm IT}(\rho) &\sim_{\rho \rightarrow 0}& \frac{\rho}{-4
                            \log(\rho)}\, .
\end{eqnarray}

\section{Large rank behavior for the symmetric community detection}
\label{appendix:AppendixSpinodaleNetwork}

To locate the phase transitions $\Delta_{\rm IT}$ and $\Delta_{\rm
  Dyn}$ in the symmetric community detection case we make a couple of remarks about
the state evolution. 
First we remark that for $\forall x \in \mathbb{R}^+$
$b = {\cal M}_r(x)$ is a fixed point of (\ref{SE_b_Community}) for
$\Delta = {{\cal M}_r(x)}/{x}$. By definition $\Delta_{\rm Dyn}$ is
the greatest $\Delta$ for which a fixed point  exists
\begin{equation}
\Delta_{\rm Dyn}(r) = \max_{x  \in
  \mathbb{R}^+}\left\{\frac{{\cal M}_r(x)}{x} \right\}\, .\label{FormulaSpinodal}
\end{equation}
To find the $\Delta_{\rm IT}$ we notice that by taking the derivative
with respect to $Q$ and $M$ of (\ref{BetheFreeEnergy_XX_DE}) one finds a combination of (\ref{SE_Equation_XX_M}) and (\ref{SE_Equation_XX_Q}).
Therefore we have
\begin{equation}
\phi(b_2,\Delta) - \phi(b_1,\Delta) =  \frac{r-1}{2 r^2
  \Delta}\int\limits_{b_1}^{b_2} {\rm d}u {\cal
  M}_r\left(\frac{u}{\Delta}\right)-u \, .\label{ProtocolTransition}
\end{equation}
We deduce a way to compute $\Delta_{\rm IT}$ as
{\small
\begin{equation}
\Delta_{\rm IT}(r) = \max_{x  \in
  \mathbb{R}^+}\left\{\frac{{\cal M}_r(x)}{x}, {\rm s.t.}\,
  \int\limits_{0}^{x} {\rm d}u {\cal M}_r(u)=\frac{x{\cal M}_r(x)}{2}
\right\}\, . \label{FormulaStatic}
\end{equation}}
To compute $\Delta_{\rm IT}$ we evaluate ${\cal M}_r$ on a whole interval, then for each $x$ draw a
line between point $(0,0)$ and $(x,H(x))$. We then compute the area
between ${\cal M}_r$ and this line. When this area is zero then ${\cal M}_r(x)/{x}$ gives us $\Delta_{\rm IT}$.

In order to compute ${\cal M}_r$ we study the function ${\cal M}_r(x)$
where we take $x = \beta r \log(r)$, with $\beta = \Omega(1)$.
The important part of ${\cal M}_r$ is the integral
\begin{equation}
\int
\frac{\exp \left( \frac{x}{r} + \sqrt{\frac{x}{r}} u_1 \right)}
{
\exp \left( \frac{x}{r} + \sqrt{\frac{x}{r}}u_1 \right) + \sum\limits_{i = 2}^r \exp \left( \sqrt{\frac{x}{r}} u_i \right)
}
\prod\limits_{i=1}^r  {\cal D} u_i\, .
\end{equation}
The important variables to look at are (when taking $x = \beta r \log(r)$)
\begin{eqnarray}
F_1 &=& \exp \left( \frac{x}{r} + \sqrt{\frac{x}{r}}u_1 \right) = \exp \left( \beta \log(r) + \sqrt{\beta \log(r)}u_1 \right)\, , \label{M_r1}
\\
F_2 &=& \sum\limits_{i = 2}^r \exp \left( \sqrt{\frac{x}{r}}u_i  \right) = \sum\limits_{i = 2}^r \exp \left( \sqrt{\beta \log(r)}u_i  \right). \label{M_r2}
\end{eqnarray}
If the typical value of $F_1$ dominates $F_2$ as $r \rightarrow +\infty$ then ${\cal M}_r = 1$, otherwise if
$F_2$ dominates $F_1$ then ${\cal M}_r = 0$.


To estimate $F_1$, $F_2$ let us notice that with high probability the maximum value of the $u_i$ will be of order
$\sqrt{2 \log(r)}$. This is a general property of Gaussian variables that the maximum of~$r$ independent Gaussian variables of variance $\sigma^2$ is of order $\sigma \sqrt{2 \log(r)}$.
We can therefore compute the mean of $F_2$ while conditioning on
the fact that all of the $u_i$ are smaller than $\sqrt{2 \log(r)}$. This allows us to compute the typical value of $F_1$ as $
F_1 \sim r^\beta $. For $F_2$ we obtain:  
when $\beta < 2$ then
$F_2 \sim r^{ \frac{\beta}{2} + 1}$, and when 
 if $\beta > 2$ then $F_2 \sim r^{\sqrt{2 \beta}}$.
We have to look at which of the $F_1$ or $F_2$ has a higher exponent. $\beta = 2$ is the value at which these two exponent cross.
Therefore we have
\begin{equation}
{\cal M}_r(\beta r \log(r)) = {1}\left(\beta > 2 \right)\, .
\end{equation}
Now let us remind (\ref{FormulaSpinodal}) while keeping $x = \beta r
\log(r)$ to get
\begin{equation}
\Delta_{\rm Dyn}(r) \sim_{r \to \infty} \max \left\{\frac{{1}\left(\beta > 2 \right)}{ \beta r \log(r)}, \beta \in \mathbb{R}^+\right\} = \frac{1}{2 r \log(r)} \,.
\end{equation}
To get the information-theoretic transition we use (\ref{FormulaStatic}).
Let us find the $\beta$ that satisfies equation (\ref{FormulaStatic}) we get
\begin{equation}
\beta r \log(r) \left[1 - \frac{2}{\beta}\right ] = \frac{\beta r
  \log(r)}{2}\, .
\end{equation}
We deduce $\beta = 4$ and therefore 
\begin{equation}
\Delta_{\rm IT}(r) \sim_{r \to \infty} \frac{1}{4 r \log(r)}\, .
\end{equation}

%

\end{appendices}
\end{document}